\documentclass{eptcs-modified}
 \providecommand{\event}{} 
 \usepackage{breakurl}             
\usepackage{amsmath}
\usepackage{amssymb}
\usepackage{amsthm}
\usepackage{amscd}
\usepackage{graphics}
\usepackage{graphicx,color}

\long\def\greybox#1{%
    \newbox\contentbox%
    \newbox\bkgdbox%
    \setbox\contentbox\hbox to \hsize{%
        \vtop{
            \kern\columnsep
            \hbox to \hsize{%
                \kern\columnsep%
                \advance\hsize by -2\columnsep%
                \setlength{\textwidth}{\hsize}%
                \vbox{
                    \parskip=\baselineskip
                    \parindent=0bp
                    #1
                }%
                \kern\columnsep%
            }%
            \kern\columnsep%
        }%
    }%
    \setbox\bkgdbox\vbox{
        \pdfliteral{0.85 0.85 0.85 rg}
        \hrule width  \wd\contentbox %
               height \ht\contentbox %
               depth  \dp\contentbox
        \pdfliteral{0 0 0 rg}
    }%
    \wd\bkgdbox=0bp%
    \vbox{\hbox to \hsize{\box\bkgdbox\box\contentbox}}%
    \vskip\baselineskip%
}

\title{Point degree spectra of represented spaces}
\author{
Takayuki Kihara
\institute{Department of Mathematics\\ University of California, Berkeley, United States\footnote{Kihara has since moved to Nagoya University, Japan.}}
\email{kihara@i.nagoya-u.ac.jp}
\and
Arno Pauly
\institute{Clare College\\ University of Cambridge, United Kingdom\footnote{Pauly has since moved to the Universit\'e Libre de Bruxelles, Belgium.}}
\email{Arno.M.Pauly@gmail.com}
}
\def\titlerunning{Point degree spectra}
\def\authorrunning{T. Kihara \& A. Pauly}

\begin{document}
\theoremstyle{definition}
\newtheorem{theorem}{Theorem}[section]
\newtheorem{definition}[theorem]{Definition}
\newtheorem{problem}[theorem]{Problem}
\newtheorem{assumption}[theorem]{Assumption}
\newtheorem{corollary}[theorem]{Corollary}
\newtheorem{proposition}[theorem]{Proposition}
\newtheorem{lemma}[theorem]{Lemma}
\newtheorem{observation}[theorem]{Observation}
\newtheorem{fact}[theorem]{Fact}
\newtheorem{question}[theorem]{Open Question}
\newtheorem{conjecture}[theorem]{Conjecture}
\newtheorem{example}[theorem]{Example}
\newtheorem*{remark}{Remark}
\newtheorem*{claim}{Claim}
\newcommand{\dom}{\operatorname{dom}}
\newcommand{\id}{\textnormal{id}}
\newcommand{\Cantor}{{\{0, 1\}^\mathbb{N}}}
\newcommand{\Baire}{{\mathbb{N}^\mathbb{N}}}
\newcommand{\Lev}{\textnormal{Lev}}
\newcommand{\hide}[1]{}
\newcommand{\mto}{\rightrightarrows}
\newcommand{\uint}{{[0, 1]}}
\newcommand{\bft}{\mathrm{BFT}}
\newcommand{\lbft}{\textnormal{Linear-}\mathrm{BFT}}
\newcommand{\pbft}{\textnormal{Poly-}\mathrm{BFT}}
\newcommand{\sbft}{\textnormal{Smooth-}\mathrm{BFT}}
\newcommand{\ivt}{\mathrm{IVT}}
\newcommand{\cc}{\textrm{CC}}
\newcommand{\lpo}{\textrm{LPO}}
\newcommand{\llpo}{\textrm{LLPO}}
\newcommand{\aou}{AoU}
\newcommand{\Ctwo}{C_{\{0, 1\}}}
\newcommand{\C}{\textrm{C}}
\newcommand{\ic}[1]{\textrm{C}_{\sharp #1}}
\newcommand{\xc}[1]{\textrm{XC}_{#1}}
\newcommand{\me}{\name{P}.~}
\newcommand{\etal}{et al.~}
\newcommand{\eval}{\operatorname{eval}}
\newcommand{\Sierp}{Sierpi\'nski }
\newcommand{\isempty}{\operatorname{IsEmpty}}
\newcommand{\spec}{\textrm{Spec}}

\newcommand{\name}[1]{\textsc{#1}}
\newcommand{\repsp}[1]{\mathbf{#1}}
\newcommand{\repspb}[1]{\mathbf{#1}}

\newcounter{saveenumi}
\newcommand{\seti}{\setcounter{saveenumi}{\value{enumi}}}
\newcommand{\conti}{\setcounter{enumi}{\value{saveenumi}}}

\maketitle

\begin{abstract}
We introduce the point degree spectrum of a represented space as a substructure of the Medvedev degrees, which integrates the notion of Turing degrees, enumeration degrees, continuous degrees, and so on.
The notion of point degree spectrum creates a connection among various areas of mathematics including computability theory, descriptive set theory, infinite dimensional topology and Banach space theory.
Through this new connection, for instance, we construct a family of continuum many infinite dimensional Cantor manifolds with property $C$ whose Borel structures at an arbitrary finite rank are mutually non-isomorphic.
This provides new examples of Banach algebras of real valued Baire class two functions on metrizable compacta, and strengthen various theorems in infinite dimensional topology such as Pol's solution to Alexandrov's old problem.
\end{abstract}



\section{Introduction}
\label{sec:introduction}

\subsection*{Computability Theory}

In computable analysis \cite{pourel,weihrauchd}, there has for a long time been an interest in how complicated the set of codes of some element in a suitable spaces may be.
\name{Pour-El} and \name{Richards} \cite{pourel} observed that any real number, and more generally, any point in a Euclidean space, has a Turing degree.
They subsequently raised the question whether the same holds true for any computable metric space.
\name{Miller} \cite{miller2} later proved that various infinite dimensional metric spaces such as the Hilbert cube and the space of continuous functions on the unit interval contain points which lack Turing degrees, i.e.~have no simplest code w.r.t.~Turing reducibility.
A similar phenomenon was also observed in algorithmic randomness theory.
\name{Day} and \name{Miller} \cite{daymiller} showed that no neutral measure has Turing degree by understanding each measure as a point in the infinite dimensional space consisting of probability measures on an underlying space.

These previous works convince us of the need for a reasonable theory of degrees of unsolvability of points in an arbitrary represented space.
To establish such a theory, we associate a substructure of the Medvedev degrees with a represented space, which we call its \emph{point degree spectrum}.
A wide variety of classical degree structures are realized in this way, e.g., Turing degrees \cite{SoareBook}, enumeration degrees \cite{friedberg}, continuous degrees \cite{miller2}, degrees of continuous functionals \cite{Hinman73}.
What is more noteworthy is that the concept of a point degree spectrum is closely linked to infinite dimensional topology.
For instance, we shall see that for a Polish space all points have Turing degrees if and only if the small transfinite inductive dimension of the space exists.

In a broader context, there are various instances of smallness properties (i.e., $\sigma$-ideals) of spaces and sets that start making sense for points in an effective treatment; e.g., arithmetical (Cohen) genericity \cite{nies,OdiBook}, Martin-L\"of randomness \cite{nies}, and effective Hausdorff dimension \cite{lutz2}.
In all these cases, individual points can carry some amount of complexity -- e.g.~a Martin-L\"of random point is in some sense too complicated to be included in a computable $G_\delta$ set having effectively measure zero.
A recent important example \cite{PoZa12,Zap14} from forcing theory is genericity with respect to the $\sigma$-ideal generated by finite-dimensional compact metrizable spaces.
Our work provides an effective notion corresponding to topological invariants such as small inductive dimension or metrizability, and e.g.~allows us to say that certain points are too complicated to be (computably) a member of a (finite-dimensional) Polish space.

Additionally, the actual importance of point degree spectrum is not merely conceptual, but also applicative.
Indeed, unexpectedly, our notion of point degree spectrum turned out to be a powerful tool in descriptive set theory and infinite dimensional topology, in particular, in the study of restricted Borel isomorphism problems, as explained in more depth below.

\subsection*{Descriptive Set Theory}

A {\em Borel isomorphism problem} (see \cite{CenMau84,Maul76,Hrba78,Har78}) asks to find a nontrivial isomorphism type in a certain class of Borel spaces (i.e., topological spaces together with their Borel $\sigma$-algebras).
An {\em $\alpha$-th level Borel/Baire isomorphism} between $\repsp{X}$ and $\repsp{Y}$ is a bijection $f$ such that $E\subseteq \repsp{X}$ is of additive Borel/Baire class $\alpha$ if and only if $f[E]\subseteq \repsp{Y}$ is of additive Borel/Baire class $\alpha$.
These restricted Borel isomorphisms are introduced by \name{Jayne} \cite{Jayne74}, in Banach space theory, to obtain certain variants of the
\name{Banach-Stone Theorem} and the \name{Gelfand-Kolmogorov Theorem} for Banach algebras of the forms $\mathcal{B}_\alpha^*(\repsp{X})$ for realcompact spaces $\repsp{X}$.
Here, $\mathcal{B}_\alpha^*(\repsp{X})$ is the Banach algebra  of bounded real valued Baire class $\alpha$ functions on a space $\repsp{X}$ with respect to the supremum norm and the pointwise operation \cite{Bade73,Dash74,Jayne74}.
The first and second level Borel/Baire isomorphic classifications have been studied by several authors (see \cite{JayRog79a,JayRog79b}).
However, it is not certain even whether there is an uncountable Polish space whose $G_{\delta\sigma}$-structure is neither isomorphic to the real line nor to the Hilbert cube:

\begin{problem}[The Second-Level Borel Isomorphism Problem]\label{mainproblem1}
Are all uncountable Polish spaces second-level Borel isomorphic either to $\mathbb{R}$ or to $\mathbb{R}^\mathbb{N}$?
\end{problem}

\name{Jayne}'s result \cite{Jayne74} shows that this is equivalent to asking the following problem on Banach algebras.

\begin{problem}[see also \name{Motto Ros} \cite{MRos13}]\label{mainproblem2}
If $\repsp{X}$ is an uncountable  Polish space.
Then does there exist $n\in\mathbb{N}$ such that $\mathcal{B}_n^*(\repsp{X})$ is linearly isometric (or ring isomorphic) either to $\mathcal{B}_n^*([0,1])$ or to $\mathcal{B}_n^*([0,1]^\mathbb{N})$?
\end{problem}

The very recent successful attempts to generalize the Jayne-Rogers theorem and the Solecki dichotomy (see \cite{MRos13,PaSa12} and also \cite{KiNg} for a computability theoretic proof) revealed that two Polish spaces are second-level Borel isomorphic if and only if they are $\sigma$-homeomorphic.
Here, a topological space $\repsp{X}$ is $\sigma$-homeomorphic to $\repsp{Y}$ (written as $\repsp{X}\cong_\sigma^\mathfrak{T}\repsp{Y}$) if there are countable covers $\{\repsp{X}_i\}_{i\in\omega}$ and $\{\repsp{Y}_i\}_{i\in\omega}$ of $\repsp{X}$ and $\repsp{Y}$ such that $\repsp{X}_i$ is homeomorphic to $\repsp{Y}_i$ for every $i\in\omega$.
Therefore, the second-level Borel isomorphism problem can be reformulated as the following equivalent problem.

\begin{problem}[\name{Motto Ros} et al.~\cite{schlicht}]\label{prob:third}
Is any Polish space $\repsp{X}$ either $\sigma$-embedded into $\mathbb{R}$ or $\sigma$-homeomorphic to $\mathbb{R}^\mathbb{N}$?
\end{problem}

Unlike the classical Borel isomorphism problem, which was able to be reduced to the same problem on zero-dimensional Souslin spaces, the second-level Borel isomorphism problem is inescapably tied to {\em infinite dimensional topology} \cite{vMBook2}, since all transfinite dimensional uncountable Polish spaces are mutually second-level Borel isomorphic.

The study of $\sigma$-homeomorphic maps in topological dimension theory dates back to a classical work by \name{Hurewicz-Wallman} \cite{hurewicz} characterizing transfinite dimensionality.
\name{Alexandrov} \cite{Alek51} asked whether there exists a weakly infinite dimensional compactum which is not $\sigma$-homeomorphic to the real line.
\name{Roman Pol} \cite{pol2} solved this problem by constructing such a compactum.
Roman Pol's compactum is known to satisfy a slightly stronger covering property, called property $C$ \cite{AdGr78,Anc85,Hav74}.

Our notion of {\em degree spectrum} on Polish spaces serves as an invariant under second-level Borel isomorphism.
Indeed, an invariant which we call {\em degree co-spectrum}, a collection of Turing ideals realized as lower Turing cones of points of a Polish space, plays a key role in solving the second-level Borel isomorphism problem.
By utilizing these computability-theoretic concepts, we will construct a continuum many pairwise incomparable $\sigma$-homeomorphism types of compact metrizable $C$-spaces, that is:
\begin{itemize}
\item[] There is a collection $(\repsp{X}_\alpha)_{\alpha<2^{\aleph_0}}$ of continuum many compact metrizable $C$-spaces such that, whenever $\alpha\not=\beta$, $\repsp{X}_\alpha$ cannot be written as a countable union of homeomorphic copies of subspaces of $\repsp{X}_\beta$.
\end{itemize}
As mentioned above, this also shows that there are continuum many second-level Borel isomorphism types of compact metric spaces.
More generally, a {\em finite-level Borel embedding of $\repsp{X}$ into $\repsp{Y}$} is a finite-level Borel isomorphism between $\repsp{X}$ and a subset of $\repsp{Y}$ of finite Borel rank.
Then, our result entails the following as a corollary:
\begin{itemize}
\item[] There is a collection $(\repsp{X}_\alpha)_{\alpha<2^{\aleph_0}}$ of continuum many compact metrizable $C$-spaces such that, whenever $\alpha\not=\beta$, $\repsp{X}_\alpha$ cannot be finite-level Borel embedded into $\repsp{X}_\beta$.
\end{itemize}
The key idea is measuring the quantity of all possible Scott ideals realized within the degree co-spectrum of a given space.
Our spaces are completely described in the terminology of computability theory (based on \name{Miller}'s work on the continuous degrees \cite{miller2}).
Nevertheless, the first of our examples turns out to be second-level Borel isomorphic to \name{Roman Pol}'s compactum (but of course, other continuum many examples cannot be second-level Borel isomorphic to \name{Pol}'s compactum).
Hence, our solution can also be viewed as a refinement of \name{Roman Pol}'s solution to \name{Alexandrov}'s problem.

\subsection*{Summary of Results}

This work is part of a general development to study the descriptive theory of represented spaces \cite{pauly-overview-arxiv}, together with approaches such as synthetic descriptive set theory proposed in \cite{paulydebrecht,pauly-descriptive}.
In Section \ref{sec:pointdegreespectra}, we introduce the notion of point degree spectrum, and clarify the relationship with $\sigma$-continuity.
In Section \ref{sec:intermediate}, we introduce the notion of an $\omega$-left-CEA operator in the Hilbert cube as an infinite dimensional analogue of an $\omega$-CEA operator (in the sense of classical computability theory), and show that the graph of a universal $\omega$-left-CEA operator is an individual counterexample to Problems \ref{mainproblem1}, \ref{mainproblem2}, and \ref{prob:third}.
In Section \ref{sec:piecewisehomeo}, we describe a general procedure to construct uncountably many mutually different compacta under $\sigma$-homeomorphism.
In Section \ref{sec:intermediate-dimension}, we clarify the relationship between a universal $\omega$-left-CEA operator and \name{Roman Pol}'s compactum.
In Section \ref{sec:internalcharacterization}, we characterize represented spaces with effectively-fiber-compact representations (which are relevant for approaches to complexity theory along the lines of \name{Weihrauch} 's \cite{weihrauchf}) as precisely the computable metric spaces.
In Section \ref{sec:quasi-Polish}, we also look at the degree structures of nonmetrizable spaces and transfer techniques in the other direction: Employing geometric reasoning about $\mathbb{R}^n$, we prove results about semirecursive enumerations degrees.
The methods used in Sections \ref{sec:internalcharacterization}--\ref{sec:quasi-Polish} do not depend on those developed in Sections \ref{sec:intermediate}--\ref{sec:piecewisehomeo}.

\subsection*{Future work}
The methods introduced in this paper, in particular the notion of the point degree spectrum and the associated connection between topology and recursion theory, have already inspired and enabled several other studies. In \cite{GreKih}, \name{Gregoriades}, \name{Kihara} and \name{Ng} are attacking the generalized Jayne Rogers conjecture from descriptive set theory. A core aspect of this work is whether certain degree-theoretic results like the Shore-Slaman join theorem and Friedberg jump inversion theorem hold for the point degree spectra of Polish spaces.

Building upon Section \ref{sec:internalcharacterization}, \name{Andrews}, \name{Igusa}, \name{Miller} and \name{Soskova} \cite{miller4} used an effective metrization argument to show that the point degree spectrum of the Hilbert cube coincides with the almost-total enumeration degrees, which in turn is used to show the purely recursion-theoretic consequence that \emph{PA above} is definable in the enumeration degrees.

\name{Kihara}, \name{Lempp}, \name{Ng} and \name{Pauly} \cite{edegrees} have embarked on the systematic endeavour to classify the point degree spectra of second-countable spaces from \emph{Counterexample in Topology} \cite{CTopBook}. This has already proven to be a rich source for the fine-grained study of the enumeration degrees, as both previously studied substructures as well as new ones of interest to recursion theorists appear in this fashion.

Based on the results both in the present paper, and in the extension mentioned here, we are confident that both directions of the link between topology and recursion theory established here have significant potential for applications.

\section{Preliminaries}

\subsection{Computability Theory}

\subsubsection{Basic Notations}

We use the standard notations from modern computability theory and computable analysis.
We refer the reader to \cite{OdiBook,OdiBook1,SoareBook} for the basics on computability theory, and to \cite{pourel,weihrauchd,pauly-synthetic-arxiv} for the basics on computable analysis.

By $f:\subseteq X\to Y$, we mean a function from a subset of $X$ into $Y$.
Such a function is called a partial function.
We fix a pairing function $(m,n)\mapsto\langle m,n\rangle$, which is a computable bijection from $\mathbb{N}^2$ onto $\mathbb{N}$ such that $\langle m,n\rangle\mapsto m$ and $\langle m,n\rangle\mapsto n$ are also computable.
For $x,y\in\mathbb{N}^\mathbb{N}$, the join $x\oplus y\in\mathbb{N}^\mathbb{N}$ is defined by $(x\oplus y)(2n)=x(n)$ and $(x\oplus y)(2n+1)=y(n)$.
An oracle is an element of $\Cantor$ or $\mathbb{N}^\mathbb{N}$.
By the notation $\Phi_e^z$ we denote the computation of the $e$-th Turing machine with oracle $z$.
We often view $\Phi_e^z$ as a partial function on $\Cantor$ or $\mathbb{N}^\mathbb{N}$.
More precisely, $\Phi_e^z(x)=y$ if and only if given an input $n\in\mathbb{N}$ with oracle $x\oplus z$, the $e$-th Turing machine computation halts and outputs $y(n)$.
The terminology ``c.e.''\ stands for ``computably enumerable.''
For an oracle $z$, by ``$z$-computable'' and ``$z$-c.e.,'' we mean ``computable relative to $z$'' and ``c.e.\ relative to $z$.''
For an oracle $x$, we write $x'$ for the Turing jump of $x$, that is, the halting problem relative to $x$.
Generally, for a computable ordinal $\alpha$, we use $x^{(\alpha)}$ to denote the $\alpha$-th Turing jump of $x$.
Here, regarding the basics on computable ordinals and transfinite Turing jumps, see \cite{ChYuBook,SacksBook}.

We will repeatedly use the following fact, known as the {\em Kleene recursion theorem} or the {\em Kleene fixed point theorem}.
\begin{fact}[The Kleene Recursion Theorem; see {\cite[Theorem II.2.10]{OdiBook}}]\label{fact:kleene-recursion}
Given an oracle $z$ and a computable function $f:\mathbb{N}\to\mathbb{N}$, one can effectively find an index $e\in\mathbb{N}$ such that $\Phi_e^z$ and $\Phi^z_{f(e)}$ are the same partial function.
\end{fact}

\subsubsection{Represented spaces}
A \emph{represented space} is a pair $\repsp{X} = (X, \delta_X)$ of a set $X$ and a partial surjection $\delta_X : \subseteq \Baire \to X$.
Informally speaking, $\delta_X$ (called a representation) gives names of elements in $X$ by using infinite words.
It enables tracking of a function $f$ on abstract sets by a function on infinite words (called a realizer of $f$).
This is crucial for introducing the notion of computability on abstract sets because we already have the notion of computability on infinite words.

Formally, a function between represented spaces is a function between the underlying sets. For $f : \repsp{X} \to \repsp{Y}$ and $F : \subseteq \Baire \to \Baire$, we call $F$ a realizer of $f$, iff $\delta_Y(F(p)) = f(\delta_X(p))$ for all $p \in \dom(f\delta_X)$, i.e.~if the following diagram commutes:
 $$\begin{CD}
\Baire @>F>> \Baire\\
@VV\delta_\repsp{X}V @VV\delta_\repsp{Y}V\\
\repsp{X} @>f>> \repsp{Y}
\end{CD}$$
A map between represented spaces is called {\em computable} ({\em continuous}), iff it has a computable (continuous) realizer.
In other words, a function $f$ is computable (continuous) if there is a computable (continuous) function $F$ on infinite words such that, given a name $p$ of a point $x$, $F(p)$ returns a name of $f(x)$.
We also use the same notation $\Phi_e^z$ to denote a function on represented spaces realized by the $e$-th partial $z$-computable function.
 Similarly, we call a point $x \in \repsp{X}$ {\em computable}, iff there is some computable $p \in \Baire$ with $\delta_{X}(p) = x$, that is, $x$ has a computable name.
In this way, we think of a represented space as a kind of space equipped with the notion of computability.

If a set $X$ is already topologized, the above notion of continuity can be inconsistent with topological continuity.
To eliminate such an undesired situation, we shall consider a restricted class of representations which are consistent with a given topological structure, so-called {\em admissible representations}.
We will not go into the details of admissibility here, but just mention that if a $T_0$-space has a countable cs-network (a.k.a.\ a countable sequential pseudo-base), then it always has an admissible representation (see \name{Schr\"oder} \cite{schroder}).

A particularly relevant subclass of represented spaces are the computable Polish spaces, which are derived from complete computable metric spaces by forgetting the details of the metric, and just retaining the representation (or rather, the equivalence class of representations under computable translations). Forgetting the metric is relevant when it comes to compatibility with definitions in effective descriptive set theory as shown in \cite{pauly-gregoriades-arxiv}.

\begin{example}\label{example:representation}
The following are examples of admissible representations.
\begin{enumerate}
\item The representation of $\mathbb{N}$ is given by $\delta_\mathbb{N}(0^n10^\mathbb{N}) = n$.
 It is straightforward to verify that the computability notion for the represented space $\mathbb{N}$ coincides with classical computability over the natural numbers.
\item A {\em computable metric space} is a tuple $\repsp{M} = (M, d, (a_n)_{n \in \mathbb{N}})$ such that $(M,d)$ is a metric space and $(a_n)_{n \in\mathbb{N}}$ is a dense sequence in $(M,d)$ such that the relation
\[ \{(t,u,v,w) \: |\: \nu_{\mathbb{Q}}(t) < d(a_u, a_v) <\nu_{\mathbb{Q}}(w) \}\]
is recursively enumerable.
The {\em Cauchy representation} $ \delta_{\repsp{M}} \: : \subseteq \: \Baire \to M $ associated with the computable metric space $ \repsp{M} =  (M, d, (a_n)_{n \in \mathbb{N}}) $ is defined by
\[ \delta_{\repsp{M}}(p) = x \: : \: \Longleftrightarrow  \begin{cases}
      d(a_{p(i)}, a_{p(k)}) \leq 2^{-i} \text{ for } i < k\\
   \text{and } x = \lim\limits_{i\rightarrow \infty}a_{p(i)}
  \end{cases} \]
\item Another, more general subclass are the quasi-Polish spaces introduced by \name{de~Brecht} \cite{debrecht6}. A represented space $\repsp{X} = (X, \delta_\repsp{X})$ is \emph{quasi-Polish}, if it is countably based, admissible and $\delta_\repsp{X} : \Baire \to \repsp{X}$ is total. These include the computable Polish spaces as well as $\omega$-continuous domains.
\item Generally, a topological $T_0$-space $\repsp{X}$ with a countable base $\mathcal{B}=\langle{B_n\rangle}_{n\in\mathbb{N}}$ is naturally represented by defining $\delta_{(\repsp{X},\mathcal{B})}(p)=x$ iff $p$ enumerates the code of a neighborhood basis for $x$, that is, ${\rm range}(p)=\{n\in\mathbb{N}:x\in B_n\}$.
One can also use a network to give a representation of a space as suggested above.
\end{enumerate}
\end{example}

We always assume that $\Cantor$, $\mathbb{R}^n$, and $[0,1]^\mathbb{N}$ are admissibly represented by the Cauchy representations obtained from their standard metics.
A real $x\in\mathbb{R}$ is {\em left-c.e.}\ if there is a computable sequence $(q_n)_{n\in\mathbb{N}}$ of rationals such that $x=\sup_nq_n$.
Generally, a real $x\in\mathbb{R}$ is {\em left-c.e.\ relative to $y\in\repsp{X}$} if there is a partial computable function $f:\subseteq\repsp{X}\to\mathbb{Q}^\mathbb{N}$ such that $x=\sup_nf(y)(n)$.
If $(M,d,(a_n)_{n\in\mathbb{N}})$ is a computable metric space, there is a computable list $(B_e)_{e\in\mathbb{N}}$ of open balls of the form $B(a_n;q)$, where $B(a_n;q)$ is the open ball of radius $q$ centered at $a_n$.
We say that a set $P\subseteq M$ is $\Pi^0_1$ if there is a c.e.\ set $W\subseteq\mathbb{N}$ such that $P=M\setminus\bigcup_{e\in W}B_e$.
By $\Pi^0_1(z)$, we mean $\Pi^0_1$ relative to an oracle $z$, which is defined using a $z$-c.e.\ set $W$ instead of a c.e.\ set.
See also Section \ref{subsec:pre-repsp} for an abstract definition of $\Pi^0_1$ which is applicable to general spaces.

\subsubsection{Degree structures}

The Medvedev degrees $\mathfrak{M}$ \cite{medvedev} are a cornerstone of our framework. These are obtained by taking equivalence classes from Medvedev reducibility $\leq_M$, defined on subsets $A$, $B$ of Baire space $\Baire$ via $A \leq_M B$ iff there is a computable function $F : B \to A$. Important substructures of $\mathfrak{M}$ also relevant to us are the Turing degrees $\mathcal{D}_T$, the continuous degrees $\mathcal{D}_r$ and the enumeration degrees $\mathcal{D}_e$, these satisfy $\mathcal{D}_T \subsetneq \mathcal{D}_r \subsetneq \mathcal{D}_e \subsetneq \mathfrak{M}$.

Turing degrees are obtained from the usual Turing reducibility $\leq_T$ defined on points $p, q \in \Baire$ with $p \leq_T q$ iff there is a computable function $F : \subseteq \Baire \to \Baire$ with $F(q) = p$. We thus see $p \leq_T q \Leftrightarrow \{p\} \leq_M \{q\}$, and can indeed understand the Turing degrees to be a subset of the Medvedev degrees. The continuous degrees were introduced by \name{Miller} in \cite{miller2}. Enumeration degrees have received a lot of attention in computability theory, and were originally introduced by \name{Friedberg} and \name{Rogers} \cite{friedberg} (see also \cite[Chapter XIV]{OdiBook1}). In both cases, we can provide a simple definition directly as a substructure of the Medvedev degrees later on.

A further reducibility notion is relevant, although we are not particularly interested in its degree structure.
This is Muchnik reducibility $\leq_w$ \cite{muchnik}, defined again for sets $A, B \subseteq \Baire$ via $A \leq_w B$ iff, for any $p \in B$, there is $q \in A$ such that $q \leq_T p$.
Clearly $A \leq_M B$ implies $A \leq_w B$, but the converse is false in general.
A detailed investigation on the difference between Medvedev and Muchnik degrees can be found in \cite{kihara3,kihara3b,Stuk07}.



\subsubsection{Preliminaries for Sections \ref{sec:internalcharacterization}--\ref{sec:quasi-Polish}}\label{subsec:pre-repsp}

We briefly present some fundamental concepts on represented spaces following \cite{pauly-synthetic-arxiv}.
An important feature of the category of admissibly represented spaces and continuous functions is that it is cartesian closed, that is, closed under exponential, based on the \textrm{UTM}-theorem.
More explicitly, given represented spaces $\repsp{X}$ and $\repsp{Y}$, recall that every partial continuous function $f:\subseteq\repsp{X}\to\repsp{Y}$ is of the form $\Phi_e^z$ for some index $e\in\mathbb{N}$ and oracle $z$, and then, one can think of the pair $\langle e,z\rangle$ as a name of $f$.
That is, the map $\langle e,z\rangle\mapsto \Phi_e^z$ yields a representation of the space $\mathcal{C}(\repsp{X},\repsp{Y})$ of continuous functions between $\repsp{X}$ and $\repsp{Y}$ so that function evaluation and the other usual notions are computable.

In the following, we will want to make use of a special represented space, the Sierpi\'nski space $\mathbb{S} = (\{\bot, \top\}, \delta_\mathbb{S})$. The representations is given by $\delta_\mathbb{S}(0^\mathbb{N}) = \bot$ and $\delta_\mathbb{S}(p) = \top$ for $p \neq 0^\mathbb{N}$.
We then have the space $\mathcal{O}(\repsp{X}) \cong \mathcal{C}(\repsp{X},\mathbb{S})$ of open subsets of a represented space $\repsp{X}$ by identifying a set with its characteristic function, and the usual set-theoretic operations on this space are computable, too. We write $\mathcal{A}(\repsp{X})$ for the space of closed subsets, where names are names of the open complement.
Traditionally in computability theory, a computable element of the hyperspace $\mathcal{O}(\repsp{X})$ is called a {\em $\Sigma^0_1$ set}, a {\em $\Sigma^0_1$ class} or a {\em c.e.~open set}, and a computable element of the hyperspace $\mathcal{A}(\repsp{X})$ is called a {\em $\Pi^0_1$ set}, a {\em $\Pi^0_1$ class} or a {\em co-c.e.~closed set}.

The canonic function $\kappa_\repsp{X} : \repsp{X} \to \mathcal{O}(\mathcal{O}(\repsp{X}))$ mapping $x$ to $\{U \in \mathcal{O} \mid x \in U\}$ is always computable. If it has a computable inverse, then we call $\repsp{X}$ {\em computably admissible}. Admissibility in this sense was introduced by \name{Schr\"oder} \cite{schroder,schroder5}. As mentioned above, the computably admissible represented spaces are those that can be understood fully as topological spaces.

\subsection{Topology and Dimension}

\subsubsection{Isomorphism and Classification}\label{section:classification-problems}


We are now interested in isomorphisms of a particular kind, this always means a bijection in that function class, such that the inverse is also in that function class.
For instance, consider the following morphisms.
For a function $f:\repsp{X}\to\repsp{Y}$,
\begin{enumerate}
\item $f$ is {\em $\sigma$-computable} ({\em $\sigma$-continuous}, resp.)\ if there are sets $(X_n)_{n \in \mathbb{N}}$ such that $\repsp{X} = \bigcup_{n \in \mathbb{N}} X_n$ and each $f|_{X_n}$ is computable (continuous, resp.)
\item $f$ is {\em $\mathbf{\Gamma}$-piecewise continuous} if there are $\mathbf{\Gamma}$-sets $(X_n)_{n \in \mathbb{N}}$ such that $\repsp{X} = \bigcup_{n \in \mathbb{N}} X_n$ and each $f|_{X_n}$ is continuous.
\item $f$ is {\em $n$-th level Borel measurable} if $f^{-1}[A]$ is $\mathbf{\Sigma}^0_{n+1}$ for every $\mathbf{\Sigma}^0_{n+1}$ set $A\subseteq\repsp{X}$.
\end{enumerate}
In particular, $f$ is {\em second-level Borel measurable} iff $f^{-1}[A]$ is $G_{\delta\sigma}$ for every $G_{\delta\sigma}$ set $A\subseteq\repsp{X}$.
We also say that $f$ is {\em finite-level Borel measurable} if it is $n$-th level Borel measurable for some $n\in\mathbb{N}$.
Note that $\sigma$-continuity is also known as {\em countable continuity}.

Note that if $\repsp{X}$ and $\repsp{Y}$ have uniformly proper representations (this includes all computable metric spaces, see Subsection \ref{subsec:internalcharacterization} for details), then the $\repsp{X}_{n}$ in the definition of $\sigma$-continuity may be assumed to be $\mathbf{\Pi}_2^0$-sets.
Moreover, by recent results from descriptive set theory (see \cite{KiNg,MRos13,PaSa12}), we have the following implication for functions on Polish spaces:
\[\mbox{$\mathbf{\Pi}^0_2$-piecewise continuous $\Rightarrow$ second-level Borel measurable $\Rightarrow$ countably continuous}\]

Consequently, the second-level Borel isomorphic classification and the $\sigma$-continuous isomorphic classification of Polish spaces are exactly the same.
More precisely, three classification problems, Problems \ref{mainproblem1}, \ref{mainproblem2} and \ref{prob:third} in Section \ref{sec:introduction} are equivalent.

Hereafter, for notation, let $\cong$ be computable isomorphism, $\cong^\mathfrak{T}$ continuous isomorphism (i.e., homeomorphism), $\cong_{\sigma}$ be isomorphism by $\sigma$-computable functions and $\cong_{\sigma}^\mathfrak{T}$ is $\sigma$-continuous isomorphism.
We also use the terminologies such as {\em $\sigma$-homeomorphism} and {\em $\sigma$-embedding} to denote $\sigma$-continuous isomorphism and $\sigma$-continuous embedding.

For any of these notions, we write $\repsp{X} \leq \repsp{Y}$ with the same decorations on $\leq$ if $\repsp{X}$ is isomorphic to a subspace of $\repsp{Y}$ (i.e., $\repsp{X}$ is embedded into $\repsp{Y}$) in that way. If $\repsp{X} \leq \repsp{Y}$ and $\repsp{X}$ is not isomorphic to $\repsp{Y}$ in the designated way, then we also write $\repsp{X} < \repsp{Y}$, again with the suitable decorations on $<$. If neither $\repsp{X} \leq \repsp{Y}$ nor $\repsp{Y} \leq \repsp{X}$, we write $\repsp{X} \ | \ \repsp{Y}$ (again, with the same decorations).
The Cantor-Bernstein argument shows the following.

\begin{observation}\label{obs:Cantor-Bernstein}
Let $\repsp{X}$ and $\repsp{Y}$ be represented spaces.
Then, $\repsp{X}\cong_\sigma\repsp{Y}$ if and only if $\repsp{X}\leq_\sigma\repsp{Y}$ and $\repsp{Y}\leq_\sigma\repsp{X}$
\end{observation}

\subsubsection{Topological Dimension theory}
\label{subsec:dimension-theory-intro}

As general source for topological dimension theory, we point to \name{Engelking} \cite{EngBook}.
See also \name{van Mill} \cite{vMBook2} for infinite dimensional topology.
A topological space $\repsp{X}$ is \emph{countable dimensional} if it can be written as a countable union of finite dimensional subspaces.
Recall that a Polish space is countable dimensional if and only if it is {\em transfinite dimensional}, that is, its transfinite small inductive dimension is less than $\omega_1$ (see \cite[pp.~50--51]{hurewicz}).
One can see that a Polish space $\repsp{X}$ is countable dimensional if and only if $\repsp{X}\leq_\sigma^\mathfrak{T}\Cantor$.

To investigate the structure of uncountable dimensional spaces, \name{Alexandrov} introduced the notion of weakly/strongly infinite dimensional space.
We say that $C$ is a {\em separator} (usually called a {\em partition} in dimension theory) of a pair $(A,B)$ in a space $\repsp{X}$ if there are two pairwise disjoint open sets $A'\supseteq A$ and $B'\supseteq B$ such that $A'\sqcup B'=\repsp{X}\setminus C$.
A family $\{(A_i,B_i)\}_{i\in\Lambda}$ of pairwise disjoint closed sets in $\repsp{X}$ is {\em essential} if whenever $C_i$ is a separator of $(A_i,B_i)$ in $\repsp{X}$ for every $i\in\mathbb{N}$, $\bigcap_{i\in\mathbb{N}}C_i$ is nonempty.
A space $X$ is said to be {\em strongly infinite dimensional} if it has an essential family of infinite length.
Otherwise, $X$ is said to be {\em weakly infinite dimensional}.

We also consider the following covering property for topological spaces.
Let $\mathcal{O}[{\repsp{X}}]$ be the collection of all open covers of a topological space $\repsp{X}$, and $\mathcal{O}_2[{\repsp{X}}]=\{\mathcal{U}\in\mathcal{O}[X] : |\mathcal{U}|=2\}$, i.e.~the collection of all covers by two open sets.
Then, $\repsp{X}\in\mathcal{S}_c(\mathcal{A},\mathcal{B})$ if for any sequence $(\mathcal{U}_n)_{n\in\mathbb{N}}\in\mathcal{A}[\repsp{X}]^\mathbb{N}$, there is a sequence $(\mathcal{V}_n)_{n\in\mathbb{N}}$ of pairwise disjoint open sets such that $\mathcal{V}_n$ refines $\mathcal{U}_n$ for each $n\in\mathbb{N}$ and $\bigcup_{n\in\mathbb{N}}\mathcal{V}_n\in\mathcal{B}[\repsp{X}]$.

Note that a topological space $\repsp{X}$ is weakly infinite dimensional if and only if $\repsp{X}\in\mathcal{S}_c(\mathcal{O}_2,\mathcal{O})$.
We say that $\repsp{X}$ is a {\em $C$-space} \cite{AdGr78,Hav74} or {\em selectively screenable} \cite{Bab05} if $\repsp{X}\in\mathcal{S}_c(\mathcal{O},\mathcal{O})$.
We have the following implications:
\[\mbox{countable dimensional }\Rightarrow\mbox{ $C$-space }\Rightarrow\mbox{ weakly infinite dimensional}.\]

\name{Alexandrov}'s old problem was whether there exists a weakly infinite dimensional compactum $\repsp{X}>_\sigma^\mathfrak{T}\Cantor$.
This problem was solved by \name{R.~Pol} \cite{pol2} by constructing a compact metrizable space of the form $R\cup L$ for a strongly infinite dimensional totally disconnected subspace $R$ and a countable dimensional subspace $L$.
Such a compactum is a $C$-space, but not countable-dimensional.
Namely, \name{R.~Pol}'s theorem says that there are at least two $\sigma$-homeomorphism types of compact metrizable $C$-spaces.

There are previous studies on the structure of continuous isomorphism types (Fr\'echet dimension types) of various kinds of infinite dimensional compacta, e.g., strongly infinite dimensional Cantor manifolds (see \cite{Cha99,ChaEP99}).
For instance, by combining the Baire category theorem and the result by Chatyrko-Pol \cite{ChaEP99}, one can show that there are continuum many {\em first-level} Borel isomorphisms types of strongly infinite dimensional Cantor manifolds.
However, there is an enormous gap between first- and second-level, and hence, such an argument never tells us anything about second-level Borel isomorphism types.
Concerning weakly infinite dimensional Cantor manifolds, \name{El{\.z}bieta Pol} \cite{EPol96} (see also \cite{Cha99}) constructed a compact metrizable $C$-space in which no separator of nonempty subspaces can be hereditarily weakly infinite dimensional.
We call such a space a {\em Pol-type Cantor manifold}.


\section{Point Degree Spectra}\label{sec:pointdegreespectra}

\subsection{Generalized Turing Reducibility}\label{subsection:generalized-Turing-reducibility}

Recall that the notion of a represented space involves the notion of computability.
More precisely, every point in a represented space is coded by an infinite word, called a {\em name}.
Then, we estimate how complicated a given point is by considering the {\em degree of difficulty of calling a name} of the point.
Of course, it is possible for each point to have many names, and this feature yields the phenomenon that there is a point with no easiest names with respect to Turing degree.

Formally, we associate analogies of Turing reducibility and Turing degrees with an arbitrary represented space in the following manner.

\begin{definition}
Let $\repsp{X}$ and $\repsp{Y}$ be represented spaces.
We say that $y\in\repsp{Y}$ is {\em point-Turing reducible} to $x\in\repsp{X}$ if there is a partial computable function $f:\subseteq\repsp{X}\to\repsp{Y}$ such that $f(x)=y$, that is, $\delta^{-1}_Y(y)$ is Medvedev reducible to $\delta^{-1}_X(x)$.
In this case, we write $y^\repsp{Y}\leq_Mx^\repsp{X}$, or simply, $y\leq_Mx$.
\end{definition}

Roughly speaking, by the condition $y\leq_Mx$ we mean that if one knows a name of $x$, one can call a name of $y$, in a uniformly computable manner.
This pre-ordering relation $\leq_M$ clearly yields an equivalence relation $\equiv_M$ on points $x^\repsp{X}$ of represented spaces, and we then call each equivalence class $[x^\repsp{X}]_{\equiv_M}$ the {\em point-Turing degree} of $x\in\repsp{X}$, denoted by $\deg(x^\repsp{X})$.
In other words,
\[\deg(x^\repsp{X}) = [\delta_{X}^{-1}(x)]_{\equiv_M} = \mbox{``the Medvedev degree of the set of all $\delta_{X}$-names of $x$.''}\]
Then, we introduce the notion of point degree spectrum of a represented space as follows.

\begin{definition}\label{def:point-degree-spectra}
For a represented space $\repsp{X}$ and a point $x\in\repsp{X}$, define
\[\textrm{Spec}(\repsp{X}) =\{ \deg(x^\repsp{X}) \mid x \in \repsp{X}\} \subseteq \mathfrak{M}.\]
We call $\textrm{Spec}(\repsp{X})$ the \emph{point degree spectrum} of $\repsp{X}$.
Given an oracle $p$, we also define the relativized point degree spectrum by $\deg^p(x^\repsp{X})=[\{p\}\times\delta^{-1}_\repsp{X}(x)]_{\equiv_M}$ and ${\rm Spec}^p(\repsp{X})=\{\deg^p(x^\repsp{X}):x\in\repsp{X}\}$.
\end{definition}


Clearly, one can identify the Turing degrees $\mathcal{D}_T$, the continuous degrees $\mathcal{D}_r$ and the enumeration degrees $\mathcal{D}_e$ with degree spectra of some spaces as follows:
\begin{itemize}
\item $\spec(\Cantor) = \spec(\Baire) = \spec(\mathbb{R}) = \mathcal{D}_T$,
\item (\name{Miller} \cite{miller2}) $\spec(\uint^\mathbb{N}) = \spec(\mathcal{C}(\uint,\uint)) = \mathcal{D}_r$,
\item $\spec(\mathcal{O}(\mathbb{N})) = \mathcal{D}_e$, where $\mathcal{O}(\mathbb{N})$ is the space of all subsets of $\mathbb{N}$ where a basic open set is the set of all supersets of a finite subset of $\mathbb{N}$.
\end{itemize}

As any separable metric space embeds into the Hilbert cube $\uint^\mathbb{N}$, we find in particular that $\spec(\repsp{X}) \subseteq \mathcal{D}_r$ for any computable metric space $\repsp{X}$.
As any second-countable $T_0$ spaces embeds into the Scott domain $\mathcal{O}(\mathbb{N})$, we also have that $\spec(\repsp{X}) \subseteq \mathcal{D}_e$ for any second-countable $T_0$ space $\repsp{X}$.
In the latter case, the point degree of $x\in\repsp{X}$ corresponds to the enumeration degree of neighborhood basis as in Example \ref{example:representation}.
The Turing degrees will be characterized in Section \ref{subsec:spectra-dimensionth} in the context of topological dimension theory.

In computable model theory, the {\em degree spectrum} of a countable structure $S$ is defined as the collection of Turing degrees of isomorphic copies of $S$ coded in $\mathbb{N}$ (see \cite{HKSS02,Rich81}).
The notion of {\em degree spectrum on a cone} (i.e., degree spectrum relative to an oracle) also plays an important role in (computable) model theory (see \cite{Monta13,Mont14}).
One can define the space of countable structures as done in invariant descriptive set theory; however, from this perspective, a countable structure is a point, and therefore, the degree spectrum of a structure is a kind of the degree spectrum of a point rather than that of a space.

Given a point $x\in\repsp{X}$, we define ${\rm Spec}(x^\repsp{X})$ as the set of all oracles $z\in\{0,1\}^\mathbb{N}$ which can compute a name of $x$, and ${\rm Spec}^p(x^\repsp{X})$ as its relativization by an oracle $p\in\{0,1\}^\mathbb{N}$.
Then, the {\em weak point degree spectrum} ${\rm Spec}_w(\repsp{X})$ is the collection of all degree spectra of points of $x\in\repsp{X}$, and ${\rm Spec}^p_w(\repsp{X})$ is its relativization by an oracle $p$, that is,
\begin{align*}
&{\rm Spec}(x^\repsp{X})=\{z\in\{0,1\}^\mathbb{N}:x\leq_Mz\},& &{\rm Spec}^p(x^\repsp{X})=\{z\in\{0,1\}^\mathbb{N}:x\leq_M(z,p)\},\\
&{\rm Spec}_w(\repsp{X})=\{{\rm Spec}(x^\repsp{X}):x\in\repsp{X}\},& &{\rm Spec}^p_w(\repsp{X})=\{{\rm Spec}^p(x^\repsp{X}):x\in\repsp{X}\}.
\end{align*}

Note that this notion can be described in terms of Muchnik redcibility, that is, we can think of the degree spactrum of $x\in\repsp{X}$ as:
\[{\rm Spec}(x^\repsp{X})\approx[\delta_X(x)]_{\equiv_w}=\mbox{``the Muchnik degree of the set of all $\delta_X$-names of $x$.''}\]

\begin{observation}\label{obs:muchnik-medvedev}
If $\repsp{X}$ and $\repsp{Y}$ are admissibly represented second-countable $T_0$-spaces, then there is an oracle $p$ such that for all $q\geq_Tp$,
\[{\rm Spec}^q(\repsp{X})\subseteq{\rm Spec}^q(\repsp{Y})\iff{\rm Spec}^q_w(\repsp{X})\subseteq{\rm Spec}^q_w(\repsp{Y})\]
\end{observation}

\begin{proof}
The point degree spectrum of an admissibly represented second-countable space can be thought of as a substructure of the enumeration degrees $\mathcal{D}_e$ as mentioned above.
It is known that enumeration reducibility coincides with its non-uniform version (see \cite{selman} or \cite[Theorem 4.2]{miller2}).
This means that $x^\repsp{X}\leq_My^\repsp{Y}$ if and only if every $\delta_Y$-name of $y$ computes a $\delta_X$-name of $x$ (see also \cite[Corollary 4.3]{miller2}).
\end{proof}

We now focus on degree spectra of separable metrizable spaces, that is, continuous degrees.
The following lemma shows -- in \name{Miller}'s words -- that the continuous degrees are \emph{almost} Turing degrees.
To be more precise, any continuous degree is relativized into a Turing degree by all Turing degrees except the smaller ones.

\begin{lemma}[\name{Miller}]\label{lem:JM-almosttotal}
For any non-total continuous degree $q \in \mathcal{D}_r \setminus \mathcal{D}_T$ we find that for all $p \in \mathcal{D}_T$, $(p,q) \in \mathcal{D}_T$ iff $p \nleq_M q$.
\begin{proof}
Let $r=(r(n))_{n \in \mathbb{N}} \in \uint^\mathbb{N}$ be a representative of a non-total continuous degree $q \in \mathcal{D}_r \setminus \mathcal{D}_T$.
Let $I$ be the set of all $y\in\Cantor$ such that $y\leq_Mr$, which is a countable set.
Choose a real $x$ whose Turing degree is incomparable with $I$.
In particular, $x$ is algebraically transcendent with all reals in $I$.
So, there is an $x$-computable homeomorphism sending $r$ to a sequence of irrationals.
Hence, given any name of $(x, r)$, we first obtain $x$, and by using $x$, transform $r$ into irrationals, and then we get the least Turing degree name of $(x,r)$.
\end{proof}
\end{lemma}

\subsection{Degree Spectra and Dimension Theory}\label{subsec:spectra-dimensionth}

One of the main tools in our work is the following characterization of the point degree spectra of represented spaces.

\begin{theorem}
\label{theo:spectrum-main}
The following are equivalent for admissibly represented spaces $\repsp{X}$ and $\repsp{Y}$:
\begin{enumerate}
\item ${\rm Spec}^r(\repsp{X})={\rm Spec}^r(\repsp{Y})$ for some oracle $r\in\Cantor$.
\item $\mathbb{N}\times\repsp{X}$ is $\sigma$-homeomorphic to $\mathbb{N}\times\repsp{Y}$, i.e., $\mathbb{N}\times\repsp{X}\cong_\sigma^\mathfrak{T}\mathbb{N}\times\repsp{Y}$.
\end{enumerate}

Moreover, if $\repsp{X}$ and $\repsp{Y}$ are Polish, then the following assertions (3) and (4) are also equivalent to the above assertions (1) and (2).
\begin{enumerate}
\item[3.] $\mathbb{N}\times\repsp{X}$ is second-level Borel isomorphic to $\mathbb{N}\times\repsp{Y}$.
\item[4.] The Banach algebra $\mathcal{B}_2^*(\mathbb{N}\times\repsp{X})$ is linearly isometric (ring isomorphic and so on) to $\mathcal{B}_2^*(\mathbb{N}\times\repsp{Y})$.
\end{enumerate}
\end{theorem}

One can also see that the following assertions are equivalent:

\begin{enumerate}
\item[2$^\prime$.] $\mathbb{N}\times\repsp{X}$ is $G_\delta$-piecewise homeomorphic to $\mathbb{N}\times\repsp{Y}$.
\item[3$^\prime$.] $\mathbb{N}\times\repsp{X}$ is $n$-th level Borel isomorphic to $\mathbb{N}\times\repsp{Y}$ for some $n\geq 2$.
\item[4$^\prime$.] The Banach algebra $\mathcal{B}_n^*(\mathbb{N}\times\repsp{X})$ is linearly isometric (ring isomorphic and so on) to $\mathcal{B}_n^*(\mathbb{N}\times\repsp{Y})$ for some $n\geq 2$.
\end{enumerate}

By our argument in Section \ref{section:classification-problems}, the assertions (2$^\prime$) is equivalent to (2).
Obviously the assertions (3) and (4) imply (3$^\prime$) and (4$^\prime$), respectively.
The equivalence between (3) and (4) (and the equivalence between (3$^\prime$) and (4$^\prime$)) has already been shown by Jayne \cite{Jayne74} for second-countable (or more generally, realcompact) spaces $\repsp{X}$ and $\repsp{Y}$.
The implication from the assertion (3$^\prime$) to (2) is, as mentioned in Section \ref{section:classification-problems}, recently proved by \cite{MRos13,PaSa12}, and more recently, an alternative computability-theoretic proof is given by \cite{KiNg} using our framework of point degree spectra of Polish spaces.
Consequently, all assertions from (2) to (4$^\prime$) are equivalent.

To see the equivalence between (1) and (2), we characterize the point degree spectra of represented spaces in the context of countably-continuous isomorphism.

\begin{lemma}\label{lemma:spectrum-main}
The following are equivalent for represented spaces $\repsp{X}$ and $\repsp{Y}$:
\begin{enumerate}
\item $\textrm{Spec}(\repsp{X}) \subseteq \textrm{Spec}(\repsp{Y})$
\item $\repsp{X}\leq_\sigma\mathbb{N}\times\repsp{Y}$, i.e., $\repsp{X}$ is a countable union of subspaces that are computably isomorphic to subspaces of $\repsp{Y}$.
\end{enumerate}
\end{lemma}

\begin{proof}
We first show that the assertion (1) implies (2).
By assumption, for any $x \in \repsp{X}$ we find $\delta_{\repsp{X}}^{-1}(x) \equiv_M \delta_{\repsp{Y}}^{-1}(y_x)$ for some $y_x \in \repsp{Y}$. Let for $\repsp{Y}$ any $i, j \in \mathbb{N}$, let $\repsp{X}_{ij}$ be the set of all points where the reductions are witnessed by $\Phi_i$ and $\Phi_j$, and let $\repsp{Y}_{ij} = \{y_x \mid x \in \repsp{X}_{ij}\} \subseteq \repsp{Y}$, where recall that $\Phi_e$ is the $e$-th partial computable function. Then $\Phi_i$, $\Phi_j$ also witness $\repsp{X}_{ij} \cong \repsp{Y}_{ij}$, and obviously $\repsp{X} = \bigcup_{\langle i, j\rangle \in \mathbb{N}} \repsp{X}_{ij}$.

Conversely, the point spectrum is preserved by computable isomorphism and $\textrm{Spec}\left (\bigcup_{n \in \mathbb{N}} \repsp{X}_n\right ) = \bigcup_{n \in \mathbb{N}} \textrm{Spec}(\repsp{X}_n)$, so the claim follows.
\end{proof}

\begin{proof}[Proof of Theorem \ref{theo:spectrum-main} (1) $\Leftrightarrow$ (2).]
It follows from relativizations of Lemma \ref{lemma:spectrum-main} and Observation \ref{obs:Cantor-Bernstein}.
Here, it is easy to see that the assertion (2) is equivalent to $\mathbb{N}\times\repsp{X}\leq_\sigma\mathbb{N}\times\repsp{Y}$.
\end{proof}

This simple argument completely solves a mystery about the occurrence of non-Turing degrees in proper infinite dimensional spaces.
Concretely speaking, by Lemma \ref{lemma:spectrum-main}, we can characterize the Turing degrees in terms of topological dimension theory as follows\footnote{The same observation was independently made by \name{Hoyrup}. \name{Brattka} and \name{Miller} had conjectured that dimension would be the crucial demarkation line for spaces with only Turing degrees (all personal communication).}.

\begin{corollary}
The following are equivalent for a separable metrizable space $\repsp{X}$ endowed with an admissible representation:
\begin{enumerate}
\item $\textrm{Spec}^p(\repsp{X}) \subseteq \mathcal{D}_T$ for some oracle $p\in\Cantor$
\item $\repsp{X}$ is countable dimensional.
\end{enumerate}
\end{corollary}

By a dimension-theoretic fact (see Section \ref{subsec:dimension-theory-intro}), if $\repsp{X}$ is Polish, transfinite dimensionality is also equivalent to the condition for $\repsp{X}$ in which any point has a Turing degree relative to some oracle.

Now, by Theorem \ref{theo:spectrum-main}, $\sigma$-homeomorphic classification can be viewed as a kind of degree theory dealing with the {\em order structure on degree structures} (on a cone).
Thus, from the viewpoint of degree theory, it is natural to ask whether Post's problem (there is an intermediate degree structure strictly between the bottom $\Cantor$ and the top $[0,1]^\mathbb{N}$), the Friedberg-Muchnik theorem (there is an incomparable degree structures), the Sacks density theorem (given two comparable, but different degree structures, there is an intermediate degree structure strictly between them), and so on, is true for degrees of degrees of uncountable Polish spaces.

More details of the structure of degree spectra of Polish space will be investigated in Sections \ref{sec:intermediate} and \ref{sec:piecewisehomeo}, and those of quasi-Polish space will be in Section \ref{sec:quasi-Polish}.

\section{Intermediate Point Degree Spectra}\label{sec:intermediate}

\subsection{Intermediate Polish Spaces}

Let $\mathfrak{P}$ be the set of all uncountable Polish spaces.
In this section, we investigate the structure of $\mathfrak{P}/\cong_{\sigma}^\mathfrak{T}$, i.e.~either of the equivalence classes w.r.t.~$\sigma$-homeomorphisms, or equivalently, the structure of point degree spectra of uncountable Polish spaces up to relativization.

It is well-known that for every uncountable Polish space $X$:
\[\Cantor\leq_{c}^\mathfrak{T} X\leq_{c}^\mathfrak{T}[0,1]^\mathbb{N},\]
where, recall that $\leq_{c}^\mathfrak{T}$ is the topological embeddability relation (i.e., the ordering of Fr\'echet dimension types).
In this section, we focus on Problem \ref{prob:third} asking whether there exists a Polish space $\repsp{X}$ satisfying the following:
\[\Cantor <_{\sigma}^\mathfrak{T}\repsp{X}<_{\sigma}^\mathfrak{T}[0,1]^\mathbb{N}.\]

One can see that there is no difference between the structures of $\sigma$-homeomorphism types of uncountable Polish spaces and uncountable compact metric spaces.

\begin{fact}\label{thm:equivalence}
Every Polish space is $\sigma$-homeomorphic to a compact metrizable space.
\end{fact}

\begin{proof}
All spaces of a given countable cardinality are clearly $\sigma$-homeomorphic, and there are compact metrizable spaces of all countable cardinalities.

So let $\repsp{X}$ be an uncountable Polish space.
\name{Lelek} \cite{lelek} showed that every Polish space $\repsp{X}$ has a compactification $\gamma \repsp{X}$ such that $\gamma \repsp{X}\setminus \repsp{X}$ is countable-dimensional.
Clearly $\repsp{X}\leq_{c}\gamma \repsp{X}$.
Then, we have $\gamma \repsp{X}\setminus \repsp{X}\leq_{\sigma}^\mathfrak{T}\Cantor\leq_{\sigma}^\mathfrak{T}\repsp{X}$, since $\repsp{X}$ is uncountable Polish and $\gamma \repsp{X}\setminus \repsp{X}$ is countable-dimensional.
Consequently, $\repsp{X},\gamma \repsp{X}\setminus \repsp{X}\leq_{\sigma}^\mathfrak{T}\repsp{X}$, and this implies $\gamma \repsp{X}=\repsp{X}\cup(\gamma \repsp{X}\setminus \repsp{X})\leq_{\sigma}^\mathfrak{T}\repsp{X}$.
\end{proof}

\subsection{The Graph Space of a Universal $\omega$-Left-CEA Operator}

Now, we provide a concrete example having an intermediate degree spectrum.
We say that a point $(r_n)_{n\in\mathbb{N}}\in[0,1]^\mathbb{N}$ is {\em $\omega$-left-CEA in} or an {\em $\omega$-left-pseudojump of} $x\in \Cantor$ if $r_{n+1}$ is left-c.e.~in $\langle x,r_0,r_1,\dots,r_n\rangle$ uniformly in $n\in\mathbb{N}$.
In other words, there is a computable function $\Psi:\Cantor\times[0,1]^{<\omega}\times\mathbb{N}^2\to\mathbb{Q}_{\geq 0}$ such that
\[r_{n}=\sup_{s\to\infty}\Psi(x,r_0,\dots,r_{n-1},n,s)\]
for every $x,n,s$, where $\mathbb{Q}_{\geq 0}$ denotes the set of all nonnegative rationals.
Whenever $r_n\in[0,1]$ for all $n\in\mathbb{N}$, such a computable function $\Psi$ generates an operator $J_\Psi^\omega:\Cantor\to[0,1]^\mathbb{N}$ with $J_\Psi^\omega(x)=(r_0,r_1,\dots)$, which is called an {\em $\omega$-left-CEA operator}.

\begin{proposition}\label{prop:universaloperator}
There is an effective enumeration $(J^\omega_e)_{e\in\mathbb{N}}$ of all $\omega$-left-CEA operators.
\end{proposition}

\begin{proof}
It is not hard to see that $y\in[0,1]$ is left-c.e.~in $x\in\Cantor\times[0,1]^k$ if and only if there is a c.e.~set $W\subseteq\mathbb{N}\times \mathbb{Q}$ such that
\[y=J^{k}_W(x):=\sup\{\min\{|p|,1\}:x\in B^k_i\mbox{ for some }(i,p)\in W\},\]
where $B_i^k$ is the $i$-th rational open ball in $[0,1]^k$.
Thus, we have an effective enumeration of all left-c.e.~operators $J:\Cantor\times[0,1]^k\rightarrow[0,1]$ by putting $J^k_e=J^k_{W_e}$, where $W_e$ is the $e$-th c.e. subset of $\mathbb{N}\times \mathbb{Q}$.
Then, we define
\[J^\omega_e(x)=(x,J^0_{\langle e,0\rangle}(x),J^1_{\langle e,1\rangle}(x,J^0_{\langle e,0\rangle}(x)),\dots),\]
that is, $J^\omega_e$ is the $\omega$-left-CEA operator generated by the uniform sequence $(J^k_{\langle e,k\rangle})_{k\in\mathbb{N}}$ of left-c.e.~operators.
Clearly, $(J^\omega_e)_{e\in\mathbb{N}}$ is an effective enumeration of all $\omega$-left-CEA operators.
\end{proof}

Hence, we may define a {\em universal $\omega$-left-CEA operator} by $J^\omega(e,x)=J^\omega_e(x)$.

\begin{definition}
The {\em $\omega$-left-computably-enumerable-in-and-above space} $\omega\repspb{CEA}$ is a subspace of $\mathbb{N}\times \Cantor\times[0,1]^\mathbb{N}$ defined by
\begin{align*}
\repspb{\omega CEA}&=\{(e,x,r)\in\mathbb{N}\times \Cantor\times[0,1]^\mathbb{N}:r=J^\omega_e(x)\}\\
&=\mbox{``the graph of a universal $\omega$-left-CEA operator.''}
\end{align*}
\end{definition}

Note that in classical recursion theory, an operator $\Psi$ is called a {\em CEA-operator} (also known as an {\em REA-operator} or a {\em pseudojump}) if there is a c.e.~procedure $W$ such that $\Psi(A)=\langle{A,W(A)\rangle}$ for any $A\subseteq\mathbb{N}$ (see \name{Odifreddi} \cite[Chapters XII and XIII]{OdiBook1}).
An $\omega$-CEA operator is the $\omega$-th iteration of a uniform sequence of CEA-operators.
In general, computability theorists have studied $\alpha$-CEA operators for computable ordinals $\alpha$ in the theory of $\Pi^0_2$ singletons.
We will also use a generalization of the notion of a $\Pi^0_2$ singleton in Section \ref{sec:piecewisehomeo}.

We say that a continuous degree is {\em $\omega$-left-CEA} if it contains a point $r\in[0,1]^\mathbb{N}$ which is $\omega$-left-CEA in a point $z\in \Cantor$ such that $z\leq_Mr$, i.e., $J^\omega_e(z)=r$ for some $e$.
The point degree spectrum of the space $\repspb{\omega CEA}$ (as a subspace of $[0,1]^\mathbb{N}$) can be described as follows.
\[\spec(\repspb{\omega CEA})=\{\mathbf{a}\in\mathcal{D}_r:\mbox{$\mathbf{a}$ is $\omega$-left-CEA}\}.\]

Clearly,
\[\spec({\Cantor})\subseteq\spec({\repspb{\omega CEA}})\subseteq\spec({[0,1]^\mathbb{N}}).\]

The following is an analog of the well-known fact from classical computability theory that every $\omega$-CEA set is a $\Pi^0_2$-singleton (see \name{Odifreddi} \cite[Proposition XIII.2.7]{OdiBook1})

\begin{lemma}
The $\omega$-left-CEA space $\repspb{\omega CEA}$ is Polish.
\end{lemma}

\begin{proof}
It suffices to show that $\repspb{\omega CEA}$ is $\Pi^0_2$.
The stage $s$ approximation to $J_{e}^k$ is denoted by $J_{e,s}^k$, that is, $J_{e,s}^k(z)=\max\{\min\{|p|,1\}:(\exists\langle i,p\rangle\in W_{e,s})\;x\in B_i^k\}$, where $W_{e,s}$ is the stage $s$ approximation to the $e$-th computably enumerable set $W_e$.
Note that the function $(e,s,k,z)\mapsto J_{e,s}^k(z)$ is computable.
We can easily see that $(e,x,r)\in\repspb{\omega CEA}$ if and only if
\[(\forall n,k\in\mathbb{N})(\exists s>n)\;d\left(\pi_{k}(r),J^k_{e,s}(x,\pi_0(r),\pi_1(r),\dots,\pi_{k-1}(r))\right)<2^{-n},\]
where $d$ is the Euclidean metric on $[0,1]$.
\end{proof}

We devote the rest of this section to a proof of the following theorem.

\begin{theorem}\label{thm:intermediatePolish}
The space $\repspb{\omega CEA}$ has an intermediate $\sigma$-homeomorphism type, that is,
\[\Cantor <_{\sigma}^\mathfrak{T} \repspb{\omega CEA} <_{\sigma}^\mathfrak{T} [0,1]^\mathbb{N}.\]
\end{theorem}

Consequently, the space $\repspb{\omega CEA}$ is a concrete counterexample to Problem \ref{prob:third}.

\subsection{Proof of $\repspb{\omega CEA}<_\sigma^\mathfrak{T}[0,1]^\mathbb{N}$}

The key idea is to measure {\em how similar the space $\repsp{X}$ is to a zero-dimensional space} by approximating each point in a space $\repsp{X}$ by a zero-dimensional space.
Recall from Observation \ref{obs:muchnik-medvedev} that the point degree spectrum coincides with the Muchnik degree spectrum for any second-countable admissibly represented space, that is, the point-Turing degree $\deg(x)$ of a point $x\in\repsp{X}$ can be identified with its Turing upper cone, that is,
\[\deg(x)\approx{\rm Spec}(x)=\{z\in \Cantor:x\leq_Mz\}.\]

We think of the spectrum ${\rm Spec}(x)$ as {\em the upper approximation of $x\in\repsp{X}$ by the zero-dimensional space $\Cantor$}.
Now, we need the notion of {\em the lower approximation of $x\in\repsp{X}$ by the zero-dimensional space $\Cantor$}.
We introduce the {\em co-spectrum} of a point $x\in\repsp{X}$ as its Turing lower cone
\[{\rm coSpec}(x)=\{z\in \Cantor:z\leq_Mx\},\]
and moreover, we define the {\em degree co-spectrum} of a represented space $\repsp{X}$ as follows:
\[{\rm coSpec}(\repsp{X})=\{{\rm coSpec}(x):x\in\repsp{X}\}.\]
Note that the degree spectrum of a represented space fully determines its co-spectrum, while the converse is not true. For every oracle $p\in \Cantor$, we may also introduce relativized co-spectra ${\rm coSpec}^p(x)=\{z\in \Cantor:z\leq_M(x,p)\}$, and the relativized degree co-spectra ${\rm coSpec}^p(\repsp{X})$ in the same manner.

\begin{observation}\label{obs:main-cospec-inv}
Let $\repsp{X}$ and $\repsp{Y}$ be admissibly represented spaces.
If ${\rm Spec}^p(\repsp{X})={\rm Spec}^p(\repsp{Y})$, then we also have ${\rm coSpec}^p(\repsp{X})={\rm coSpec}^p(\repsp{Y})$.
Therefore, by Theorem \ref{theo:spectrum-main}, the cospectrum of an admissibly represented space up to an oracle is invariant under $\sigma$-homeomorphism.
\end{observation}

We say that a collection $\mathcal{I}$ of subsets of $\mathbb{N}$ is {\em realized as the co-spectrum of $x$} if ${\rm coSpec}(x)=\mathcal{I}$.
A countable set $\mathcal{I}\subseteq\mathcal{P}(\mathbb{N})$ is a {\em Scott ideal} if it is the standard system of a countable nonstandard model of Peano arithmetic, or equivalently, a countable $\omega$-model of the theory ${\sf WKL}_0$.
We will not go into the details of a Scott ideal (see \name{Miller} \cite[Section 9]{miller2} for more explicit definition); we will only use the fact that every jump ideal is a Scott ideal.
Here, a {\em jump ideal} $\mathcal{I}$ is a collection of subsets of natural numbers which is closed under the join $\oplus$, downward Turing reducibility $\leq_T$, and the Turing jump, that is, $p,q\in\mathcal{I}$ implies $p\oplus q\in\mathcal{I}$; $p\leq_Tq\in\mathcal{I}$ implies $p\in\mathcal{I}$; and $p\in\mathcal{I}$ implies $p'\in\mathcal{I}$.
\name{Miller} \cite[Theorem 9.3]{miller2} showed that every countable Scott ideal (hence, every countable jump ideal) is realized as a co-spectrum in $[0,1]^\mathbb{N}$.

\begin{figure}[t]
 \begin{center}
  \includegraphics[width=100mm]{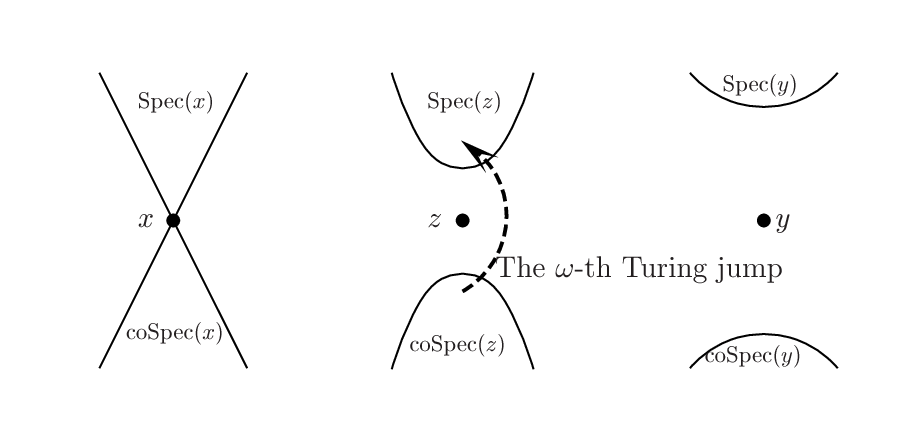}
 \end{center}
 \caption{The upper and lower approximations of $\Cantor$, $\omega\repspb{CEA}$ and $[0,1]^\mathbb{N}$}
 \label{fig:spec-cospec}
\end{figure}

\begin{example}
The spectra and co-spectra of Cantor space $\Cantor$, the space $\omega{\rm CEA}$, and the Hilbert cube $[0,1]^\mathbb{N}$ are illustrated as follows (see also Figure \ref{fig:spec-cospec}):
\begin{enumerate}
\item The co-spectrum ${\rm coSpec}(x)$ of any point $x\in\Cantor$ is principal, and {\em meets with ${\rm Spec}(x)$} exactly at $\deg_T(x)$.
The same is true up to some oracle for an arbitrary Polish spaces $\repsp{X}$ such that $\repsp{X}\cong_\sigma^\mathfrak{T}\Cantor$.
\item For any point $z\in\omega\repsp{CEA}$, the ``{\em distance}'' between ${\rm Spec}(z)$ and ${\rm coSpec}(z)$ has to be at most the $\omega$-th Turing jump (see Lemma \ref{lem:Scott-ideal}).
\item An arbitrary countable Scott ideal is realized as ${\rm coSpec}(y)$ of some point $y\in [0,1]^\mathbb{N}$.
Hence, ${\rm Spec}(y)$ and ${\rm coSpec}(y)$ can be separated by {\em an arbitrary distance}.
\end{enumerate}
\end{example}

This upper/lower approximation method clarifies the differences of $\sigma$-homeomorphism types of spaces because both relativized point-degree spectra and co-spectra are invariant under $\sigma$-homeomorphism by Theorem \ref{theo:spectrum-main} and Observation \ref{obs:main-cospec-inv}.

\begin{lemma}\label{lem:Scott-ideal}
For any oracle $p\in \Cantor$, there is a countable Scott ideal which cannot be realized as a $p$-co-spectrum of an $\omega$-left-CEA continuous degree.
\end{lemma}

\begin{proof}
Let $y=(e,x,r)\in\omega\repspb{CEA}$ be an arbitrary point.
Clearly, $x\leq_M(e,x,r)$, and this means that $x\in{\rm coSpec}(y)$ since $x\in\Cantor$.
However, $x^{(\omega)}\in{\rm Spec}(y)$ since $r$ is $\omega$-left-CEA in $x$.
Hence, ${\rm coSpec}(y)$ is not closed under the $\omega$-th Turing jump for any $y\in\omega\repspb{CEA}$.
Thus, for any oracle $p$, the jump ideal $\mathcal{A}^p=\{x\in \Cantor:(\exists n\in\mathbb{N})\;x\leq_T p^{(\omega\cdot n)}\}$ cannot be realized as a co-spectrum in $\omega\repspb{CEA}$.
\end{proof}

Consequently, the $\omega$-left-CEA space is not $\sigma$-homeomorphic to the Hilbert cube.
Note that \name{Day} and \name{Miller} \cite{daymiller} showed that every countable Scott ideal $\mathcal{I}$ is realized by a neutral measure.
Hence, we can also conclude that there is a neutral measure whose continuous degree is not $\omega$-left-CEA.





\subsection{Proof of $\Cantor<_\sigma^\mathfrak{T}\repspb{\omega CEA}$}\label{subsec:proof-main-half}

Next, we have to show that the $\omega$-left-CEA space is not countable-dimensional.
For a compact set $P\subseteq[0,1]^\mathbb{N}$, we inductively define $\min P\in P$ as follows:
\[\pi_n(\min P)=\min\pi_n[\{z\in P:(\forall i<n)\;\pi_i(z)=\pi_i(\min P)\}],\]
where $\pi_n:[0,1]^\mathbb{N}\to[0,1]$ is the projection onto the $n$-th coordinate.
We call the point $\min P$ the {\em leftmost point of $P$}.
Kreisel's basis theorem (see \cite[Proposition V.5.31]{OdiBook}) in classical computability theory says that the leftmost element of a $\Pi^0_1$ subset of $\Cantor$ or $[0,1]$ is always left-c.e.
The following lemma can be viewed as an infinite dimensional version of Kreisel's basis theorem.

\begin{lemma}\label{lem:basis}
For any oracle $p\in\Cantor$, the leftmost point of a $\Pi^0_1(p)$ subset of $[0,1]^\mathbb{N}$ is $\omega$-left-CEA in $p$.
\end{lemma}
\begin{proof}
We first note that Hilbert cube $[0,1]^\mathbb{N}$ is computably compact
in the sense that there is a computable enumeration of all finite collections $\mathcal{D}$ of basic open sets which covers the whole space, that is, $\bigcup\mathcal{D}=[0,1]^\mathbb{N}$.

Fix a $\Pi^0_1(p)$ set $P\subseteq[0,1]^\mathbb{N}$.
It suffices to show that $\pi_{n+1}(\min P)$ is left-c.e.~in $\langle\pi_{i}(\min P)\rangle_{i\leq n}$ uniformly in $n$ relative to $p$.
Given a sequence $\mathbf{a}=(a_0,a_1,\dots,a_n)$ of reals and an real $q$, we denote by $C(\mathbf{a},q)$ the set of all points in $P$ of the form $(a_0,a_1,\dots,a_n,r,\dots)$ for some $r\leq q$, that is,
\[C(\mathbf{a},q)=P\cap \bigcap_{i\leq n}\pi_i^{-1}\{a_i\}\cap\pi_{n+1}^{-1}[0,q].\]

By computable compactness of Hilbert cube, one can see that $C^*(\mathbf{a})=\{q\in [0,1]:C(\mathbf{a},q)=\emptyset\}$ is $p$-c.e.~open uniformly relative to $\mathbf{a}$ since the complement of $C(\mathbf{a},q)$ is $p$-c.e.~open uniformly relative to $\mathbf{a}$ and $q$.
Therefore, $\sup C^*(\mathbf{a})$ is $p$-left-c.e.~uniformly relative to $\mathbf{a}$.
Finally, we can easily see that  $\pi_{n+1}(\min X)$ is exactly $\sup C^*(\langle\pi_{i}(\min X)\rangle_{i\leq n})$.
\end{proof}

We use the following relativized versions of \name{Miller}'s lemmas.

\begin{lemma}[\name{Miller} {\cite[Lemma 6.2]{miller2}}]\label{lemma:Miller-6-2}
For every $p\in \Cantor$, there is a multivalued function $\Psi^p:[0,1]^\mathbb{N}\to[0,1]^\mathbb{N}$ with a $\Pi^0_1(p)$ graph and nonempty, convex images such that, for all $e\in\mathbb{N}$, $\alpha\in[0,1]^\mathbb{N}$ and $\beta\in\Psi^p(\alpha)$, if for every name $\lambda$ of $\alpha$, $\varphi_e^{\lambda\oplus p}$ is a name of $x\in [0,1]$, then $\beta(e)=x$.
\end{lemma}

Note that \name{Kakutani}'s fixed point theorem ensures the existence of a fixed point of $\Psi$.
If $\alpha$ is a fixed point of $\Psi^p$, that is, $\alpha\in\Psi^p(\alpha)$, then ${\rm coSpec}^p(\alpha)=\{\alpha(n):n\in\mathbb{N}\}$.
Therefore, such an $\alpha$ has no Turing degree relative to $p$ (see \cite[Proposition 5.3]{miller2}).

\begin{lemma}[\name{Miller} {\cite[Lemma 9.2]{miller2}}]\label{lem:Miller2}
For every $p\in \Cantor$, there is an index $e\in\mathbb{N}$ such that for any $x\in[0,1]$, there is a fixed point $\alpha$ of $\Psi^p$ such that $\alpha(e)=x$.
\end{lemma}

\begin{lemma}\label{thm:nontotalcea}
For any oracle $p\in\Cantor$, there is an $\omega$-left-CEA continuous degree which is not contained in ${\rm Spec}^p(\Cantor)$.
\end{lemma}

\begin{proof}
Let ${\rm Fix}(\Psi^p)$ be the set of all fixed points of $\Psi^p$.
Then, ${\rm Fix}(\Psi^p)$ is $\Pi^0_1(p)$ since it is the intersection of the graph of a $\Pi^0_1(p)$ set and the diagonal set.
Let $e$ be an index in Lemma \ref{lem:Miller2}.
Clearly, $A=\{\alpha\in{\rm Fix}(\Psi^p):\alpha(e)=p\}$ is again a $\Pi^0_1(p)$ subset of $[0,1]^\mathbb{N}$.
By Lemma \ref{lem:basis}, $A$ contains an element $\alpha$ which is $\omega$-left-CEA in $p$.
By the property of $A$ discussed above, $\alpha$ has no Turing degree relative to $p$.
\end{proof}

\begin{proof}[Proof of Theorem \ref{thm:intermediatePolish}]
By Lemma \ref{lem:Scott-ideal}, ${\rm coSpec}^p(\omega\repspb{CEA})\subsetneq{\rm coSpec}^p([0,1]^\mathbb{N})$ for any oracle $p$.
Moreover, by Lemma \ref{thm:nontotalcea}, ${\rm Spec}^p(\Cantor)\subsetneq{\rm Spec}^p(\omega\repspb{CEA})$ for any oracle $p$.
Therefore, by Theorem \ref{theo:spectrum-main} and Observation \ref{obs:main-cospec-inv}, we conclude $\Cantor<_\sigma^\mathfrak{T}\omega\repspb{CEA}<_\sigma^\mathfrak{T}[0,1]^\mathbb{N}$.
\end{proof}

\section{Structure of $\sigma$-Homeomorphism Types}\label{sec:piecewisehomeo}


In this section, we will show that there are continuum many $\sigma$-homeomorphism types of compact metrizable spaces.

\begin{theorem}\label{thm:maintheorem_emb}
There exists a collection $(\repsp{X}_\alpha)_{\alpha<2^{\aleph_0}}$ of continuum many compact metric spaces such that if $\alpha\not=\beta$, $\repsp{X}_\alpha$ cannot be $\sigma$-embedded into $\repsp{X}_\beta$.
\end{theorem}

We devote the rest of this section to prove Theorem \ref{thm:maintheorem_emb}.
Actually, we will show the following:
\begin{itemize}
\item [] There is an embedding of the inclusion ordering $([\omega_1]^{\leq\omega},\subseteq)$ of countable subsets of the smallest uncountable ordinal $\omega_1$ into the $\sigma$-embeddability ordering of compact metric spaces.
\end{itemize}

As a corollary, there are an uncountable chain and a continuum antichain of $\sigma$-homeomorphism types of compact metric spaces.

\subsection{Almost Arithmetical Degrees}

In Section \ref{sec:intermediate}, we used a co-spectrum as a $\sigma$-topological invariant.
More explicitly, in our proof, it was essential to examine closure properties of co-spectra to obtain an intermediate $\sigma$-homeomorphism type of Polish spaces.
In this section, we will develop a method for controlling closure properties of co-spectra.
As a result, we will construct a compact metrizable space whose co-spectra realize a given well-behaved family of {\em ``almost'' arithmetical degrees}.

First, we introduce a notion which estimates the strength of closure properties of functions up to the arithmetical equivalence.

\begin{definition}\label{def:almostarithmetical}
Let $g$ and $h$ be total Borel measurable functions from $\Cantor$ into $\Cantor$.
\begin{enumerate}
\item We inductively define $g^0(x)=x$ and $g^{n+1}(x)=g^{n}(x)\oplus g(g^{n}(x))$.
\item For every oracle $r\in \Cantor$, consider the following jump ideal defined as
\[\mathcal{J}_a(g,r)=\{z\in \Cantor:(\exists n\in\mathbb{N})\;x\leq_ag^n(r)\},\]
where $\leq_a$ denotes the arithmetical reducibility, that is, $p\leq_aq$ is defined by $p\leq_Tq^{(m)}$ for some $m\in\mathbb{N}$ (see \name{Odifreddi} \cite[Section XII.2 and Chapter XIII]{OdiBook1}).
\item A function $g$ is {\em almost arithmetical reducible to} a function $h$ (written as $g\leq_{aa}h$) if for any $r\in\Cantor$ there is $x\in\Cantor$ with $x\geq_Tr$ such that
\[\mathcal{J}_a(g,x)\subseteq\mathcal{J}_a(h,x).\]
\item Let $\mathcal{G}$ and $\mathcal{H}$ be countable sets of total functions.
We say that $\mathcal{G}$ is {\em $aa$-included in} $\mathcal{H}$ (written as $\mathcal{G}\subseteq_{aa}\mathcal{H}$) if for all $g\in\mathcal{G}$, there is $h\in\mathcal{H}$ such that $g\equiv_{aa}h$ (i.e., $g\leq_{aa}h$ and $h\leq_{aa}g$).
\end{enumerate}
\end{definition}

A function $g:\Cantor\to \Cantor$ is said to be {\em monotone} if $x\leq_Ty$ implies $g(x)\leq_Tg(y)$.
An {\em oracle $\mathbf{\Pi}^0_2$-singleton} is a total function $g:\Cantor\to \Cantor$ whose graph is $G_\delta$.
Clearly, every oracle $\mathbf{\Pi}^0_2$-singleton is Borel measurable, whereas there is no upper bound of Borel ranks of oracle $\mathbf{\Pi}^0_2$-singletons.
For instance, if $\alpha$ is a computable ordinal, then the $\alpha$-th Turing jump $j_\alpha(x)=x^{(\alpha)}$ is a monotone oracle $\mathbf{\Pi}^0_2$-singleton for every computable ordinal $\alpha$ (see \name{Odifreddi} \cite[Proposition XII.2.19]{OdiBook1}, \name{Sacks} \cite[Corollary II.4.3]{SacksBook}, and \name{Chong-Yu} \cite[Theorem 2.1.4]{ChYuBook}).
The following is the key lemma in our proof, which will be proved in Section \ref{sec:main-construction}.

\begin{lemma}[Realization Lemma]\label{thm:refle}
There is a map $\repspb{Rea}$ transforming each countable set of monotone oracle $\mathbf{\Pi}^0_2$-singletons into a Polish space such that
\[\repspb{Rea}(\mathcal{G})\leq^\mathfrak{T}_\sigma\repspb{Rea}(\mathcal{H})\;\Longrightarrow\;\mathcal{G}\subseteq_{aa}\mathcal{H}.\]
\end{lemma}

\subsection{Construction}\label{sec:main-construction}

We construct a Polish space whose co-spectrum codes almost arithmetical degrees contained in a given countable set $\mathcal{G}$ of oracle $\Pi^0_2$ singletons.
For notational simplicity, given $x\in[0,1]^\mathbb{N}$, we write $x_n$ for the $n$-th coordinate of $x$, and moreover, $x_{<n}$ and $x_{\leq n}$ for $(x_i)_{i<n}$ and $(x_i)_{i\leq n}$ respectively.

Our idea comes from the construction by \name{Miller} \cite[Lemma 9.2]{miller2}.
Our purpose is constructing a Polish space such that given $g\in\mathcal{G}$ and oracle $r$ the space has a point $x=(x_i)_{i\in\mathbb{N}}$ whose co-spectrum is not very different from $\mathcal{J}_a(g,r)$.
Then, at least, such a point should compute $g^i(r)$ for all $i\in\mathbb{N}$.
We can achieve this by requiring $x_i=g^i(r)$ for infinitely many $i\in\mathbb{N}$; however, we need to control the co-spectrum simultaneously, and therefore, we have to choose such coding locations $i$ very carefully.
The actual construction is that, from $r$ and $x_{<v}$, we will find a finite set $(\ell(u))_{u\leq v}$ of candidates of safe coding locations, and then we define $x_{\ell(u)}=g^{\ell(u)}(r)$ at a genuine safe coding location $\ell(u)$.
Then, for each $i$ with $v\leq i<\ell(u)$, we define $x_i$ from $(r,x_{<i},x_{\ell(u)})$ in a left-c.e.\ manner.
This idea yields the following definition.

\begin{definition}\label{def:reaspace}
Let $\mathcal{G}=(g_n)_{n\in\mathbb{N}}$ be a countable collection of oracle $\mathbf{\Pi}^0_2$-singletons.
The space $\repspb{\omega CEA}(\mathcal{G})$ consists of $(n,d,e,r,x)\in\mathbb{N}^3\times \Cantor\times[0,1]^\mathbb{N}$ such that for every $i$,
\begin{enumerate}
\item either $x_i=g_n^i(r)$, or
\item there are $u\leq v\leq i$ such that $x_i\in [0,1]$ is the $e$-th left-c.e.~real relative to $\langle r,x_{<i},x_{\ell(u)}\rangle$ and $x_{\ell(u)}=g_n^{\ell(u)}(r)$, where $\ell(u)=\Phi_d(u,r,x_{<v})\geq i$.
\end{enumerate}

%

Moreover, for a set $P\subseteq[0,1]^\mathbb{N}$, define $\repspb{\omega CEA}(\mathcal{G},P)$ to be the set of all points $(d,e,r,x)\in\repspb{\omega CEA}(\mathcal{G})$ with $(r,x)\in P$.
\end{definition}

\begin{lemma}
Suppose that $\mathcal{G}$ is an oracle $\mathbf{\Pi}^0_2$-singleton, and $P$ is a $\mathbf{\Pi}^0_2$ subset of $[0,1]^\mathbb{N}$.
Then, $\repspb{\omega CEA}(\mathcal{G},P)$ is Polish.
\end{lemma}

\begin{proof}
It suffices to show that $\repspb{\omega CEA}(\mathcal{G})$ is $\mathbf{\Pi}^0_2$.
The condition (1) in Definition \ref{def:reaspace} is clearly $\mathbf{\Pi}^0_2$.
Let $\forall a\exists b>a\;G(a,b,n,\ell,r,x)$ be a $\mathbf{\Pi}^0_2$ condition representing $x=g_n^\ell(r)$, and $\ell(u)[s]$ be the stage $s$ approximation of $\Phi_d(u,r,x_{<v})$.
The condition (2) is equivalent to the statement that there are $u\leq v\leq i$ such that
\begin{align*}
(\forall t\in\mathbb{N})(\exists s>t)\;&\ell(u)[s]\downarrow\geq i,\;d(x_i,J_{e,s}^{i+1}(r,x_{<i},x_{\ell(u)[s]}))<2^{-t},\\
&\mbox{ and }G(t,s,n,\ell(u)[s],r,x_{\ell(u)[s]}).
\end{align*}

Clearly, this condition is $\mathbf{\Pi}^0_2$.
\end{proof}

\begin{remark}
The space $\repspb{\omega CEA}(\mathcal{G})$ is totally disconnected for any countable set $\mathcal{G}$ of oracle $\mathbf{\Pi}^0_2$ singletons, since for any fixed $(n,d,e,r)\in\mathbb{N}^3\times \Cantor$, its extensions form a finite-branching infinite tree $T\subseteq[0,1]^{<\omega}$.
\end{remark}

Recall from Section \ref{subsec:proof-main-half} that \name{Miller} \cite[Lemma 6.2]{miller2} constructed a $\Pi^0_1$ set ${\rm Fix}(\Psi)$ such that ${\rm coSpec}(x)=\{x_i:i\in\mathbb{N}\}$ for every $x=(x_i)_{i\in\mathbb{N}}\in{\rm Fix}(\Psi)$.
By Lemma \ref{lem:Miller2}, without loss of generality, we may assume that ${\rm Fix}(\Psi)\cap\pi_0^{-1}\{r\}\not=\emptyset$ for every $r\in [0,1]$.
Now, consider the space $\repspb{Rea}(\mathcal{G})=\repspb{\omega CEA}(\mathcal{G},{\rm Fix}(\Psi))$.
To state properties of $\repspb{Rea}(\mathcal{G})$, for an oracle $\mathbf{\Pi}^0_2$-singleton $g$ and an oracle $r\in \Cantor$, we use the following Turing ideal:
\[\mathcal{J}_T(g,r)=\{z\in \Cantor:(\exists n\in\mathbb{N})\;x\leq_Tg^n(r)\}.\]

The following is the key lemma, which states that any collection of jump ideals generated by countably many oracle $\mathbf{\Pi}^0_2$-singletons has to be the degree co-spectrum of a Polish space up to the almost arithmetical equivalence!

\begin{lemma}\label{lem:keyref}Suppose that $\mathcal{G}=(g_n)_{n\in\mathbb{N}}$ is a countable set of oracle $\mathbf{\Pi}^0_2$-singletons.
\begin{enumerate}
\item For every $x\in\repspb{Rea}(\mathcal{G})$, there are $r\in \Cantor$ and $n\in\mathbb{N}$ such that
\[\mathcal{J}_T(g_n,r)\subseteq{\rm coSpec}(x)\subseteq\mathcal{J}_a(g_n,r).\]
\item For every $r\in \Cantor$ and $n\in\mathbb{N}$, there is $x\in\repspb{Rea}(\mathcal{G})$ such that
\[\mathcal{J}_T(g_n,r)\subseteq{\rm coSpec}(x).\]
\end{enumerate}
\end{lemma}

\begin{proof}[Proof of Lemma \ref{lem:keyref} (1)]
We have $(r,x)\in{\rm Fix}(\Psi)$ for every $y=(n,d,e,r,x)\in\repspb{Rea}(\mathcal{G})$.
For every $i\in\mathbb{N}$, we inductively assume that for every $j<i$, $x_j$ is arithmetical in $g_n^k(r)$ for some $k\in\mathbb{N}$.
Now, either $x_i=g_n^i(r)$ or $x_i$ is left-c.e.~in $(r,x_{<i},g_n^l(r))$ for some $l$.
In both cases, $x_i$ is arithmetical in $g_n^k(r)$ for some $k$.
Since $(r,x)\in{\rm Fix}(\Psi)$, by Lemma \ref{lemma:Miller-6-2}, ${\rm coSpec}(y)=\{r\}\cup\{x_i:i\in\mathbb{N}\}$.
This shows that ${\rm coSpec}(y)\subseteq\mathcal{J}_a(g_n,r)$.
Moreover, $x_i=g_n^i(r)$ for infinitely many $i\in\mathbb{N}$, since either $x_i=g^n_i(r)$ holds or there is $l\geq i$ such that $x_l=g^l_n(r)$ by the condition (2) in Definition \ref{def:reaspace}.
Therefore, $g_n^k(r)\leq_Mx$ for all $k\in\mathbb{N}$, that is, $\mathcal{J}_T(g_n,r)\subseteq{\rm coSpec}(y)$.
\end{proof}

To verify the assertion (2) in Lemma \ref{lem:keyref}, indeed, for any $n\in\mathbb{N}$, we will construct indices $d$ and $e$ such that for every $r\in \Cantor$, there is $x$ with $(n,d,e,r,x)\in\repspb{Rea}(\mathcal{G})$, where $x_i=g_n^i(r)$ for infinitely many $i\in\mathbb{N}$.
The $e$-th left-c.e.~procedure $J^{i+1}_e(r,x_{<i},x_{l(u)})$ is a simple procedure extending $r,x_{<i},x_{l(u)}$ to a fixed point of $\Psi$.
The function $\Phi_d$ searches for a {\em safe coding location} $c(n)$ from a given {\em name} of $x_{\leq c(n-1)}$, where $c(n-1)$ is the previous coding location.

To make sure the search of the next coding location is bounded, as in Definition \ref{def:reaspace}, we have to restrict the set of names of a $v$-tuple $x_{<v}$ to at most $v+1$ candidates.
It is known that a separable metrizable space is at most $n$-dimensional if and only if it is the union of $n+1$ many zero-dimensional subspaces (see \cite[Theorem 1.5.8]{EngBook} or \cite[Corollary 3.1.7]{vMBook2}).
We say that an admissibly represented Polish space is {\em computably at most $n$-dimensional} if it is the union of $n+1$ many subspaces that are computably homeomorphic to subspaces of $\mathbb{N}^\mathbb{N}$.

\begin{lemma}\label{lem:compu-n-dim}
Suppose that $(\repsp{X},\rho_X)$ is a computably at most $n$-dimensional admissibly represented space.
Then, there is a partial computable injection $\nu_X:\subseteq (n+1)\times\repsp{X}\to\mathbb{N}^\mathbb{N}$ such that for every $x\in\repsp{X}$, there is $k\leq n$ such that $(k,x)\in{\rm dom}(\nu_X)$ and $\rho_X\circ\nu_X(k,x)=x$.
\end{lemma}

\begin{proof}
By definition, $\repsp{X}$ is divided into $n+1$ many subspaces $S_0,\dots,S_n$ such that $S_k$ is homeomorphic to $N_k\subseteq\mathbb{N}^\mathbb{N}$ via computable maps $\tau_k$ and $\tau_k^{-1}$.
Then, the partial computable injection $\tau_k^{-1}:\subseteq\mathbb{N}^\mathbb{N}\to\repsp{X}$ has a computable realizer $\tau^*_k$, i.e., $\tau_k^{-1}=\rho_X\circ\tau^*_k$.
Define $\nu_X(k,x)=\tau^*_k\circ\tau_k(x)$ for $x\in S_k$.
Then, we have $\rho_X\circ\nu_X(k,x)=\tau_k^{-1}\circ\tau_k(x)=x$ for $x\in S_k$.
\end{proof}

The Euclidean $n$-space $\mathbb{R}^n$ is clearly computably $n$-dimensional, e.g., let $S_k$ be the set of all points $x\in\mathbb{R}^n$ such that exactly $k$ many coordinates are irrationals.
Furthermore, one can effectively find an index of $\nu_n:=\nu_{\mathbb{R}^n}$ in Lemma \ref{lem:compu-n-dim} uniformly in $n$.
Hereafter, let $\rho_i$ be the usual Euclidean admissible representation of $\mathbb{R}^i$.
Now, a coding location $c(n)$ will be obtained as a fixed point in the sense of Kleene's recursion theorem (Fact \ref{fact:kleene-recursion}).
Hence, one can effectively find such a location in the following sense.

\begin{lemma}[\name{Miller} {\cite[Lemma 9.2]{miller2}}]
Suppose that $(r,x_{<i})$ can be extended to a fixed point of $\Psi$, and fix a partial computable function $\nu$ which sends $x_{<i}$ to its name, i.e.,  $\rho_i\circ\nu(x_{<i})=(x_{<i})$.
From an index $t$ of $\nu$ and the sequence $x_{<i}$, one can effectively find a location $p=\Gamma(t,r,x_{<i})$ such that for every real $y$, the sequence $(r,x_{<i})$ can be extended to a fixed point $(r,x)$ of $\Psi$ such that $x_p=y$.
\end{lemma}

Let $t(n,k)$ be an index of the partial computable function $x\mapsto\nu_n(k,x)$.
We define $\Phi_d(u,r,x_{<v})$ to be $\Gamma(t(v,u),r,x_{<v})$ for every $u\leq v$.
Note that indices $d$ and $e$ do not depend on $g_n$.

\begin{proof}[Proof of Lemma \ref{lem:keyref} (2)]
Now, we claim that for every $r\in \Cantor$ and $n\in\mathbb{N}$, there is $x$ with $(n,d,e,r,x)\in\repspb{Rea}(\mathcal{G})$, where $x_i=g_n^i(r)$ for infinitely many $i\in\mathbb{N}$.
We follow the argument by \name{Miller} \cite[Lemma 9.2]{miller2}.
Suppose that $i$ is a coding location of $g_n^i(r)$, and $(r,x_{\leq i})$ is extendible to a fixed point of $\Psi$.
Then, there is $k\leq i+1$ such that $p=\Phi_d(k,r,x_{\leq i})$ is defined, and then we set $x_p=g_n^p(r)$.
By the property of $\Phi_d$, $(r,x_{\leq i},x_p)$ can be extended to a fixed point of $\Psi$.
Then, the $e$-th left-c.e.~procedure automatically produces $x_{\leq p}$ which is extendible to a fixed point of $\Psi$.
Note that the condition (2) in Definition \ref{def:reaspace} is ensured via $u=k$, $v=i+1$, and $l(u)=p$.
Eventually, we obtain $(r,x)\in{\rm Fix}(\Psi)$ such that $z=(n,d,e,r,x)\in\repspb{Rea}(\mathcal{G})$.

Clearly, $g_n^k(r)\in{\rm coSpec}(z)$ for every $k\in\mathbb{N}$, since ${\rm coSpec}(z)$ is a Turing ideal, and $g_n^k(r)\leq_Tg_n^{k+1}(r)$.
Consequently, $\mathcal{J}_T(g_n,r)\subseteq{\rm coSpec}(z)$.
\end{proof}

\begin{proof}[Proof of Lemma \ref{thm:refle}]
Suppose that $\repspb{Rea}(\mathcal{G})\leq^\mathfrak{T}_\sigma\repspb{Rea}(\mathcal{H})$.
Then, $\mathbb{N}\times\repspb{Rea}(\mathcal{G})\leq^\mathfrak{T}_\sigma\mathbb{N}\times\repspb{Rea}(\mathcal{H})$, and by Theorem \ref{lemma:spectrum-main} and Observation \ref{obs:main-cospec-inv}, the degree cospectrum of $\repspb{Rea}(\mathcal{G})$ is a sub-cospectrum of that of $\repspb{Rea}(\mathcal{H})$ up to an oracle $p$.
Fix enumerations $\mathcal{G}=(g_n)_{n\in\mathbb{N}}$ and $\mathcal{H}=(h_n)_{n\in\mathbb{N}}$.

\begin{claim}
For any $n$ and $u$, there are $m$ and $v$ such that $\mathcal{J}_a(g_n,u)=\mathcal{J}_a(h_m,v)$.
\end{claim}

By Lemma \ref{lem:keyref} (2), for any $n$ and $u\geq_Tp$, there is $x\in\repspb{Rea}(\mathcal{G})$ such that $\mathcal{J}_T(g_n,u)\subseteq{\rm coSpec}(x)\subseteq\mathcal{J}_a(g_n,u)$.
Then, there is $y\in\repspb{Rea}(\mathcal{H})$ such that ${\rm coSpec}^p(x)={\rm coSpec}^p(y)$.
We may assume that $p\leq_My$, otherwise $(y,p)$ has Turing degree by Lemma \ref{lem:JM-almosttotal}.
By Lemma \ref{lem:keyref} (1), there exist $m$ and $v$ such that $\mathcal{J}_T(h_m,v)\subseteq{\rm coSpec}(y)\subseteq\mathcal{J}_a(h_m,v)$.
Now, ${\rm coSpec}(x)={\rm coSpec}(y)$ holds, and note that $\mathcal{J}_T(h_m,v)\subseteq\mathcal{J}_a(g_n,u)$ implies $\mathcal{J}_a(h_m,v)\subseteq\mathcal{J}_a(g_n,u)$.
This verifies the claim.

For a fixed $n$, $\beta_n(u)$ chooses $m$ fulfilling the above claim for some $v$.
It is not hard to see that there is $m(n)$ such that $\beta_n(u)=m(n)$ for cofinally many $u$.
Then, for cofinally many $v$, there is $u$ such that $\mathcal{J}_a(g_n,u\oplus v)=\mathcal{J}_a(h_{m(n)},u\oplus v)$ by monotonicity.
Therefore, $g_n\equiv_{aa}h_{m(n)}$.
Consequently, $\mathcal{G}\subseteq_{aa}\mathcal{H}$.
\end{proof}

\begin{proof}[Proof of Theorem \ref{thm:maintheorem_emb}]
Let $S$ be a countable subset of $\omega_1$.
Note that $\sup S$ is countable by regularity of $\omega_1$.
Then, there is an oracle $p$ such that $\sup S<\omega_1^{{\rm CK},p}$, where $\omega_1^{{\rm CK},p}$ is the smallest noncomputable ordinal relative to $p$.
Now, the $\alpha$-th Turing jump operator $j^p_\alpha$ for $\alpha<\omega_1^{{\rm CK},p}$ is defined via a $p$-computable coding of $\alpha$.
By Spector's uniqueness theorem (see \name{Sacks} \cite[Corollary II.4.6]{SacksBook} or \name{Chong-Yu} \cite[Section 2.3]{ChYuBook}), the Turing degree of $j^p_\alpha(x)$ for $x\geq_Tp$ is independent of the choice of coding of $\alpha$, and so is $\mathcal{J}_a(j^p_\alpha,x)$.
Therefore, we simply write $j_\alpha$ for $j^p_\alpha$.

Define $\mathcal{G}_S=\{j_{\omega^{1+\alpha}}:\alpha\in S\}$.
We show that $S\subseteq T$ if and only if $\mathcal{G}_S\subseteq_{aa}\mathcal{G}_T$.
Suppose $\alpha\not=\beta$, say $\alpha<\beta$.
Clearly, $j_{\omega^{\alpha}}\leq_{aa}j_{\omega^\beta}$.
Suppose for the sake of contradiction that $j_{\omega^{\beta}}\leq_{aa}j_{\omega^\alpha}$.
Then, in particular, for every $x\leq_a\emptyset^{(\omega^\beta\cdot t)}$ with $t\in\mathbb{N}$, we must have $\emptyset^{(\omega^\beta\cdot(t+1))}\leq_ax^{(\omega^{\alpha}\cdot m)}$ for some $m\in\mathbb{N}$.
Thus, there is $n$ such that $\emptyset^{(\omega^\beta\cdot t+\omega^\beta)}\leq_T\emptyset^{(\omega^\beta\cdot t+\omega^{\alpha}\cdot m+n)}<_T\emptyset^{(\omega^\beta\cdot t+\omega^{\alpha+1})}$.
This is a contradiction.

Now, given countable sets $S,T\subseteq\omega_1$, if $S\subseteq T$, then $\repsp{Rea}(\mathcal{G}_S)$ clearly embeds into $\repsp{Rea}(\mathcal{G}_T)$.
If $S\not\subseteq T$, then the above argument shows that $\mathcal{G}_S\not\subseteq_{aa}\mathcal{G}_T$, and therefore, by Lemma \ref{thm:refle}, we have $\repsp{Rea}(\mathcal{G}_S)\not\leq_\sigma^\mathfrak{T}\repsp{Rea}(\mathcal{G}_T)$.
Consequently, $S\mapsto\gamma\repsp{Rea}(\mathcal{G}_S)$ is an order-preserving embedding of $([\omega_1]^{\leq\omega},\subseteq)$ into the $\sigma$-embeddability order $\leq_\sigma^\mathfrak{T}$ on compact metrizable spaces, where $\gamma\repsp{X}$ is \name{Lelek}'s compactification of $\repsp{X}$ in Fact \ref{thm:equivalence}.
\end{proof}

\begin{corollary}\label{cor:main-cor}
There exists a collection $(\repsp{X}_\alpha)_{\alpha<2^{\aleph_0}}$ of continuum many compact metrizable spaces satisfying the following conditions:
\begin{enumerate}
\item If $\alpha\not=\beta$, then $\repsp{X}_\alpha$ does not finite level Borel embed into $\repsp{X}_\beta$.
\item If $\alpha\not=\beta$, then the Banach algebra $\mathcal{B}^*_n(\repsp{X}_\alpha)$ is not linearly isometric (not ring isomorphic etc.) to $\mathcal{B}^*_n(\repsp{X}_\beta)$ for all $n\in\mathbb{N}$.
\end{enumerate}
\end{corollary}

\begin{proof}
By Theorems \ref{theo:spectrum-main} and \ref{thm:maintheorem_emb}.
Here, we note that if $\repsp{X}$ is $n$-th level Borel isomorphic to $\repsp{Y}$, then $\mathbb{N}\times\repsp{X}$ is again $n$-th level Borel isomorphic to $\mathbb{N}\times\repsp{Y}$.
\end{proof}

\section{Infinite Dimensional Topology}\label{sec:intermediate-dimension}

\subsection{Pol's Compactum}

In this section, we will shed light on dimension-theoretic perspectives of the $\omega$-left-CEA space.
Note that $\repspb{\omega CEA}$ is a totally disconnected infinite dimensional space.
We first compare our space $\repspb{\omega CEA}$ and a totally disconnected infinite dimensional space $\repspb{RSW}$ which is constructed by \name{Rubin}, \name{Schori}, and \name{Walsh} \cite{schori}.
A {\em continuum} is a connected compact metric space, and a continuum is {\em nondegenerated} if it contains at least two points.

It is known that the hyperspace ${\sf CK}(\repsp{X})$ of continua in a compact metrizable space $\repsp{X}$ equipped with the Vietoris topology forms a Polish space.
Hence, we may think of ${\sf CK}(\repsp{X})$ as a represented space, which corresponds to a positive-and-negative representation of the hyperspace in computable analysis.
We consider a closed subspace $S$ of ${\sf CK}([0,1]^\mathbb{N})$ consisting of all continua connecting opposite faces $\pi_0^{-1}\{0\}$ and $\pi_0^{-1}\{1\}$.
Then, fix a total Cantor representation of $S$, i.e., a continuous surjection $\delta_{\sf CK}$ from the Cantor set $C\subseteq[0,1]$ onto $S$.
We define the {\em Rubin-Schori-Walsh space} $\repspb{RSW}$ \cite{schori} (see also \cite[Theorem 3.9.3]{vMBook2}) as follows:
\begin{align*}
\repspb{RSW}&=\{\min (\delta_{\sf CK}(p)^{[p]}):p\in C\},\\
&=\{\min A^{[p]}:\mbox{$A$ is the $p$-th continuum of $[0,1]^\mathbb{N}$ with $[0,1]\subseteq\pi_0[A]$}\},
\end{align*}
where $A^{[p]}=A\cap\pi_0^{-1}\{p\}=\{z\in A:\pi_0(z)=p\}$, and recall that $\min P$ is the leftmost point of $P$ defined in Section \ref{subsec:proof-main-half}.
For notational convenience, without loss of generality, we may assume that the $e$-th $z$-computable continuum is equal to the $\langle e,z\rangle$-th continuum, where recall that $\langle\cdot,\cdot\rangle$ is a pairing function on natural numbers.

A compactification of $\repspb{RSW}$ is well-known in the context of \name{Alexandrov}'s old problem in dimension theory.
\name{Pol}'s compactum $\repspb{RP}$ is given as a compactification in the sense of \name{Lelek} of the space $\repspb{RSW}$.
Hence, we can see that $\repspb{RP}$ and $\repspb{RSW}$ have the same point degree spectra (modulo an oracle) as in the proof of Fact \ref{thm:equivalence}.
Surprisingly, these spaces have the same degree spectra as the space $\omega\repspb{CEA}$ up to an oracle.

\begin{theorem}\label{thm:degspec-dim}
The following spaces are all $\sigma$-homeomorphic to each other.
\begin{enumerate}
\item The $\omega$-left-CEA space $\repspb{\omega CEA}$.
\item Rubin-Schori-Walsh's totally disconnected strongly infinite dimensional space $\repspb{RSW}$.
\item Roman Pol's counterexample $\repspb{RP}$ to Alexandrov's problem.
\end{enumerate}
\end{theorem}

As a corollary, Roman Pol's compactum is second-level Borel isomorphic to the $\omega$-left-CEA space $\repspb{\omega CEA}$.
To prove Theorem \ref{thm:degspec-dim}, we show two lemmata.

\begin{lemma}\label{lem:degspec-dim1}
Every point of $\repspb{RSW}$ is $\omega$-left-CEA.
\end{lemma}

\begin{proof}
By Lemma \ref{lem:basis}, $\min A^{[p]}$ is $\omega$-left-CEA in $p$, since $A^{[p]}$ is $\Pi^0_1(p)$.
Moreover, clearly, $p\leq_M\min A^{[p]}$.
Thus, $\min A^{[p]}$ is $\omega$-left-CEA.
\end{proof}

For $\omega\repsp{CEA}\leq_\sigma^\mathfrak{T}\repsp{RSW}$, we need to show that every $\omega$-left-CEA point is realized as a leftmost point of a computable continuum in a uniform manner.
Indeed, we will show the following.

\begin{lemma}\label{lem:degspec-dim2}
Suppose that $x\in[0,1]^\mathbb{N}$ is $\omega$-left-CEA in a point $z\in \Cantor$.
Then, there is a nondegenerated $z$-computable continuum $A\subseteq[0,1]^\mathbb{N}$ such that $[0,1]\subseteq\pi_0[A]$ and $\min A^{[p]}=(p,x)$ for a name $p$ of $A$.
\end{lemma}

\begin{proof}
Given $p$, we will effectively construct a name $\Psi(p)$ of a continuum $A$.
By Kleene's recursion theorem (Fact \ref{fact:kleene-recursion}), we may fix $p$ such that the $p$-th continuum is equal to the $\Psi(p)$-th continuum.

Fix an $\omega$-left-CEA operator $J$ generated by $\langle W_n\rangle_{n\in\mathbb{N}}$ such that $J(z)=x$.
Here, as in the proof of Proposition \ref{prop:universaloperator}, each $W_n$ is a c.e.~list of pairs $(i,p)$, which indicates that ``if a given $n$-tuple $(z_0,\dots,z_{n-1})$ is in the $i$-th ball $B_i^n\subseteq[0,1]^n$, then $J^n_{W_n}(z_0,\dots,z_{n-1})\geq p$.''
Since $p=\langle e,z\rangle$ for some $e\in\mathbb{N}$, we have a computable function $\pi$ with $\pi(p)=z$, and then, redefine $W_0$ to be $W_0\circ\pi$.
In this way, we may assume that $J(p)=x$.

At stage $0$, $\Psi$ constructs $A_0=[0,1]\times[0,1]^\mathbb{N}$.
At stage $s+1$, if we find some rational open ball $B_i^n\subseteq[0,1]^{n}$ and a rational $q\in\mathbb{Q}$ such that $W_{n,s}$ declares that ``if a given $n$-tuple $(z_0,\dots,z_{n-1})$ is in the $i$-th ball $B_i^n$, then $J^n_{W_n}(z_0,\dots,z_{n-1})\geq q$,'' by enumerating $(i,q)$, then $\Psi$ removes $\pi_0^{-1}[B(p;2^{-s})]\cap(B_i^n\times[0,q)\times[0,1]^\mathbb{N})$ from the previous continuum $A_{s-1}$, where $B(p;2^{-s})$ is the rational open ball with center $p$ and radius $2^{-s}$.

Now, we show $\min A^{[p]}=x:=(x_0,x_1,\dots)$.
Assume that $x_0,\dots,x_{n-1}$ is an initial segment of $\min A^{[p]}$.
We will show that $x_n=\pi_n(\min A^{[p]})=\min\pi_n[\{z\in A^{[p]}:(\forall i<n)\;\pi_i(z)=x_i\}]$.
Since $J^n_{W_n}(p,x_0,\dots,x_{n-1})=x_n$, $W_n$ declares this fact at some point, that is, for any rational $q<x_n$, there is $i$ such that $(i,q)\in W_n$ and $(p,x_0,\dots,x_{n-1})\in B^n_i$.
Therefore, $A\cap(\pi_0^{-1}[B(p;2^{-s})]\cap(B^n_i\times[0,q)\times[0,1]^\mathbb{N}))=\emptyset$.
Hence, if $y<x_n$, then no extension of $(p,x_0,\dots,x_{n-1},y)$ is contained in $A$.
Moreover, if $(p,x_0,\dots,x_{n-1})\in B^n_i$ and $q<x_n$, then $(i,q)\not\in W_n$.
Hence, $x_n=\pi_n(\min A^{[p]})$ as desired.

Now, clearly $\min A^{[p]}=(p,x)$.
Note that $\Psi$ defines a $z$-computable continuum $A$ in a uniform manner.
The computability is ensured because we only remove a subset of $\pi_0^{-1}[B(p;2^{-s})]$ after stage $s$.
For the connectivity, assume that $A\subseteq U\cup V$ for some open sets $U,V\subseteq[0,1]^\mathbb{N}$.
By compactness, one can assume that $U$ and $V$ mention only finitely many coordinates, that is, there is $n_0$ such that if $y=(y_n)_{n\in\mathbb{N}}\in U$ ($V$, resp.)\ then $(y_0,\dots,y_{n_0},z_{n_0+1},z_{n_0+2},\dots)\in U$ ($V$, resp.)\ for any $(z_{n_0+m})_{m\in\mathbb{N}}$.
Given $y=(y_n)_{n\in\mathbb{N}}$, define $y^\ast=(y_0,\dots,y_{n_0},\vec{1})$.
By our choice of $n_0$, and our definition of $A$, $y\in A\cap U$ implies $y^\ast\in A\cap U$.
By our construction of $A$, if $k\leq n_0$, then any $(y_0,\dots,y_k,\vec{1})\in A\cap U$ is connected to $(y_0,\dots,y_{k-1},1,\vec{1})\in A\cap U$ by a line segment inside $A\cap U$.
Therefore, for any point $y\in A\cap U$, $y^\ast$ is connected to $\vec{1}$ by a polygonal line inside $A\cap U$.
The same holds true for $V$.
Hence, if $A\cap U$ and $A\cap V$ are nonempty, $y\in A\cap U$ and $z\in A\cap V$ say, $y^\ast\in A\cap U$ and $z^\ast\in A\cap V$, and they are connected to $\overline{1}$, and therefore, there is a path from $y^\ast$ to $z^\ast$ in $A\cap(U\cup V)$.
By connectivity of the path, $A\cap U$ and $A\cap V$ has an intersection in the path.
This shows that $A$ cannot written as a union of disjoint open subsets.
Consequently, $A$ is connected.
\end{proof}

\begin{proof}[Proof of Theorem \ref{thm:degspec-dim}]
By Theorem \ref{theo:spectrum-main}, and Lemmata \ref{lem:degspec-dim1} and \ref{lem:degspec-dim2}.
\end{proof}

The properness of $\repspb{RSW} <_{\sigma}^\mathfrak{T} [0,1]^\mathbb{N}$ can also be obtained by some relatively recent work on infinite dimensional topology: the Hilbert cube (indeed, any strongly infinite dimensional compactum) is not $\sigma$-hereditary-disconnected (see \cite{Radul}).
However, such an argument does not go any farther for constructing second-level Borel isomorphism types, and indeed, to the best of our knowledge, no known topological technique provides us four or more second-level Borel isomorphism types.

On a side note, one can also define the graph $n\repspb{CEA}\subseteq\mathbb{N}\times\Cantor\times[0,1]^n$ of a universal $n$-left-CEA operator (as an analogy of an $n$-REA operator) in a straightforward manner. The space $n\repspb{CEA}$ is an example of a finite-dimensional Polish spaces whose infinite product has again the same dimension. The first such examples were constructed by \name{Kulesza} in \cite{kulesza}.

\begin{proposition}
The space $n\repspb{CEA}$ is a totally disconnected $n$-dimensional Polish space.
Moreover, the countable product $n\repspb{CEA}^\mathbb{N}$ is again $n$-dimensional.
\end{proposition}

\begin{proof}
Clearly, $n\repspb{CEA}$ is totally disconnected and Polish.
To check the $n$-dimensionality, we think of $n\repspb{CEA}$ as a subspace of $[0,1]^{n+1}$ by identifying $(e,x)\in\mathbb{N}\times\Cantor$ with $\iota(0^e1x)\in[0,1]$, where $\iota$ is a computable embedding of $\Cantor$ into $[0,1]$.
We claim that $n\repspb{CEA}$ intersects with all continua $A\subseteq [0,1]^{n+1}$ such that $[0,1]\subseteq\pi_0[A]$.
We have a computable function $d$ such that the $d(e)$-th $n$-left-CEA procedure $J^n_{d(e)}(x)$ for a given input $x\in\Cantor$ outputs the value $y\in[0,1]^n$ such that $(\iota(0^e1x),y)=\min A_{e,x}^{[\iota(0^e1x)]}$, where $A_{e,x}$ is the $e$-th $x$-computable continuum in $[0,1]^{n+1}$ such that $[0,1]\subseteq\pi_0[A_{e,x}]$.
By Kleene's recursion theorem (Fact \ref{fact:kleene-recursion}), there is $r$ such that $J^n_{d(r)}=J^n_r$.
Hence, $(\iota(0^r1x),J^n_r(x))\in n\repspb{CEA}\cap A_{e,x}$, which verifies the claim.
The claim implies that $n\repspb{CEA}$ is $n$-dimensional (see \name{van Mill} \cite[Corollary 3.7.5]{vMBook2}).

To verify the second assertion, consider the (computably) continuous map $g$ from the square $n\repspb{CEA}^2$ into $n\repspb{CEA}$ such that for two points $\mathbf{x}=(e,r,x_0,\dots,x_{n-1})$ and $\mathbf{y}=(d,s,y_0,\dots,y_{n-1})$ in $n\repspb{CEA}$,
\[g(\mathbf{x},\mathbf{y})=(\langle e,d\rangle,r\oplus s,(x_0+y_0)/2,\dots,(x_{n-1}+y_{n-1})/2).\]
It is not hard to verify that $g^{-1}$ is also (computably) continuous.
Hence, $n\repspb{CEA}^2$ is computably embedded into $n\repspb{CEA}$.
In particular, it is $n$-dimensional.
Then, we can conclude that $n\repspb{CEA}^\mathbb{N}$ is also $n$-dimensional (by the same argument as in \name{van Mill} \cite[Theorem 3.9.5]{vMBook2}).
\end{proof}

\subsection{Nondegenerated Continua and $\omega\repspb{CEA}$ Degrees}

We may extract computability-theoretic contents from the construction of Rubin-Schori-Walsh's strongly infinite-dimensional totally disconnected space $\repspb{RSW}$.
The standard proof of non-countable-dimensionality of $\repspb{RSW}$ (hence, the existence of a non-Turing degree in $\repspb{RSW}$) indeed implies the following computability theoretic result.

\begin{proposition}\label{thm:avoid-Turing}
There exists a nondegenerated continuum $A\subseteq[0,1]^\mathbb{N}$ in which no point has Turing degree.
\end{proposition}

\begin{proof}
Define $\repsp{H}_{\langle{i,j\rangle}}\subseteq[0,1]^\mathbb{N}$ to be the set of all points which can be identified with an element in $\Cantor$ via the witnesses $\Phi_i$ and $\Phi_j$ (as in the proof of Lemma \ref{lemma:spectrum-main}).
Then, $\bigcup_{n}\repsp{H}_{n}$ is the set of all points in $[0,1]^\mathbb{N}$ having Turing degrees.
Note that each $\repsp{H}_{n}$ is zero-dimensional since it is homeomorphic to a subspace of $\Cantor$.

Consider the hyperplane $P_n^i=[0,1]^n\times\{i\}\times[0,1]^\mathbb{N}$ for each $n\in\mathbb{N}$ and $i\in\{0,1\}$.
It is well known that $\{(P_n^0,P_n^1)\}_{n\in\mathbb{N}}$ is essential in $[0,1]^\mathbb{N}$.
Then, by using the dimension-theoretic fact (see \name{van Mill} \cite[Corollary 3.1.6]{vMBook2}), we can find a separator $L_{n}$ of $(P_{n+1}^0,P_{n+1}^1)$ in $[0,1]^\mathbb{N}$ such that $L_{n}\cap\repsp{H}_{n}=\emptyset$ since $\repsp{H}_{n}$ is zero-dimensional.

Put $L=\bigcap_nL_{n}$.
Then, $L$ contains no point having Turing degree, since $L\cap\repsp{H}_n=\emptyset$ for every $n\in\mathbb{N}$.
Moreover, $L$ contains a continuum $A$ from $P_0^0$ to $P_0^1$ (see \name{van Mill} \cite[Proposition 3.7.4]{vMBook2}).
\end{proof}


Recall that our infinite dimensional version of Kreisel's basis theorem (Lemma \ref{lem:basis}) says that every $\Pi^0_1$ subset $P$ of the Hilbert cube has a point of an $\omega$-left-CEA continuous degree.
Surprisingly, we do not need any effectivity assumption on $P$ to prove this if $P$ is a nontrivial connected compact set.

\begin{proposition}\label{thm:Pol-deg}
Every nondegenerated continuum $A\subseteq[0,1]^\mathbb{N}$ contains a point of an $\omega$-left-CEA continuous degree.
\end{proposition}

\begin{proof}
Note that there is $n\in\omega$ such that $P_n^{[0,p]}$ and $P_n^{[q,1]}$ with some rationals $p<q\in\mathbb{Q}$ intersect with $A$, since $A$ is nondegenerated, where $P_n^{[a,b]}=[0,1]^n\times[a,b]\times[0,1]^\mathbb{N}$.
Clearly, there is no separator $C$ of $P_n^{[0,p]}$ and $P_n^{[q,1]}$ with $C\cap A=\emptyset$ (i.e., the pair $(P_n^{[0,p]},P_n^{[q,1]})$ is essential in $A$), since $A$ is not zero-dimensional.
Therefore, the pair $(P_n^{p},P_n^{q})$ is essential in the compact subspace $A\cap P_n^{[p,q]}$.
Hence, $A\cap P_n^{[p,q]}$ contains a continuum $B$ intersecting with $P_n^{p}$ and $P_n^{q}$ (see van Mill \cite[Proposition 3.7.4]{vMBook2}).
Consider a computable homeomorphism $h:P_n^{[p,q]}\cong[0,1]^\mathbb{N}$ mapping $P_n^{p}$ and $P_n^{q}$ to $P_0^{0}=\pi_0^{-1}(0)$ and $P_0^{1}=\pi_1^{-1}(1)$, respectively.
Then $h[B]$ is a continuum intersecting with $\pi_0^{-1}(0)$ and $\pi_0^{-1}(0)$, and therefore $[0,1]\subseteq\pi_0[h[B]]$.
Let $s$ be a name of $h[B]$.
Then, by definition, $\min h[B]^{[s]}\in\repspb{RSW}$, which has an $\omega$-left-CEA continuous degree by Lemma \ref{lem:degspec-dim1}.
In particular, $h[B]$ contains a point of an $\omega$-left-CEA continuous degree, and so does $A$ since $h$ is a computable homeomorphism and $B\subseteq A$.
\end{proof}

As a corollary, we can see that every compactum $A\subseteq[0,1]^\mathbb{N}$ of positive dimension contains a point of an $\omega$-left-CEA continuous degree.
Our proof of Theorem \ref{thm:avoid-Turing} is essentially based on the fact that for any sequence of zero-dimensional spaces $\{X_i\}_{i\in\mathbb{N}}$, there exists a continuum avoiding all $X_i$'s.
Contrary to this fact, Theorem \ref{thm:Pol-deg} says that $\{X_i\}_{i\in\mathbb{N}}$ cannot be replaced with a sequence of totally disconnected spaces.
We say that a space is $\sigma$-totally-disconnected if it is a countable union of totally disconnected subspaces.
Note that the complement of a $\sigma$-totally-disconnected subset of the Hilbert cube is infinite dimensional.

\begin{corollary}
There exists a $\sigma$-totally-disconnected set $X\subseteq [0,1]^\mathbb{N}$ such that any compact subspace of the complement $Y=[0,1]^\mathbb{N}\setminus X$ is zero-dimensional.
\end{corollary}

\begin{proof}
Define $X_{\langle i,j\rangle}$ to be the set of all points which can be identified with an element in $\omega\repspb{CEA}$ via the witnesses $\Phi_i$ and $\Phi_j$.
Then, $X_{\langle i,j\rangle}$ is totally disconnected since it is homeomorphic to a subspace of $\omega\repspb{CEA}$.
Clearly, no point $Y=[0,1]^\mathbb{N}\setminus\bigcup_{i,j\in\mathbb{N}}X_{\langle i,j\rangle}$ has an $\omega$-left-CEA continuous degree.
Assume that $Z$ is a compact subspace of $Y$ of positive dimension.
Then $Z$ has a nondegenerated subcontinuum $A$.
However, by Theorem \ref{thm:Pol-deg}, $A$ contains a point of an $\omega$-left-CEA continuous degree.
\end{proof}


\subsection{Weakly Infinite Dimensional Cantor Manifolds}

Recall that a Pol-type Cantor manifold is a compact metrizable $C$-space which cannot be disconnected by a hereditarily weakly infinite dimensional compact subspaces.
By combining a known construction in infinite dimensional topology, we can slightly extend Theorem \ref{thm:maintheorem_emb} as follows.

\begin{proposition}\label{thm:maintheorem_emb-cm}
There exists a collection $(\repsp{X}_\alpha)_{\alpha<2^{\aleph_0}}$ of continuum many Pol-type Cantor manifolds satisfying the following conditions:
\begin{enumerate}
\item if $\alpha\not=\beta$, $\repsp{X}_\alpha$ does not $\sigma$-embed into $\repsp{X}_\beta$.
\item If $\alpha\not=\beta$, then $\repsp{X}_\alpha$ does not finite level Borel embed into $\repsp{X}_\beta$.
\item If $\alpha\not=\beta$, then the Banach algebra $\mathcal{B}^*_n(\repsp{X}_\alpha)$ is not linearly isometric (not ring isomorphic etc.) to $\mathcal{B}^*_n(\repsp{X}_\beta)$ for all $n\in\mathbb{N}$.
\end{enumerate}
\end{proposition}

\begin{lemma}\label{lem:Pol-manifold}
For any $\mathcal{G}$, there exists a Pol-type Cantor manifold $\repsp{Z}(\mathcal{G})$ such that $\omega\repspb{CEA}\oplus\repspb{Rea}(\mathcal{G})\equiv_\sigma^\mathfrak{T}\repsp{Z}(\mathcal{G})$.
\end{lemma}

\begin{proof}
Recall from Theorem \ref{thm:degspec-dim} that $\omega\repspb{REA}$ is $\sigma$-homeomorphic to a strongly infinite dimensional space $\repspb{RSW}$.
Let $\repsp{R}_0$ and $\repsp{R}_1$ be homeomorphic copies of $\repspb{RSW}$, and let $\repsp{X}$ be a compactification of $\repsp{R}_0\oplus \repsp{R}_1\oplus\repspb{Rea}(\mathcal{G})$ in the sense of \name{Lelek} (recall from Fact \ref{thm:equivalence}).
Then, $\repsp{X}$ is $\sigma$-homeomorphic to $\omega\repspb{CEA}\oplus\repspb{Rea}(\mathcal{G})$.

We follow the construction of \name{El\.{z}bieta Pol} \cite[Example 4.1]{EPol96}.
Now, $\repsp{R}_0$ has a hereditarily strongly infinite dimensional subspace $\repsp{Y}$ \cite{Rub80}.
Choose a point $p\in\repsp{Y}$ and a closed set $F\subseteq\repsp{Y}$ containing $p$ such that every separator between $p$ and ${\rm cl}_\repsp{X}F$ is strongly infinite dimensional as in \cite[Example 4.1 (A)]{EPol96}.

Define $\repsp{K}=\repsp{X}/{\rm cl}_\repsp{X}F$ as in \cite[Example 4.1 (A)]{EPol96}.
To see that $\repsp{K}$ is $\sigma$-homeomorphic to $\repsp{X}$, we note that ${\rm cl}_{\repsp{X}}F\cap(\repsp{R}_1\cup\repspb{Rea}(\mathcal{G}))=\emptyset$ since $\repsp{R}_0$, $\repsp{R}_1$ and $\repspb{Rea}(\mathcal{G})$ are separated in $\repsp{X}$.
Therefore, ${\rm cl}_{\repsp{X}}F$ is covered by the union of $\repsp{R}_0$ (which is homeomorphic to $\repsp{R}_1$) and a countable dimensional space.
Define $\repsp{Z}$ as a Pol-type Cantor manifold in \cite[Example 4.1 (C)]{EPol96}.
Then, $\repsp{Z}(\mathcal{G}):=\repsp{Z}$ is the union of a finite dimensional space and countably many copies of $\repsp{K}$.
Consequently, $\repsp{Z}(\mathcal{G})$ is $\sigma$-homeomorphic to $\repspb{Rea}(\mathcal{G})$.
\end{proof}

\begin{proof}[Proof of Proposition \ref{thm:maintheorem_emb-cm}]
Combine Theorem \ref{thm:maintheorem_emb}, Corollary \ref{cor:main-cor}, and Lemma \ref{lem:Pol-manifold}.
\end{proof}


\subsection{An Ordinal Valued $\sigma$-Topological Invariant}

In this section, we will try to extract a topological content from our construction in Section \ref{sec:piecewisehomeo}.
However, although we have constructed continuum many mutually different spaces, it is difficult to discern dimension-theoretic differences among these spaces.
For instance, all of our spaces have the same transfinite Steinke dimensions \cite{ArChPu00,Radul}, game dimensions \cite{FedOsi08}, and so on (see \name{Chatyrko} and \name{Hattori} \cite{ChaHat13} for the thorough treatment of the notion of various kinds of transfinite dimensions).

We now focus on an $\aleph_1$ chain of $\sigma$-homeomorphism types of Polish spaces:
\[\mathbb{R}^n<^\mathfrak{T}_\sigma\repspb{Rea}(\{j_1\})<^\mathfrak{T}_\sigma\repspb{Rea}(\{j_\omega\})<^\mathfrak{T}_\sigma\repspb{Rea}(\{j_{\omega^2}\})<^\mathfrak{T}_\sigma\repspb{Rea}(\{j_{\omega^3}\})<^\mathfrak{T}_\sigma\dots\]

Our key observation was that closure properties of Scott ideals reflect $\sigma$-homeomorphism types of Polish spaces.
The purpose here is to provide a topological understanding of our method.

\subsubsection{Universal Dimension}

We first generalize an idea of understanding the topological dimension in the context of a universal function.
A continuous function $f:\repsp{X}\to\repsp{Y}$ is {\em universal} if for any continuous function $g:\repsp{X}\to\repsp{Y}$, there is $x\in\repsp{X}$ such that $f(x)=g(x)$.
It is known that a separable metrizable space $\repsp{X}$ is $n$-dimensional if and only if there is a universal function $f:\repsp{X}\to[0,1]^n$, and $\repsp{X}$ is strongly infinite dimensional if and only if there is a universal function $f:\repsp{X}\to[0,1]^\mathbb{N}$.
Then, how do we introduce an intermediate notion?
Our answer is thinking about partial {\em countably universal} functions and Baire functions of infinite rank.

Let $\repsp{X}$ and $\repsp{Y}$ be topological spaces, and assume that $\mathcal{G}$ is a collection of partial functions of type $\mathbb{N} \times \repsp{X} \to \repsp{Y}$.
A partial continuous function $f:\subseteq\mathbb{N}\times\repsp{X}\to\repsp{Y}$ is {\em countably $\mathcal{G}$-universal} if for any $g\in \mathcal{G}$,
\[(\exists x\in\repsp{X})(\forall n\in\mathbb{N})[(n,x)\in{\rm dom}(g)\;\Longrightarrow\;(\exists m\in\mathbb{N})\;f(m,x)=g(n,x)].\]

Our definition of countable universality ensures that the notion behaves well w.r.t.~$\sigma$-homeomorphisms. To express this, we now allow the spaces $\repsp{X}$ and $\repsp{Y}$ in the type of $\mathcal{G}$ to vary, and demand that $\mathcal{G}$ is closed under composition with continuous partial functions from both sides, as well as under products with the identity. Then we find that:

\begin{proposition}
\label{prop:countablyuniversalfunctions}
If $\repsp{X}$ $\sigma$-embeds into $\repsp{X}'$, then
\begin{enumerate}
\item If there exists a countably $\mathcal{G}$-universal partial function $f : \subseteq\mathbb{N}\times\repsp{X}\to\repsp{Y}$, then there exists a countably $\mathcal{G}$-universal partial function $f' : \subseteq\mathbb{N}\times\repsp{X}'\to\repsp{Y}$.
\item If there exists a countably $\mathcal{G}$-universal partial function $f' : \subseteq\mathbb{N}\times\repsp{Z}\to\repsp{X}'$, then there exists a countably $\mathcal{G}$-universal partial function $f : \subseteq\mathbb{N}\times\repsp{Z}\to\repsp{X}$.
\end{enumerate}
\begin{proof}
\begin{enumerate}
\item Let the $\sigma$-embedding of $\repsp{X}$ into $\repsp{X}'$ be witnessed by partial continuous embeddings $(\iota_n)_{n \in \mathbb{N}}$ with disjoint domains. Now we define $f' : \subseteq\mathbb{N}\times\repsp{X}'\to\repsp{Y}$ via $f'(\langle n, m\rangle, x) = f(m, \iota_n^{-1}(x))$, and obtain again a partial continuous function. Let $g' : \subseteq \mathbb{N} \times \repsp{X}' \to \repsp{Y}$ be in $\mathcal{G}$. Then $g : \subseteq \mathbb{N} \times \repsp{X} \to \repsp{Y}$ defined via $g(\langle n, m\rangle, x) = g'(m,\iota_n(x))$ is the composition of $g$ and a partial continuous function, hence $g \in \mathcal{G}$. Thus, by assumption there exists some $x_0 \in \mathbf{X}$ such that for all $\langle n, m\rangle \in \mathbb{N}$ such that $g(\langle n, m\rangle, x_0)$ is defined there exists some $l \in \mathbb{N}$ with $g(\langle n, m\rangle, x_0) = f(l,x_0)$.

    As the $(\iota_n)_{n \in \mathbb{N}}$ witness a $\sigma$-embedding, there exists a unique $n_0$ with $x_0 \in \dom(\iota_{n_0})$. We claim that the point $\iota_{n_0}(x_0)$ witnesses that $f'$ is countably $\mathcal{G}$-universal given $g'$. For any $m \in \mathbb{N}$ such that $g'(m,\iota_{n_0}(x_0))$ is defined, we find that $g'(m,\iota_{n_0}(x_0)) = g(\langle n_0, m\rangle, x_0)$. By assumption there exists some $l \in \mathbb{N}$ with $g(\langle n_0, m\rangle, x_0) = f(l,x_0)$. Moreover, $f(l,x_0) = f'(\langle n_0, l\rangle, x_0)$. Hence, $\langle n_0, l\rangle$ works as the witness for $m \in \mathbb{N}$.

\item Very similar to (1.): Let the $\sigma$-embedding of $\repsp{X}$ into $\repsp{X}'$ be witnessed by partial continuous embeddings $(\iota_n)_{n \in \mathbb{N}}$ with disjoint domains. We define $f : \subseteq\mathbb{N}\times\repsp{Z} \to \repsp{X}$ via $f(\langle n, m\rangle, x) = \iota_n^{-1}(f'(m, x))$, and obtain again a partial continuous function. Let $g : \subseteq \mathbb{N} \times \repsp{Z} \to \repsp{X}$ be in $\mathcal{G}$.  Then $g' : \subseteq \mathbb{N} \times \repsp{Z} \to \repsp{X}'$ defined via $g'(\langle n, m\rangle, x) = \iota_n(g'(m,x))$ is the composition of the product of the identity on $\mathbb{N}$ and $g$, and a partial continuous function, hence $g' \in \mathcal{G}$. Thus, by assumption there exists some $x_0 \in \mathbf{Z}$ such that for all $\langle n, m\rangle \in \mathbb{N}$ such that $g'(\langle n, m\rangle, x_0)$ is defined there exists some $l \in \mathbb{N}$ with $g'(\langle n, m\rangle, x_0) = f'(l,x_0)$.

    By definition $g'(\langle n, m\rangle, x_0) = \iota_n(g(m,x_0))$, so as the $(\iota_n)_{n \in \mathbb{N}}$ witness a $\sigma$-embedding, whenever $g(m,x_0)$ is defined, there exists a unique $n_m$ such that $g'(\langle n_m, m\rangle, x_0)$ is defined. By construction, we have that $\iota_{n_m}(g(m,x_0)) = f'(l, x_0)$. As $\iota_{n_m}$ is an embedding, it follows that $g(m,x_0) = \iota_{n_m}^{-1}(f'(l,x_0)) = f(\langle n_m, l\rangle, x_0)$ by definition of $f$.
\end{enumerate}
\end{proof}
\end{proposition}

Let $j_\alpha : \Cantor \to \Cantor$ be the $\alpha$-th Turing jump, and let $\mathcal{B}_\alpha(\repsp{X})$ be the set of functions $f :\subseteq \mathbb{N} \times \repsp{X} \to \Cantor$ such that there are continuous functions $H: \subseteq \mathbb{N} \times \repsp{X} \to \Cantor$ and $K : \subseteq \Cantor \to \Cantor$ such that $f(n,x) = K(j_\alpha(H(n,x)))$.



\begin{definition}
The {\em universal dimension} ${\rm udim}(\repsp{X})$ of $\repsp{X}$ is the supremum of countable ordinals $\alpha<\omega_1$ such that there is a countably $\mathcal{B}_\beta$-universal function $f:\subseteq\mathbb{N}\times\repsp{X}\to\Cantor$ for any $\beta<\omega^\alpha$, where, $\omega^0=1$.
If such $\alpha$ does not exist, then ${\rm udim}(\repsp{X})=-1$.
 \end{definition}

 It follows from Proposition \ref{prop:countablyuniversalfunctions} (1) that the universal dimension is invariant under $\sigma$-homeomorphism:

\begin{corollary}
If $\repsp{X}$ $\sigma$-embeds into $\repsp{Y}$, then ${\rm udim}(\repsp{X})\leq{\rm udim}(\repsp{Y})$.
\end{corollary}

\begin{proposition}\label{prop:udim}
If $\repsp{X}$ is a countable dimensional Polish space, then ${\rm udim}(\repsp{X})\leq 0$.
For every countable ordinal $\alpha>0$, there is a Pol-type Cantor manifold $\repsp{X}$ such that ${\rm udim}(\repsp{X})=\alpha$.
 \end{proposition}

 To show Proposition \ref{prop:udim}, we need the following effective interpretation of universal dimension.

\begin{lemma}\label{lem:jdimequiv1}
For any admissibly represented space $\repsp{X}$, the product $\mathbb{N}\times\repsp{X}$ admits a countably $\mathcal{B}_\alpha$-universal function if and only if relative to some oracle $r$, for all $z\in\Cantor$ there is a point $x\in\repsp{X}$ such that $z\in{\rm coSpec}^r(x)$ and ${\rm coSpec}^r(x)$ is closed under the $\alpha$-th Turing jump.
\end{lemma}

\begin{proof}
We will let $\Phi_e^r$ denote the $e$-th partial $r$-computable function from $\repsp{X}$ into $\Cantor$. We point out that if the $e$-th Turing functional outputs distinct elements of Cantor space on oracle $r$ and distinct names of $x$ as input, then $x \notin \dom \Phi_e^r(x)$

Suppose that $f:\subseteq\mathbb{N}\times\repsp{X}\to\Cantor$ is countably $\mathcal{B}_\alpha$-universal.
Since $\alpha$ is countable, and $f$ is continuous, some oracle $r$ satisfies that $\alpha<\omega_1^{{\rm CK},r}$, and that $f$ is computable relative to $r$.
Then, clearly $\{f(n,x):n\in\mathbb{N},\;(n,x)\in{\rm dom}(f)\}\subseteq{\rm coSpec}^r(x)$ holds for all $x\in\repsp{X}$.

Fix $z$.
Define $g(e,x)=z\oplus (j_\alpha\circ\Phi_e^r(x))$. Then $g \in \mathcal{B}_\alpha$. Let $x\in\repsp{X}$ be a point witnessing universality of $f$ for given $g$, that is, for any $e\in\mathbb{N}$, there is $n\in\mathbb{N}$ such that $f(n,x)=g(e,x)$.
Since $z\leq_Tg(e,x)$ for any $e$ and $x$ where $g(e,x)$ is defined, this implies that $z\leq_Tf(n,x)\in{\rm coSpec}^r(x)$.
Thus, $z\in {\rm coSpec}^r(x)$ since the $r$-cospecturm is closed under Turing reducibility.

Now, for any $p\in{\rm coSpec}^r(x)$, there is an index $e$ such that $p=\Phi_e^r(x)$.
By universality of $f$, there is $n$ such that $f(n,x)=g(e,x)=z\oplus p^{(\alpha)}$.
This shows that $p^{(\alpha)}\in{\rm coSpec}^r(x)$ since $p^{(\alpha)}\leq_T z\oplus p^{(\alpha)}=f(n,x)\in{\rm coSpec}^r(x)$.
Consequently, ${\rm coSpec}^r(x)$ is closed under the $\alpha$-th Turing jump.

Conversely, suppose that the condition in Lemma \ref{lem:jdimequiv1} holds for $r$.
We define $f(e,x)=\Phi_e^r(x)$ and claim that $f$ is universal.

For any $g \in \mathcal{B}_\alpha$ there exists some oracle $z \geq_T r$ such that the witnesses $H$, $K$ for membership of $g$ in $\mathcal{B}_\alpha$ can be chosen as computable relative to $z$. Clearly, $H(n,x)\oplus z\in{\rm coSpec}^z(x)$ and $K(j_\alpha(H(n,x)) \leq_T(H(n,x)\oplus z)^{(\alpha)}$.

By our assumption, there is a point $x\in\repsp{X}$ with $z\in{\rm coSpec}^r(x)$ such that the $r$-cospectrum of $x$ is closed under the $\alpha$-th Turing jump.
Note that $r\leq_T z\in{\rm coSpec}^r(x)$ implies ${\rm coSpec}^r(x)={\rm coSpec}^z(x)$, and therefore, the $z$-cospectrum of $x$ is also closed under the $\alpha$-th Turing jump.
This implies that $(H(n,x)\oplus z)^{(\alpha)}\in{\rm coSpec}^z(x)$.
Hence, we have $g(n,x) \in{\rm coSpec}^r(x)$ since the $r$-cospectrum is closed under Turing reducibility.
This means that $\Phi_e^r(x)=g(n,x)$ for some $e$, that is, $f(e,x)=g(n,x)$.
\end{proof}

\begin{proof}[Proof of Proposition \ref{prop:udim}]
If $\repsp{X}$ is countable dimensional, by Theorem \ref{theo:spectrum-main}, there is $r$ such that any point $x\in\repsp{X}$ has a Turing degree relative to $r$.
Therefore, ${\rm coSpec}^r(x)$ is not closed under the Turing jump.
By Lemma \ref{lem:jdimequiv1}, $\mathbb{N}\times\repsp{X}$ does not admit $\mathcal{B}_1$-universal function.
Consequently, ${\rm udim}(\repsp{X})\leq 0$.

For the second assertion, we first see that the jump-dimension of $\repspb{Rea}(j_{\omega^{\alpha}})$ is $\alpha+1$.
We have ${\rm jdim}(\repspb{Rea}(j_{\omega^{\alpha}}))\geq\alpha+1$ because for any $z$, there is $x\in\repspb{Rea}(j_{\omega^{\alpha}})$ such that $\mathcal{J}_T(j_{\omega^\alpha},z)\subseteq{\rm coSpec}(x)\subseteq\mathcal{J}_a(j_{\omega^\alpha},z)$ by Lemma \ref{lem:keyref}.
If $y\in{\rm coSpec}(x)$, then $y\in \mathcal{J}_a(j_{\omega^\alpha},z)$.
Therefore, $y^{(\omega^\alpha+n)}\in \mathcal{J}_T(j_{\omega^\alpha},z)\subseteq{\rm coSpec}(x)$ for all $n\in\mathbb{N}$.
Hence, ${\rm coSpec}^z(x)={\rm coSpec}(x)$ is closed under the $\beta$-th Turing jump for all $\beta<\omega^{\alpha+1}$.
To see ${\rm jdim}(\repspb{Rea}(j_{\omega^{\alpha}}))<\alpha+2$, we note for any $x\in\repspb{Rea}(j_{\omega^{\alpha}})$ that $\mathcal{J}_T(j_{\omega^\alpha},z)\subseteq{\rm coSpec}(x)\subseteq\mathcal{J}_a(j_{\omega^\alpha},z)$ for some $z$ by Lemma \ref{lem:keyref}.
Then, $z\in{\rm coSpec}(x)$, but ${\rm coSpec}(x)$ is covered by the Turing ideal generated by $z^{(\omega^{\alpha+1})}$.
If $\alpha$ is a limit ordinal, then consider $\repsp{X}=\repspb{Rea}(\{j_{\omega^{\beta}}\}_{\beta<\alpha})$.
\end{proof}


 \subsubsection{Jump Dimension}

We next consider a variant of countable universality to introduce another $\sigma$-topological invariant.
 A partial continuous function $f:\subseteq\mathbb{N}\times\repsp{X}\to\repsp{Y}$ is {\em countably $\mathcal{G}$-avoiding} if for any $g\in \mathcal{G}$,
\[(\exists x\in\repsp{X})(\forall n\in\mathbb{N})(\exists m\in\mathbb{N})[(\forall k\in\mathbb{N})\;(n,k,x)\in{\rm dom}(g)\;\Longrightarrow\;f(m,x)\not=g(n,k,x)].\]

\begin{proposition}\label{prop:jdim-countable-dimensional}
A Polish space $\repsp{X}$ is countable dimensional if and only if $\mathbb{N}\times\repsp{X}$ admits no countably $\mathcal{C}$-avoiding $\repsp{I}$-valued function, where $\mathcal{C}$ is the class of continuous functions.
\end{proposition}

Before proving Proposition \ref{prop:jdim-countable-dimensional}, we consider other classes of functions.
Let $\mathcal{B}_\alpha^\omega(\repsp{X})$ denote the set of partial functions of the form $g:\subseteq\mathbb{N}^2\times\repsp{X}\to\Cantor$ such that there are partial continuous functions $H:\subseteq\mathbb{N}\times\repsp{X}\to\Cantor$ and $K : \subseteq \mathbb{N} \times \Cantor \to \Cantor$ such that $g(n,k,x)=K(k,j_\alpha(H(n,x)))$ for any $(n,k,x)\in{\rm dom}(g)$.

\begin{definition}
The {\em jump dimension} ${\rm jdim}(\repsp{X})$ of $\repsp{X}$ is the supremum of countable ordinals $\alpha<\omega_1$ such that there is a countably $\mathcal{B}_\beta^{\omega}$-avoiding function $f:\subseteq\mathbb{N}\times\repsp{X}\to\Cantor$ for any $\beta<\omega^\alpha$.
\end{definition}


That the jump dimension is invariant under $\sigma$-homeomorphism can be easily verified analogously to Proposition \ref{prop:countablyuniversalfunctions} (1):

\begin{observation}
If $\repsp{X}$ $\sigma$-embeds into $\repsp{Y}$, then ${\rm jdim}(\repsp{X})\leq{\rm jdim}(\repsp{Y})$.
\end{observation}

The jump dimension gives an alternative way of understanding our construction in previous sections.

\begin{proposition}\label{thm:jdim-countable-dimension}
For every countable ordinal $\alpha>0$, there is a Pol-type Cantor manifold $\repsp{X}$ such that ${\rm jdim}(\repsp{X})=\alpha$.
 \end{proposition}

To show Propositions \ref{prop:jdim-countable-dimensional} and \ref{thm:jdim-countable-dimension}, we need the following effective interpretation of jump-dimension.
We say that $\mathcal{I}\subseteq\Cantor$ is {\em $\alpha$-principal} if there is $p\in\mathcal{I}$ such that $q\leq_Tp^{(\alpha)}$ for all $q\in\mathcal{I}$.

\begin{lemma}\label{lem:jdimequiv}
For an admissibly represented separable metrizable space $\repsp{X}$, the product $\mathbb{N}\times\repsp{X}$ admits a countably $\mathcal{B}_\alpha^\omega$-avoiding function if and only if relative to some oracle $r$, for all $z\in\Cantor$ there is a point $x\in\repsp{X}$ such that $z\in{\rm coSpec}^r(x)$ and ${\rm coSpec}^r(x)$ is not $\alpha$-principal.
\end{lemma}

\begin{proof}
Assume that $f:\subseteq\mathbb{N}\times\repsp{X}\to\Cantor$ is countably $\mathcal{B}_\alpha^\omega$-avoiding.
Since $\repsp{X}$ and $\Cantor$ are separable metrizable, $\alpha$ is countable, and $f$ is continuous, some oracle $r$ satisfies that $\repsp{X}$ and $\repsp{I}$ are $r$-computably embedded into Hilbert cube, $\alpha<\omega_1^{{\rm CK},r}$, and $f$ is computable relative to $r$.

Then, clearly $\{f(n,x):n\in\mathbb{N},\;(n,x)\in{\rm dom}(f)\}\subseteq{\rm coSpec}^r(x)$ holds for all $x\in\repsp{X}$.
Fix $z$.
Define $g(e,d,x)=\Phi_d\circ j_\alpha\circ\Phi_e^{r\oplus z}(x)$, and note that $g \in \mathcal{B}_\alpha^\omega$. Let $x\in\repsp{X}$ be a point witnessing the avoiding property of $f$ for given $g$.
Now, every $p\in{\rm coSpec}^r(x)$ is of the form $\Phi_e^{r\oplus z}(x)$ for some $e\in\mathbb{N}$.
By the avoiding property of $f$, there is $n\in\mathbb{N}$ such that $f(n,x)\not=g(e,d,x)$ for any $d\in\mathbb{N}$.
In other words, we have $f(n,x)\not=\Phi_d(p^{(\alpha)})$ for all $d\in\mathbb{N}$, i.e., $f(n,x)\not\leq_Tp^{(\alpha)}$.
This shows that ${\rm coSpec}^r(x)$ is not $\alpha$-principal since $f(n,x)\in{\rm coSpec}^r(x)$ for any $n$.

We claim that $z\in{\rm coSpec}^r(x)$, i.e., $z\leq_M(x,r)$.
Otherwise, by Lemma \ref{lem:JM-almosttotal}, $(x,r,z)$ has a Turing degree since $(x,r)$ has a continuous degree.
Let $y\in\Cantor$ be such that $(x,r,z)\equiv_My$, and let $e$ be such that $\Phi_e^{r\oplus z}(x)=y$.
In particular, $(x,r,z)\leq_M\Phi_e^{r\oplus z}(x)$, and therefore, if $p\leq_M(x,r)$, then $p\leq_T\Phi_e^{r\oplus z}(x)$.
Since $f(m,x)\leq_M(x,r)$, this shows that for any $m$, $f(m,x)=g(e,d,x)$ for some $d$.
This contradicts our assumption that $f$ countably $\mathcal{B}_\alpha^\omega$-avoiding via $x$.
Consequently, $z\in{\rm coSpec}^r(x)$.

Conversely, suppose that the condition in Lemma \ref{lem:jdimequiv} holds for $r$.
We define $f(e,x)=\Phi_e^r(x)$.
Then, we claim that $f$ is $\mathcal{B}_\alpha^\omega$-avoiding.
Given a $g \in \mathcal{B}_\alpha^\omega$, there is an oracle $z\geq_Tr$ such that there are $z$-computable $H$, $K$ with $g(n,k,x) = K(k,j_\alpha(H(n,x)))$. Then $g(n,k,x)\leq_T(h(n,x)\oplus z)^{(\alpha)}$ for any $x\in\repsp{X}$ and $n,k\in\mathbb{N}$.

By our assumption, there is a point $x\in\repsp{X}$ with $z\in{\rm coSpec}^r(x)$ such that the $r$-cospectrum of $x$ is not $\alpha$-principal.
Note that $r\leq_T z\in{\rm coSpec}^r(x)$ implies ${\rm coSpec}^r(x)={\rm coSpec}^z(x)$, and therefore, the $z$-cospectrum of $x$ is not $\alpha$-principal.
This implies that for any $n$, there is $p\in{\rm coSpec}^z(x)$ such that $p\not\leq_T(h(n,x)\oplus z)^{(\alpha)}$.
In particular, $p\not\leq_Tg(n,k,x)$ for any $k$.
Let $e$ be an index such that $p=\Phi_e^r(x)$, that is, $p=f(e,x)$.
This concludes that $f(e,x)\not=g(n,k,x)$ for any $k\in\mathbb{N}$.
\end{proof}

\begin{proof}[Proofs of Propositions \ref{prop:jdim-countable-dimensional} and \ref{thm:jdim-countable-dimension}]
For Proposition \ref{prop:jdim-countable-dimensional}, by Miller's result \cite[Proposition 5.3]{miller2}, we can deduce that a point $x\in[0,1]^\mathbb{N}$ has a Turing degree relative to $r$ if and only if ${\rm coSpec}^r(x)$ is principal (i.e., $0$-principal).
Hence, if $\repsp{X}$ is countable dimensional, all cospectra are $0$-principal up to some oracle.
Therefore, $\repsp{X}$ is not $0$-avoiding.
Conversely, suppose that $\repsp{X}$ is not countable dimensional.
We claim that for all $z\in\Cantor$, there is $x\in\repsp{X}$ such that $z\leq_Mx$ and $x$ has no Turing degree.
Otherwise, $(x,z)$ has a Turing degree by Lemma \ref{lem:JM-almosttotal}.
In this case, ${\rm Spec}^z(\repsp{X})\subseteq\mathcal{D}_T$.
This implies that $\repsp{X}$ is countable dimensional.
Now, our claim clearly implies the desired condition by Lemma \ref{lem:jdimequiv}.

For Proposition \ref{thm:jdim-countable-dimension}, combine Lemma \ref{lem:jdimequiv} and the argument in the proof of Proposition \ref{prop:udim}.
\end{proof}

\begin{example}
\begin{enumerate}
\item The universal/jump-dimension of Hilbert cube $[0,1]^\mathbb{N}$ is $\omega_1$.
This is because every countable Scott ideal is realized as a cospectrum in the Hilbert cube \cite[Theorem 9.3]{miller2} and by Lemma \ref{lem:jdimequiv}.
\item The universal/jump-dimension of $\repspb{Rea}(\mathcal{G})$ cannot be $\omega_1$ for every countable set of $\mathcal{G}$ of oracle $\mathbf{\Pi}^0_2$ singletons.
This is because every oracle $\mathbf{\Pi}^0_2$ singleton is Borel measurable.
Therefore, there is a countable ordinal $\alpha$ which bounds all Borel ranks of functions contained in $\mathcal{G}$ since $\aleph_1$ is regular.
Thus, for any $g\in\mathcal{G}$, we have $g(r)\leq_T(r\oplus z)^{(\alpha)}$ for some oracle $z$.
One can see that the cospectrum of a point in $\repspb{Rea}(\mathcal{G})$ is $(\alpha\cdot\omega)$-principal.
\item We have $0\leq{\rm jdim}({\rm \omega CEA})\leq 1$.
This is because the cospectrum of a point in $\omega{\rm CEA}$ is $\omega$-principal by the proof of Lemma \ref{lem:Scott-ideal}.
\item There is a strongly infinite dimensional Polish space such that ${\rm udim}(\repsp{X})={\rm jdim}(\repsp{X})=1$ (e.g., $\repsp{X}=\omega\repsp{CEA}(\{j\})$).
We do not know whether there exists a strongly infinite dimensional compact metric space with ${\rm udim}(\repsp{X})<\omega_1$ (or ${\rm jdim}(\repsp{X})<\omega_1$).
\end{enumerate}
\end{example}

\section{Internal Characterization of Degree Structures}\label{sec:internalcharacterization}

\subsection{Characterizing Continuous Degrees Through a Metrization Theorem}
\label{subsec:internalcharacterization}
In this section, we will provide a rather strong metrization theorem, namely that any computably admissible space with an effectively fiber-compact representation can be computably embedded in a computable metric space. Our result is a slightly stronger version of a result by \name{Schr\"oder} that an admissible space with a proper representation is metrizable \cite{schroder6}. This also gives us a characterization of the continuous degrees inside the Medvedev degrees that does not refer to represented spaces at all.

For some closed set $A \subseteq \Cantor$, let $T(A) \subseteq \Cantor$ be the set of trees for $A$, where each infinite binary tree is identified with an element of Cantor space. Now let $\delta : \subseteq \Cantor \to \repsp{X}$ be an {\em effectively fiber-compact representation}, i.e.~let $x \mapsto \delta^{-1}(\{x\}) : \repsp{X} \to \mathcal{A}(\Cantor)$ be computable. Then $T(\delta^{-1}(\{x\})) \leq_M \delta^{-1}(\{x\})$. If $\delta$ is computably admissible, we also have $\delta^{-1}(\{x\}) \leq_M T(\delta^{-1}(\{x\}))$. Note that being effectively fiber-compact is equivalent to being effectively proper, as the union of compactly many compact sets is compact. It is known that any computable metric space has a computably admissible effectively fiber-compact representation (e.g.~\cite{weihrauchf}).
We shall prove that the converse holds, too.

\begin{theorem}\label{thm:internal_characterization}
A represented space $\repsp{X}$ admits a  computably admissible effectively fiber-compact representation iff $\repsp{X}$ embeds computably into a computable metric space.
\end{theorem}

\begin{corollary}
\label{corr:contdegreecharac}
$A \subseteq \Cantor$ has continuous degree iff there is $B \in \mathcal{A}(\Cantor)$ such that $A \equiv_M B \equiv_M T(B)$.
\end{corollary}

To prove Theorem \ref{thm:internal_characterization}, we need the following two lemmata and a result by \name{Weihrauch}.
\begin{lemma}
\label{lemma:applysequence}
Let $\repsp{X}$ admit an effectively fiber-compact representation. Then there is a space $\repsp{Y}$ such that:
\begin{enumerate}
\item $\repsp{X} \hookrightarrow \repsp{Y}$ (as a closed subspace),
\item $\repsp{Y}$ has an effectively fiber-compact representation,
\item $\repsp{Y}$ has a computable dense sequence,
\item if $\repsp{X}$ is computably admissible, so is $\repsp{Y}$.
\end{enumerate}
\begin{proof}
{\bf Construction of $\repsp{Y}$:} We start with some preliminary technical notation. Let $\operatorname{Wrap} : \Cantor \to \Cantor$ be defined by $\operatorname{Wrap}(p)(2i) = p(i)$ and $\operatorname{Wrap}(p)(2i+1) = 0$. Let $\operatorname{Prefix} :\subseteq \Cantor \to \{0,1\}^*$ be defined by $\operatorname{Prefix}(p) = w$ iff $p = 0w(1)0w(2)0\ldots011q$ for some $q \in \Cantor$. Note that $\dom(\operatorname{Prefix}) \cap \dom(\operatorname{Wrap}^{-1}) = \emptyset$ and $\dom(\operatorname{Prefix}) \cup \dom(\operatorname{Wrap}^{-1}) = \Cantor$.

Let the presumed representation of $\repsp{X}$ be $\delta_\repsp{X} : \subseteq \Cantor \to X$. Our construction of $\repsp{Y}$ will utilize a notation $\nu_\repsp{Y} : \{0,1\}^* \to Y'$ as auxiliary part, this notation (or alternatively, equivalence relation on $\{0,1\}^*$) will be dealt with later. We set $Y = X \cup Y'$ (in particular, we add only countably many elements to $\repsp{X}$) and then define $\delta_\repsp{Y}$ via $\delta_\repsp{Y}(p) := \delta_\repsp{X}(\operatorname{Wrap}^{-1}(p))$ if $p \in \dom(\delta_\repsp{X}\circ\operatorname{Wrap}^{-1})$ and $\delta_\repsp{Y}(p) = \nu_\repsp{Y}(\operatorname{Prefix}(p))$ if $p \in \dom(\operatorname{Prefix})$.

In order to define $\nu_\repsp{Y}$, we do need to refer to the effective fiber-compactness of $\delta_\repsp{X}$. From the function realizing $x \mapsto \delta_\repsp{X}^{-1}(\{x\}) : \repsp{X} \to \mathcal{A}(\Cantor)$ we can obtain an indexed family of finite trees $(T_w)_{w \in \{0,1\}^*}$ with the following properties:
\begin{enumerate}
\item Each $T_w$ has height $|w|$.
\item If $w \prec u$, then $T_u \cap \{0,1\}^{\leq |w|} = T_w$.
\item $w \in T_w$.
\item For any $p \in \dom(\delta_\repsp{X})$, some $q \in \Cantor$ is an infinite path through $\bigcup_{n \in \mathbb{N}} T_{p_{\leq n}}$ iff $\delta_\repsp{X}(q) = \delta_\repsp{X}(p)$.
\end{enumerate}
Now we set $\nu_\repsp{Y}(w) = \nu_\repsp{Y}(u)$ iff $T_w = T_u$. Note in particular that $T_w = T_u$ is a decidable property.

{\bf Proof of the properties:} To see that $\repsp{X} \hookrightarrow \repsp{Y}$ it suffices to note that both $\operatorname{Wrap}$ and $\operatorname{Wrap}^{-1}$ are computable. That $\repsp{X}$ embeds as a closed subspace follows from $\dom(\operatorname{Wrap}^{-1})$ being closed in $\Cantor$.

Next we shall see that $\delta_\repsp{Y}$ is effectively fiber-compact by reversing the step from the function $x \mapsto \delta_\repsp{X}^{-1}(\{x\}) : \repsp{X} \to \mathcal{A}(\Cantor)$ to the family $(T_w)_{w \in \{0,1\}^*}$. First, we define a version of $\operatorname{Wrap}$ for finite trees via $\operatorname{T-Wrap}(T) = \{0w(1)0w(2)\ldots w(|w|)\mid w \in T\} \cup \{0w(1)0w(2)\ldots w(|w|)0\mid w \in T\}$. Given some set $W \subseteq \{0,1\}^*$, let the induced tree of height be defined via $T(W,n) = \{u \exists w \in W \ u \prec w\} \cup \{u \in \{0,1\}^n \mid \exists w \in W \ \wedge w \prec u\}$. Then we define a derived family $(T'_w)_{w \in \{0,1\}^*}$ by $T'_{0w(1)0w(2)\ldots0w(|w|)} = T'_{0w(1)0w(2)\ldots0w(|w|)0} = \operatorname{T-Wrap}(T_w)$ and $T'_{0w(1)\ldots w(|w|)1v} = T(\{u \mid T_u = T_w\},|w|+|v|)$. This construction too satisfies that if $w \prec u$, then $T'_u \cap \{0,1\}^{\leq \operatorname{height}(T'(w))} = T'_w$. Thus, the function that maps $p$ to the set of all infinite pathes through $\bigcup_{n \in \mathbb{N}} T'_{p_{\leq n}}$ does define some function $t : \Cantor \to \mathcal{A}(\Cantor)$, and one can verify readily that $t(p) = \delta_\repsp{Y}^{-1}(\delta_\repsp{Y}(p))$ whenever $p \in \dom(\delta_\repsp{Y}(p))$.

It is clear that $\repsp{Y}$ has a computable dense sequence: Fix some standard enumeration $\nu : \mathbb{N} \to \{0,1\}^*$, and consider $(y_n)_{n \in \mathbb{N}}$ with $y_n = \delta_\repsp{Y}(0\nu(n)(1)\ldots0\nu(n)(|\nu(n)|)1^\omega)$.

It remains to show that if $\delta_\repsp{X}$ is admissible, so is $\delta_\repsp{Y}$. It is this step which requires the identification of some points via $\nu_\repsp{Y}$, and through this, also depends on $\delta_\repsp{X}$ being effectively-fiber-compact. Given some tree encoding some $\delta_\repsp{Y}^{-1}(\{x\})$, we need to be able to compute a path through it. As long as the tree seems to have a path without repeating $1$'s, we lift the corresponding map for $\delta_\repsp{X}$. If $x = \nu_\repsp{Y}(w)$ for some $w \in \{0,1\}^*$, we notice eventually, and can extend the current path in a computable way by virtue of the identifications.
\end{proof}
\end{lemma}

The preceding lemma produces spaces with a somewhat peculiar property: The designated dense sequence is an open subset of the space, unlike the usual examples. In \cite{gregoriades3}, \name{Gregoriades} has explored a general construction yielding Polish spaces with such properties (cf.~\cite[Theorem 2.5]{gregoriades3}), which in particular serves to prevent effective Borel isomorphisms between spaces.

\begin{lemma}
\label{lemma:regular}
Let $\repsp{X}$ admit a computably admissible effectively fiber-compact representation. Then $\repsp{X}$ is computably regular.
\begin{proof}
The properties of the representations mean that we can consider $\repsp{X}$ as a subspace of $\mathcal{A}(\Cantor)$ containing only pair-wise disjoint sets. Let $A \in \mathcal{A}(\mathcal{A}(\Cantor))$ be a closed subset in $\repsp{X}$. Note that we can compute $\bigcup A \in \mathcal{A}(\Cantor)$, as every infinite path computes the relevant tree. Furthermore, given $x \in \repsp{X} \subseteq \mathcal{A}(\Cantor)$ and $A$, we can compute $x \cap \bigcup A \in \mathcal{A}(\Cantor)$. As $\repsp{X}$ only contains pair-wise disjoint points, this set is empty if and only if $x$ and $A$ are disjoint. As $\Cantor$ is compact, the corresponding tree will have to die out at some finite level, which means that the trees for $x$ and $A$ are disjoint below this level. Let $I$ be the vertices at this level belonging to $x$. We may now define two open sets $U_I, U_{I^C} \in \mathcal{O}(\mathcal{A}(\Cantor))$ by letting $U_{X}$ for $X \in \{I, I^C\}$ accept its input sets $A$ as soon as $A \cap X\Cantor = \emptyset$ is verified. Then $U_{I} \cap U_{I^C} = \{\emptyset\}$, thus $\repsp{X} \cap U_I$ and $\repsp{X} \cap U_{I^C}$ are disjoint open sets. Moreover, we find $A \subseteq U_{I}$ and $x \in U_{I^C}$, so the two open sets are those we needed to construct for computable regularity.
\end{proof}
\end{lemma}

\begin{proof}[Proof of Theorem \ref{thm:internal_characterization}]
The $\Leftarrow$-direction is present e.g.~in \cite{weihrauchf}. We can use Lemma \ref{lemma:applysequence} to make sure w.l.o.g.~that $\repsp{X}$ has a computable dense sequence. By Lemma \ref{lemma:regular}, the space is computably regular. As shown in \cite{grubba3,weihrauchm}, a computably regular space with a computable dense sequence admits a compatible metric.
\end{proof}

\name{Miller} showed that the Turing degrees below any non-total continuous degree form a Scott ideal \cite{miller2}, heavily drawing on topological arguments. However, based on Corollary \ref{corr:contdegreecharac} we see that the statement itself can be phrased entirely in the language of trees, points and Medvedev reducibility. So far, we do not know of a direct proof involving only these concepts:

\begin{proposition}
Let $A \subseteq \Cantor$ be such that $A \equiv_M T(A)$ and that there is no $r \in \Cantor$ with $A \equiv_M \{r\}$. Then $T(B) \leq_M \{p\} <_M A$ for $p \in \Cantor$, $B \subseteq \Cantor$ implies $B \leq_M A$.
\end{proposition}

\subsection{Enumeration Degrees and Overtness}
\label{subsec:fiberovert}
The often overlooked dual notion to compactness is \emph{overtness} (see \cite{taylor,taylor2}). Intuitively, overtness makes existential quantification well-behaved: a space $\repsp{X}$ is overt if $E_\repsp{X}:\mathcal{O}(\repsp{X})\to\mathcal{S}$ is continuous, where $E_\repsp{X}(U)=\top$ iff $U$ is nonempty.  Therefore, if $\repsp{X}$ is overt and $P \subseteq \repsp{X} \times \repsp{Y}$ is open, then $\{y \in \repsp{Y} \mid \exists x \in \repsp{X} \ (x,y) \in P\}$ is open, too. Classically, this is a trivial notion, however, the situation is different from an effective point of view.

One may identify an overt subspace $A$ of $\repsp{X}$ with $E_A$, or equivalently, its overtness witness $\{U\in\mathcal{O}(\repsp{X}):A\cap U\not=\emptyset\}$ as a point in the represented space $\mathcal{O}(\mathcal{O}(\repsp{X}))$.
Via this identification, we obtain the hyperspace $\mathcal{V}(\repsp{X})$ of representatives $\overline{A}$ of all overt subspaces $A$ of $\repsp{X}$ (see also \cite{pauly-synthetic-arxiv}).
Note that this corresponds to the lower Vietoris topology on the hyperspace of closed sets.
A computable point in $\mathcal{V}(\repsp{X})$ is also called a {\em c.e.~closed} set in computable analysis.

Now we call a representation $\delta : \subseteq \Baire \to \repsp{X}$ \emph{effectively fiber-overt}, iff $\overline{\delta^{-1}} : \repsp{X} \to \mathcal{V}(\Baire)$ is computable. A straightforward argument shows that this is equivalent to $\delta$ being effectively open, i.e.~$U \mapsto \delta[U] : \mathcal{O}(\Baire) \to \mathcal{O}(\repsp{X})$ being computable. Now we see that every space with an effectively fiber-overt representation inherits an effective countable basis from $\Baire$, while on the other hand, the standard representations of countably based spaces introduced in Example \ref{example:representation} are all effectively fiber-overt. Thus we see that while effectively fiber-compact representations characterize metrizability, effectively fiber-overt representations characterize second-countability.

\section{Point Degree Spectra of Quasi-Polish Spaces}\label{sec:quasi-Polish}

\subsection{Lower Reals and Semirecursive Enumeration Degrees}

Let us move on to the $\sigma$-homeomorphic classification of quasi-Polish spaces \cite{debrecht6}.
We now focus on the following chain of quasi-Polish spaces:
\[\mathcal{D}_T\subsetneq\mathcal{D}_r\subsetneq\mathcal{D}_e\mbox{, and }\Cantor<_\sigma^\mathfrak{T}[0,1]^\mathbb{N}<_\sigma^\mathfrak{T}\mathcal{O}(\mathbb{N}).\]

Here, the proper inclusion $[0,1]^\mathbb{N}<_\sigma^\mathfrak{T}\mathcal{O}(\mathbb{N})$ follows from relativizing \name{Miller}'s observation in \cite{miller2} that no quasi-minimal degree has continuous degree.

In quasi-Polish case, the notion of the specialization order is quite useful.
Indeed, \name{Motto Ros} has already used the specialization order to give an alternative way to show the properness of $[0,1]^\mathbb{N}<_\sigma^\mathfrak{T}\mathcal{O}(\mathbb{N})$.
Recall that the {\em specialization order} $\prec_\repsp{X}$ on a topological space $\repsp{X}$ is defined via $x \prec_\repsp{X} y :\Leftrightarrow x \in \overline{\{y\}}$.
(Then $[0,1]^\mathbb{N}<_\sigma^\mathfrak{T}\mathcal{O}(\mathbb{N})$ follows from the observation that the specialization order on $\mathcal{O}(\mathbb{N})$ coincides with subset-inclusion, while the $T_1$ separation property asserts that no two elements are comparable w.r.t.~$\prec_\repsp{X}$, i.e.~that the specialization order is a single antichain.)

We will show that, inside a single specialization order type, there are continuum many incomparable $\sigma$-homeomorphism types of quasi-Polish spaces which do not $\sigma$-embed into the Hilbert cube.

\begin{theorem}\label{thm:quasi-polish-main}
There is a map $\repsp{Q}$ transforming each countable set $S\subseteq\omega_1$ into a nonmetrizable quasi-Polish space $\repsp{Q}(S)$ such that for any countable sets $S,T\subseteq\omega_1$,
\begin{align*}
&\repsp{Q}(S)\not\leq_\sigma^\mathfrak{T}[0,1]^\mathbb{N},\qquad(\repsp{Q}(S),\preceq_{\repsp{Q}(S)})\simeq(\repsp{Q}(T),\preceq_{\repsp{Q}(T)}),\\
&S\not\subseteq T\;\Longrightarrow\;\repsp{Q}(S)\not\leq_\sigma^\mathfrak{T}\repsp{Q}(T).
\end{align*}
\end{theorem}

Let $\mathbb{R}_<$ be the real line endowed with the lower topology, that is, its topology is generated by open intervals of the form $(p,\infty)$.
One can easily see that $\mathbb{R} \ |_{\sigma}^\mathfrak{T} \ \mathbb{R}_<$ by comparing their specialization orders.
From the computability theoretic viewpoint, the property $\mathbb{R}\not\leq^\mathfrak{T}_\sigma\mathbb{R}_<$ can be strengthened as follows.

\begin{lemma}[Co-spectrum Preservation]\label{thm:cospectrum-preservation}
Let $\repsp{X}$ be an admissibly represented Polish space.
Then,
\[{\rm coSpec}(\repsp{X}\times\mathbb{R}_{<})\subseteq{\rm coSpec}(\repsp{X})\cup{\rm coSpec}(\Cantor).\]

In particular, if such an $\repsp{X}$ is uncountable, then there is an oracle $r\in\Cantor$ such that
\[{\rm coSpec}^r(\repsp{X}\times\mathbb{R}_{<})={\rm coSpec}^r(\repsp{X}).\]
\end{lemma}

\begin{lemma}\label{lem:cospectrum-preservation}
Let $\repsp{X}$ admit an effectively fiber-overt representation $\delta_\repsp{X}$ (cf.~Subsection \ref{subsec:fiberovert}), $x\in \repsp{X}$, $y\in\mathbb{R}_<$, and $z\in\Cantor$.
If $z\leq_M(x,y)$, then either $z\leq_M x$ or $-y\leq_M x$ holds.
\end{lemma}

\begin{proof}
Let computable $f :\subseteq \repsp{X} \times \mathbb{R}_< \to \Cantor$ witness the reduction $z\leq_M(x,y)$.
By extending the domain of $f$ if necessary, it can be identified with a c.e.~open set $U\subseteq\Cantor\times\mathbb{Q}\times\mathbb{N}\times\{0,1\}$ satisfying that $f(x,y)(n)=i$ if and only if the following two condition holds:

\begin{enumerate}
\item For any $p \in \delta_\repsp{X}^{-1}(x)$ there is some rational $s < y$ such that $(p, s, n, i) \in U$.
\item For any $p \in \delta_\repsp{X}^{-1}(x)$ and any rational $s<y$, $(p, s, n, 1-i) \not\in U$.
\end{enumerate}

As $\delta_\repsp{X}$ is effectively fiber-overt, the set $U' := \{(x',t,n,i) \mid \exists p \in \delta_\repsp{X}^{-1}(x') \ (p,t,n,i) \in U\}$ is also computable as an open subset of $\repsp{X} \times \mathbb{Q} \times \mathbb{N} \times \{0,1\}$. Now we can distinguish two cases:

\begin{enumerate}
\item For any $\varepsilon>0$, there exist rationals $t<s<y+\varepsilon$ such that $(x,t,n,i) \in U'$ and $(x,s,n,1-i) \in U'$ for some $n\in\mathbb{N}$ and $i\in\{0,1\}$.
\item Otherwise, there exists $\varepsilon>0$ such that for all $t<y+\varepsilon$, if $(x,t,n,i) \in U'$ for some $n\in\mathbb{N}$ and $i\in\{0,1\}$, then we must have $i = f(x,y)(n)$.
\end{enumerate}

Note that if $(x, t, n, i) \in U'$ and $(x,s,n,1-i) \in U'$ for $t < s$, then we automatically have $y \leq s$.
Therefore, in the first case we can compute $-y \in \mathbb{R}_<$ as the supremum of $-s$ over all witnesses $s$, thus find $-y \leq_M x$. In the second case, there will be some rational number $y_0$ with $y \leq y_0 < y + \varepsilon$. Using $y_0$ in place of $y$ leaves the value $f$ is producing unchanged, thus we have that $z \leq_M x$.
\end{proof}

\begin{proof}[Proof of Lemma \ref{thm:cospectrum-preservation}]
Suppose that $y\in\mathbb{R}_<$ and $x\in\repsp{X}$.
If $-y\not\leq_Mx$, then ${\rm coSpec}(x,y)={\rm coSpec}(x)$ by Lemma \ref{lem:cospectrum-preservation}.
Otherwise, $(x,y,-y)\equiv_M(x,y)$.
If $y\leq_Mx$, then clearly, ${\rm coSpec}(x,y)={\rm coSpec}(x)$.
Otherwise, $(y,-y)\not\leq_Mx$.
Obviously, $(y,-y)$ has Turing degree.
By Lemma \ref{lem:JM-almosttotal}, we have $(x,y,-y)\in\mathcal{D}_T$.
Hence, ${\rm coSpec}(x,y)\in{\rm coSpec}(\Cantor)$.
For the latter half of Lemma \ref{thm:cospectrum-preservation}, if $\repsp{X}$ is uncountable, then there is an $r$-computable embedding of $\Cantor$ into $\repsp{X}$ for some oracle $r$.
\end{proof}

\begin{proof}[Proof of Theorem \ref{thm:quasi-polish-main}]
Let $\mathcal{G}_S$ be the countable set of monotone oracle $\mathbf{\Pi}^0_2$ singletons constructed in the proof of Theorem \ref{thm:maintheorem_emb}.
Define $\repsp{Q}(S):=\repspb{Rea}(\mathcal{G}_S)\times\overline{\mathbb{R}}_<$, where $\overline{\mathbb{R}}_<:=\mathbb{R}_<\cup\{\infty\}$ is a quasi-completion of $\mathbb{R}_<$.
First note that $\mathbb{R}_<$ does not $\sigma$-embed into the Hilbert cube, since by Lemma \ref{lem:cospectrum-preservation}, for any $y\in\mathbb{R}_<$ and oracle $x$, if $z\in{\rm coSpec}^x(y)$ then either $z\leq_Tx$ or $-y\leq_Mx$, that is, $-y$ is left-c.e.\ relative to $x$.
Given an oracle $x$, there are only countably many $y$ such that $-y\leq_Mx$, and therefore, this means that given $x$, the $x$-co-spectrum of $y$ consists only of $x$-computable points for all but countably many $y$.
If $y$ is not $x$-computable, such a $y$ is called quasi-minimal relative to $x$.
However, the Hilbert cube does not contain such a point as shown in \cite[Corollary 7.3]{miller2}.
Since quasi-minimality is a degree-theoretic property, this concludes $\mathbb{R}_<\not\leq^\mathfrak{T}_\sigma[0,1]^\mathbb{N}$ by Theorem \ref{theo:spectrum-main}.

Concerning specialization orders, $\repspb{Rea}(\mathcal{G}_S)$ is metrizable (hence $T_1$), its specialization order is a single antichain of cardinality continuum.
It is easy to see that the specialization order on a product space $A\times B$ is the product of the specialization orders on $A$ and $B$.
Thus, the specialization order on $\repsp{Q}(S)$ is order-isomorphic to that on $\repsp{Q}(T)$.

Finally, by Lemma \ref{thm:cospectrum-preservation}, the quasi-Polish space $\repsp{Q}(S)=\repspb{Rea}(\mathcal{G}_S)\times\overline{\mathbb{R}}_<$ has the same cospectrum as $\repspb{Rea}(\mathcal{G}_S)$.
By the proofs of Theorem \ref{thm:maintheorem_emb} and Lemma \ref{thm:refle}, if $S\not\subseteq T$, then the cospectrum of $\repspb{Rea}(\mathcal{G}_S)$ is not a sub-cospectrum of $\repspb{Rea}(\mathcal{G}_T)$ relative to all oracles.
Therefore, by Observation \ref{obs:main-cospec-inv}, we have $\repsp{Q}(S)\not\leq_\sigma^\mathfrak{T}\repsp{Q}(T)$.
\end{proof}

As a consequence of Lemma \ref{thm:cospectrum-preservation}, any lower real can compute only a $\Delta^0_2$ real:
\[{\rm coSpec}(\mathbb{R}_{<})=\{\{x\in \Cantor:x\leq_Ty\}:y\mbox{ is right-c.e.}\}\]

Indeed, Lemma \ref{lem:cospectrum-preservation} provides a very simple and natural construction of a quasi-minimal enumeration degree.

\begin{corollary}[see also {\name{Arslanov}, \name{Kalimullin} \& \name{Cooper} \cite[Theorem 4]{kalimullin}}]
Suppose that $z\in\mathbb{R}$ is neither left-c.e nor right-c.e.
Then, the enumeration degree of the cut $\{q\in\mathbb{Q}:q<z\}$ is quasi-minimal.
\end{corollary}

On the one hand, we deduced the property $[0,1]^\mathbb{N}<^\mathfrak{T}_\sigma\mathcal{O}(\mathbb{N})$ from the topological argument concerning the specialization order on the lower real $\mathbb{R}_<$.
On the other hand, \name{Miller}'s original proof used the existence of a quasi-minimal enumeration degree to show ${\rm Spec}([0,1]^\mathbb{N})\subsetneq{\rm Spec}(\mathcal{O}(\mathbb{N}))$.
Surprisingly, however, the previous argument clarifies that these two seemingly unrelated approaches are essentially equivalent.

Note that the point degree spectrum of the lower real $\mathbb{R}_<$ is indeed strongly connected with the notion of a semirecursive set in the context of the enumeration degrees.
Recall from \cite{jockusch2} that a set $A \subseteq \mathbb{N}$ is called {\em semirecursive}, if there is a computable function $f : \mathbb{N} \times \mathbb{N} \to \mathbb{N}$ such that for all $n, m \in \mathbb{N}$ we find $f(n,m) \in \{n,m\}$, and if $n \in A$ or $m \in A$, then $f(n,m) \in A$. We call an enumeration degree $q \in \mathcal{D}_e$ semirecursive, if it is the degree of a semirecursive point in $\mathcal{O}(\mathbb{N})$.

\name{Jockusch} \cite{jockusch2} pointed out that every left-cut (i.e., every lower real $x\in\mathbb{R}_<$) is semirecursive, and conversely, \name{Ganchev} and \name{Soskova} \cite{GanSos15} showed that every semirecursive enumeration degree contains a left-cut.
Consequently, the point degree spectra of the lower real $\mathbb{R}_<$ can be characterized as follows:
\[{\rm Spec}(\mathbb{R}_<)=\{\mathbf{d}\in\mathcal{D}_e:\mathbf{d}\mbox{ is a semirecursive enumeration degree}\}.\]

\subsection{Higher Dimensional Lower Cubes}

We can also consider the higher dimensional lower real cubes $\mathbb{R}_<^n$.
Surprisingly, the spectra of $\mathbb{R}_<^n$ form a proper hierarchy as follows.

\begin{theorem}\label{thm:lowerrealhierarchy}
If $\repsp{X}$ is a second-countable $T_1$ space, then $\mathbb{R}_<^{n+1}\ |_\sigma^\mathfrak{T} \ \repsp{X}\times\mathbb{R}_<^{n}$ for every $n$.
\end{theorem}

To show the above theorem, we use the following order theoretic lemma.
Let $\Lambda^n=(\{0,1\}^n,\leq)$ be a partial order on $\{0,1\}^n$ obtained as the $n$-th product of the ordering $0<1$.

\begin{lemma}\label{lem:lowerrealhierarchy}
For every countable partition $(P_i)_{i\in\omega}$ of the $n$-dimensional hypercube $[0,1]^n$ (endowed with the standard product order), there is $i\in\omega$ such that $P_i$ has a subset which is order isomorphic to the product order $\Lambda^n$.
\end{lemma}

\begin{proof}
We use Vaught's ``non-meager'' quantifier $\exists^*x\varphi(x)$, which states that the set $\{x:\varphi(x)\}$ is not meager in $[0,1]$ (with respect to the standard Euclidean topology).
We claim that for every countable partition $(P_i)_{i\in\omega}$ of $[0,1]^n$, there is $i\in\omega$ such that
\[\exists^*x_1\exists^*x_2\dots\exists^*x_n\;(x_1,x_2,\dots,x_n)\in P_i\]

Inductively assume that the above claim is true for $n-1$.
If the above claim does not hold for $n$, then by the Baire category theorem, there are comeager many $x_1$ such that
\[\neg\exists^*x_2\dots\exists^*x_n\;(x_1,x_2,\dots,x_n)\in\bigcup_iP_i.\]

However, for any such $x_1$, by the induction hypothesis, the $x_1$-sections of $P_i$'s do not cover the $x_1$-section of $[0,1]^n$.
In particular, $\bigcup_iP_i$ cannot cover the $n$-hypercube $[0,1]^n$, which verifies the claim.

Now, let $S$ be a nonmeager set consisting of all $x_1$'s in the above claim.
Note that since there are non-meager many $x_1\in S$, there is a nonempty open set $U$ such that for any nonempty open set $V\subseteq U$, one can find uncountably many such $x_1\in V\cap S$.
Otherwise, $S$ is covered by the closure of the union of the collection $\mathcal{B}$ of all rational open balls $B$ such that $B\cap S$ is countable.
Therefore, $S$ is divided into the union of the nowhere dense set $\partial\bigcup\mathcal{B}$ and the countable set $\bigcup_{B\in\mathcal{B}}B\cap S$, which contradicts the fact that $S$ is nonmeager.
We fix such a nonempty open set $U$.

Now, for any $x_1\in S$, we may inductively assume that the $x_1$-th section of $P_i$ has a subset $L(x_1)$ which is order isomorphic to $\Lambda^{n-1}$.
Let $\hat{L}(x_1)$ be the region bounded by $L(x_1)$, which is homeomorphic to $[0,1]^{n-1}$.
We may also inductively assume that $P_i$ is dense in $\hat{L}(x_1)$.
Therefore, since $\hat{L}(x_1)$ for any $x_1\in S$ has positive $(n-1)$-dimensional Lebesgue measure, for any nonempty open set $V\subseteq U$ one can find $x_1^0<x_1^1$ in $V\cap S$ such that the intersection $\pi\circ\hat{L}(x_1^0)\cap\pi\circ\hat{L}(x_1^1)$ also has positive $(n-1)$-dimensional Lebesgue measure, where $\pi:[0,1]^n\to[0,1]^{n-1}$ is the projection defined by $\pi(x_1,x_2,\dots,x_n)=(x_2,\dots,x_n)$.
By density of $P_i$, one can find a smaller $(n-1)$-cubes $L^*(x^0_1),L^*(x^1_1)\subseteq\pi\circ\hat{L}(x_1^0)\cap\pi\circ\hat{L}(x_1^1)$ such that $(\{x^0_1\}\times L^*(x^0_1))\cup(\{x^1_1\}\times L^*(x^1_1))\subseteq P_i$ is order isomorphic to $\Lambda^{n}$.
\end{proof}

\begin{proof}[Proof of Theorem \ref{thm:lowerrealhierarchy}]
Note that the specialization order on the space $\mathbb{R}_<^{n+1}$ is exactly the same as the standard product order on $\mathbb{R}^{n+1}$.
By Lemma \ref{thm:lowerrealhierarchy}, for every countable partition $(P_i)_{i\in\omega}$ of $\mathbb{R}_<^{n+1}$, there is $i\in\omega$ such that the specialization order on $P_i$ has a subset which is order isomorphic to the product order $\Lambda^{n+1}$ whose order dimension is $n+1$.
If $P_i$ is embedded into the specialization order on $\repsp{X}\times\mathbb{R}_<^n$, then the embedded image of an isomorphic copy of $\Lambda^{n+1}$ has to be contained in a connected component of the order of $\repsp{X}\times\mathbb{R}_<^n$.
However, the specialization order on $\repsp{X}\times\mathbb{R}_<^n$ is now ${\rm card}(\repsp{X})$ many copies of that on $\mathbb{R}_<^n$ since $\repsp{X}$ is $T_1$.
Therefore, every connected component of the specialization order on $\repsp{X}\times\mathbb{R}_<^n$ is isomorphic to the product order on $\mathbb{R}^n$ whose order dimension is $n$.
Hence, $\mathbb{R}_<^{n+1}$ cannot be $\sigma$-embedded into $\repsp{X}\times\mathbb{R}_<^n$.

Conversely, suppose that $\repsp{X}\times\mathbb{R}_<^n$ is $\sigma$-embedded into $\mathbb{R}_<^{n+1}$.
By Lemma \ref{thm:lowerrealhierarchy}, for every countable partition $(P_i)_{i\in\omega}$ of $\repsp{X}\times\mathbb{R}_<^{n}$, there must exist $i\in\omega$ such that $P_i$ contains a uncountable family $(\Lambda^n_\alpha)_{\alpha\in\aleph_1}$ of pairwise incomparable suborders of $P_i$ which are order isomorphic to $\Lambda^n$.

Let $L_\alpha$ be the embedded image of $\Lambda^n_\alpha$ in $\mathbb{R}_<^{n+1}$, and $\hat{L}_\alpha$ be the region bounded by $L_\alpha$, which is homeomorphic to $[0,1]^n$.
As in the proof of Lemma \ref{thm:lowerrealhierarchy}, we may also assume that the embedded image $P_i^*\subseteq\mathbb{R}_<^{n+1}$ of $P_i$ is dense in $\hat{L}_\alpha$ for any $\alpha<\aleph_1$.
For any $\alpha<\aleph_1$, the projection $\pi_k[\hat{L}_\alpha]=\{(x_0,\dots,x_{k-1},x_{k+1},\dots,x_{n}):(x_0,\dots,x_n)\in\hat{L}_\alpha\}$ of $\hat{L}_\alpha$ for some $k\leq n$ has positive $n$-dimensional Lebesgue measure.
Fix $k<n+1$ such that $\pi_k[\hat{L}_\alpha]$ has positive $n$-dimensional Lebesgue measure for uncountably many $\alpha$.
Then, there are $\alpha\not=\beta$ such that $\pi_k[\hat{L}_\alpha]\cap\pi_k[\hat{L}_\beta]$ also has positive $(n+1)$-dimensional Lebesgue measure.
It is not hard to see that it contradicts our assumption that $L_\alpha$ and $L_\beta$ are incomparable.
\end{proof}

Note that Theorem \ref{thm:lowerrealhierarchy} has immediate computability-theoretic corollaries:

\begin{corollary}
For every $n \in \mathbb{N}$ there are enumeration degrees $p_n$, $q_n$ such that
\begin{itemize}
\item $p_n$ is the product of $n + 1$ semirecursive degrees, but not of $n$ semirecursive degrees and a Turing degree.
\item $q_n$ is the product of $n$ semirecursive degrees and a Turing degree, but not of $n-1$ semirecursive degrees and a Turing degree, or of $n + 1$ semirecursive degrees.
\end{itemize}
\end{corollary}

\subsection{The co-spectrum of a Universal Quasi-Polish Space}

Recall that the co-spectrum of the universal Polish space $[0,1]^\mathbb{N}$ consists of all principal countable Turing ideals and all countable Scott ideals.
However, there are many non-principal countable Turing ideals that are not Scott ideals, e.g., countable $\omega$-models of ${\sf WWKL}+\neg{\sf WKL}$, ${\sf RT}^2_2+\neg{\sf WKL}$ and so on.
We now see that every countable Turing ideal is realized as a co-spectrum of the universal quasi-Polish space $\mathcal{O}(\mathbb{N})$ by modifying the standard forcing construction of quasi-minimal enumeration degrees.

\begin{theorem}\label{thm:cospec-univscottdom}
Every countable Turing ideal is realized as a co-spectrum in the universal quasi-Polish space $\mathcal{O}(\mathbb{N})$.
In particular, ${\rm coSpec}^r([0,1]^\mathbb{N})\subsetneq{\rm coSpec}^r(\mathcal{O}(\mathbb{N}))$ for every $r\in\mathbb{N}$.
\end{theorem}

\begin{proof}
It suffices to show that, for any sequence $(x_i)_{i\in\mathbb{N}}$ of reals and oracle $r$, there is $A\in\mathcal{O}(\mathbb{N})$ whose $r$-co-spectrum ${\rm coSpec}^r(A)$ is equal to all $y\in \Cantor$ such that $y\leq_Tr\oplus\bigoplus_{m\leq n}x_m$ for some $n\in\mathbb{N}$.
Without loss of generality, we may assume that $x_0=r$.
Suppose $\bot\not\in\mathbb{N}$, and let $\mathbb{N}_\bot=\mathbb{N}\cup\{\bot\}$.
We say that a sequence $\sigma\in\mathbb{N}_\bot$ strongly extends $\tau\in\mathbb{N}_\bot$ if $\tau$ is an initial segment of $\sigma$ as a $\mathbb{N}_\bot$-valued sequence.
A sequence $\sigma\in\mathbb{N}_\bot$ extends $\tau\in\mathbb{N}_\bot$ if $\sigma$ extends $\tau$ as a partial function on $\mathbb{N}$, where the equality $\sigma(n)=\bot$ is interpreted as meaning that $\sigma(n)$ is undefined, that is, $n\not\in{\rm dom}(\sigma)$.

Every partial function $\varphi:\subseteq\mathbb{N}\to\mathbb{N}$ generates a tree $T_\varphi\subseteq\mathbb{N}_\bot^{<\omega}$ by
\[T_\varphi=\{\sigma\in\mathbb{N}_\bot^{<\omega}:(\forall n<|\sigma|)\;\varphi(n)\downarrow\;\rightarrow\;\sigma(n)=\varphi(n)\}.\]


Let $\mathbb{P}$ be the collection of pairs $(\sigma,\varphi)$ of a string $\sigma\in \mathbb{N}_\bot^{<\omega}$ and a partial function $\varphi$ such that $\sigma\in T_\varphi$ and ${\rm dom}(\varphi)$ is of the form $D(A)=\{(m,n):n\in\mathbb{N}\;\&\;m\in A\}$ for some finite set $A\subseteq\mathbb{N}$.
We write $(\tau,\psi)\leq(\sigma,\varphi)$ if $\tau$ strongly extends $\sigma$, $\psi$ extends $\varphi$, and $\psi\upharpoonright|\sigma|=\varphi\upharpoonright|\sigma|$.

By induction, we assume that $(\sigma_0,\varphi_0)$ is the pair of an empty string and an empty function, and $(\sigma_s,\varphi_s)\in\mathbb{P}$ has already been defined.
Moreover, we inductively assume that the tree $T_{\varphi_s}$ is computable in $\bigoplus_{2t<s}x_t$.
We now have ${\rm dom}(\varphi_s)=D(A_s)$ for some $s\in\mathbb{N}$ by the definition of $\mathbb{P}$.
If $s=2e$ for some $e\in\mathbb{N}$, then choose sufficiently large $m_{s+1}\not\in A_s$ with $m_{s+1}>|\sigma_{s}|$.
Then, put $\sigma_{s+1}=\sigma_s$, and define $\varphi_{s+1}(m_{s+1},n)=x_e(n)$ for every $n\in\mathbb{N}$.
Clearly, the tree $T_{\varphi_{s+1}}$ is computable in $\bigoplus_{2t\leq s}x_t$.

If $s=2e+1$ for some $e\in\mathbb{N}$, we look for a string $\tau\in T_{\varphi_s}$ strongly extending $\sigma_s$ which forces the $e$-th computation $\Psi_e$ to be inconsistent, that is, two different values $\Psi_e(\tau)(n)=i$ and $\Psi_e(\tau)(n)=j$ for some $n$ and $i\not=j$ are enumerated.
If there is such a $\tau$, define $\sigma_{s+1}=\tau$ and $\varphi_{s+1}=\varphi_s$.

If there is no such a $\tau$, we look for strings $\eta,\theta\in T_{\varphi_s}$ strongly extending $\sigma_{s}$ such that the $e$-th computations $\Psi_e$ on $\eta$ and $\theta$ split and are consistent, that is, the consistent computations $\Psi_e(\eta)(n)=i$ and $\Psi_e(\theta)(n)=j$ for some $n$ and $i\not=j$ are enumerated.
In this case, for a sufficiently large $k>\max{|\eta|,|\theta|}$, define $\sigma_{s+1}$ to be the rightmost node of $T_{\varphi_s}$ strongly extending $\sigma_s$, where we declare that $\bot$ is the rightmost element in $\mathbb{N}_\bot$ in the sense that $n<\bot$ for every $n\in\mathbb{N}$.
Note that $\eta$ (resp.~$\theta$) (non-strongly) extends $\sigma_{s+1}\upharpoonright|\eta|$ (resp.~$\sigma_{s+1}\upharpoonright|\theta|$) since $\sigma_{s+1}$ chooses as many $\bot$'s as possible.
Then, define $\varphi_{s+1}=\varphi_s$.

Otherwise, define $\sigma_{s+1}=\sigma_s$ and $\varphi_{s+1}=\varphi_s$.
Finally, we obtain a partial function $\Phi$ on $\mathbb{N}$ by combining $\{\varphi_s\}_{s\in\mathbb{N}}$.

As in the usual argument, we will show that $\Phi$ is quasi-minimal above the collection $\{\bigoplus_{m\leq n}x_m:n\in\mathbb{N}\}$.
Clearly, $\bigoplus_{m\leq e}x_m$ is computable in $\Phi$ by our strategy at stage $2e$.

To show quasi-minimality of $\Phi$, consider the $e$-th computation $\Psi_e$.
If we find an inconsistent computation on some $\tau$ at stage $s=2e+1$, then clearly, $\Psi_e(\Phi)$ does not define an element of $\Cantor$.
If we find a consistent $e$-splitting $\eta$ and $\theta$ on an input $n$ at stage $s=2e+1$, $\Psi_e(\Phi)(n)$ is undefined, since otherwise $\Psi_e(\Phi)(n)=k$ implies $\Psi_e(\eta)=\Psi_e(\theta)=k$.
Otherwise, for every $n\in\mathbb{N}$, if $\Psi_e(\Phi)(n)$ is defined, then it is consistent, and uniquely determined inside $T_{\varphi_{s}}$.
Therefore, $\Psi_e(\Phi)(n)=k$ if and only if there is $\tau\in T_{\varphi_s}$ strongly extending $\sigma_s$ such that $\Psi_e(\tau)(n)=k$.
Consequently, $\Psi_e(\Phi)$ is computable in $\bigoplus_{m\leq e}x_m$, since $T_{\varphi_s}$ is a pruned $\bigoplus_{m\leq e}x_m$-computable tree by induction.
\end{proof}

\begin{corollary}
For any separable metrizable space $\repsp{X}$, we have $\repsp{X}\times\mathbb{R}_<<_\sigma^\mathfrak{T}\mathcal{O}(\mathbb{N})$.
\end{corollary}

\begin{proof}
By Observation \ref{obs:main-cospec-inv}, Lemma \ref{lem:cospectrum-preservation} and Theorem \ref{thm:cospec-univscottdom}.
\end{proof}

\section*{Acknowledgements}

The work has benefited from the Marie Curie International Research Staff Exchange Scheme \emph{Computable
Analysis}, PIRSES-GA-2011- 294962.
The first author was partially supported by a Grant-in-Aid for JSPS fellows.
The authors are also grateful to Masahiro Kumabe, Joseph Miller, Luca Motto Ros, Philipp Schlicht, and Takamitsu Yamauchi for their insightful comments and discussions.

\bibliographystyle{eptcs}
\bibliography{pointspectra}

\def\cprime{$'$}
\begin{thebibliography}{10}
\providecommand{\bibitemstart}[1]{\bibitem{#1}}
\providecommand{\bibitemend}{}
\providecommand{\bibliographystart}{}
\providecommand{\bibliographyend}{}
\providecommand{\url}[1]{\texttt{#1}}
\providecommand{\urlprefix}{Available at }
\providecommand{\bibinfo}[2]{#2}
\bibliographystart

\bibitemstart{AdGr78}
\bibinfo{author}{David~F. Addis} \& \bibinfo{author}{John~H. Gresham}
  (\bibinfo{year}{1978}): \emph{\bibinfo{title}{A class of infinite-dimensional
  spaces. {I}. {D}imension theory and {A}lexandroff's problem}}.
\newblock {\sl \bibinfo{journal}{Fund. Math.}}
  \bibinfo{volume}{101}(\bibinfo{number}{3}), pp. \bibinfo{pages}{195--205}.
\bibitemend

\bibitemstart{Alek51}
\bibinfo{author}{P.~S. Aleksandrov} (\bibinfo{year}{1951}):
  \emph{\bibinfo{title}{The present status of the theory of dimension}}.
\newblock {\sl \bibinfo{journal}{Uspehi Matem. Nauk (N.S.)}}
  \bibinfo{volume}{6}(\bibinfo{number}{4(45)}), pp. \bibinfo{pages}{43--68}.
\bibitemend

\bibitemstart{Anc85}
\bibinfo{author}{Fredric~D. Ancel} (\bibinfo{year}{1985}):
  \emph{\bibinfo{title}{The role of countable dimensionality in the theory of
  cell-like relations}}.
\newblock {\sl \bibinfo{journal}{Trans. Amer. Math. Soc.}}
  \bibinfo{volume}{287}(\bibinfo{number}{1}), pp. \bibinfo{pages}{1--40}.
\bibitemend

\bibitemstart{ArChPu00}
\bibinfo{author}{F.~G. Arenas}, \bibinfo{author}{V.~A. Chatyrko} \&
  \bibinfo{author}{M.~L. Puertas} (\bibinfo{year}{2000}):
  \emph{\bibinfo{title}{Transfinite extension of {S}teinke's dimension}}.
\newblock {\sl \bibinfo{journal}{Acta Math. Hungar.}}
  \bibinfo{volume}{88}(\bibinfo{number}{1-2}), pp. \bibinfo{pages}{105--112}.
\bibitemend

\bibitemstart{kalimullin}
\bibinfo{author}{M.~M. Arslanov}, \bibinfo{author}{I.~Sh. Kalimullin} \&
  \bibinfo{author}{S.~B. Kuper} (\bibinfo{year}{2003}):
  \emph{\bibinfo{title}{Splitting properties of total enumeration degrees}}.
\newblock {\sl \bibinfo{journal}{Algebra Logika}}
  \bibinfo{volume}{42}(\bibinfo{number}{1}), pp. \bibinfo{pages}{3--25, 125}.
\bibitemend

\bibitemstart{Bab05}
\bibinfo{author}{Liljana Babinkostova} (\bibinfo{year}{2005}):
  \emph{\bibinfo{title}{Selective screenability game and covering dimension}}.
\newblock {\sl \bibinfo{journal}{Topology Proc.}}
  \bibinfo{volume}{29}(\bibinfo{number}{1}), pp. \bibinfo{pages}{13--17}.
\bibitemend

\bibitemstart{Bade73}
\bibinfo{author}{William~G. Bade} (\bibinfo{year}{1973}):
  \emph{\bibinfo{title}{Complementation problems for the {B}aire classes}}.
\newblock {\sl \bibinfo{journal}{Pacific J. Math.}} \bibinfo{volume}{45}, pp.
  \bibinfo{pages}{1--11}.
\bibitemend

\bibitemstart{debrecht6}
\bibinfo{author}{Matthew de~Brecht} (\bibinfo{year}{2013}):
  \emph{\bibinfo{title}{Quasi-{P}olish spaces}}.
\newblock {\sl \bibinfo{journal}{Ann. Pure Appl. Logic}}
  \bibinfo{volume}{164}(\bibinfo{number}{3}), pp. \bibinfo{pages}{356--381}.
\bibitemend

\bibitemstart{CenMau84}
\bibinfo{author}{Douglas Cenzer} \& \bibinfo{author}{R.~Daniel Mauldin}
  (\bibinfo{year}{1984}): \emph{\bibinfo{title}{Borel equivalence and
  isomorphism of coanalytic sets}}.
\newblock {\sl \bibinfo{journal}{Dissertationes Math. (Rozprawy Mat.)}}
  \bibinfo{volume}{228}, p.~\bibinfo{pages}{28}.
\bibitemend

\bibitemstart{ChaHat13}
\bibinfo{author}{V.~Chatyrko} \& \bibinfo{author}{Ya. Hattori}
  (\bibinfo{year}{2013}): \emph{\bibinfo{title}{Small scattered topological
  invariants}}.
\newblock {\sl \bibinfo{journal}{Mat. Stud.}}
  \bibinfo{volume}{39}(\bibinfo{number}{2}), pp. \bibinfo{pages}{212--222}.
\bibitemend

\bibitemstart{Cha99}
\bibinfo{author}{Vitalij~A. Chatyrko} (\bibinfo{year}{1999}):
  \emph{\bibinfo{title}{On locally {$r$}-incomparable families of
  infinite-dimensional {C}antor manifolds}}.
\newblock {\sl \bibinfo{journal}{Comment. Math. Univ. Carolin.}}
  \bibinfo{volume}{40}(\bibinfo{number}{1}), pp. \bibinfo{pages}{165--173}.
\bibitemend

\bibitemstart{ChaEP99}
\bibinfo{author}{Vitalij~A. Chatyrko} \& \bibinfo{author}{El{\.z}bieta Pol}
  (\bibinfo{year}{2000}): \emph{\bibinfo{title}{Continuum many {F}r\'echet
  types of hereditarily strongly infinite-dimensional {C}antor manifolds}}.
\newblock {\sl \bibinfo{journal}{Proc. Amer. Math. Soc.}}
  \bibinfo{volume}{128}(\bibinfo{number}{4}), pp. \bibinfo{pages}{1207--1213}.
\bibitemend

\bibitemstart{ChYuBook}
\bibinfo{author}{Chi~Tat Chong} \& \bibinfo{author}{Liang~and Yu}
  (\bibinfo{year}{2015}): \emph{\bibinfo{title}{Recursion theory}}, {\sl
  \bibinfo{series}{De Gruyter Series in Logic and its
  Applications}}~\bibinfo{volume}{8}.
\newblock \bibinfo{publisher}{De Gruyter, Berlin}.
\newblock \urlprefix\url{http://dx.doi.org/10.1515/9783110275643}.
\newblock \bibinfo{note}{Computational aspects of definability, With an
  interview with Gerald E. Sacks}.
\bibitemend

\bibitemstart{Dash74}
\bibinfo{author}{F.~K. Dashiell, Jr.} (\bibinfo{year}{1974}):
  \emph{\bibinfo{title}{Isomorphism problems for the {B}aire classes}}.
\newblock {\sl \bibinfo{journal}{Pacific J. Math.}} \bibinfo{volume}{52}, pp.
  \bibinfo{pages}{29--43}.
\bibitemend

\bibitemstart{daymiller}
\bibinfo{author}{Adam~R. Day} \& \bibinfo{author}{Joseph~S. Miller}
  (\bibinfo{year}{2013}): \emph{\bibinfo{title}{Randomness for non-computable
  measures}}.
\newblock {\sl \bibinfo{journal}{Trans. Amer. Math. Soc.}}
  \bibinfo{volume}{365}(\bibinfo{number}{7}), pp. \bibinfo{pages}{3575--3591}.
\bibitemend

\bibitemstart{miller4}
\bibinfo{author}{Characterizing the~continuous degrees} (\bibinfo{year}{2017}).
\newblock \emph{\bibinfo{title}{Uri Andrews and Gregory Igusa and Joseph S.
  Miller and Mariya I. Soskova}}.
\newblock \bibinfo{howpublished}{submitted for publication}.
\newblock \urlprefix\url{http://www.math.wisc.edu/~jmiller/Papers/codable.pdf}.
\bibitemend

\bibitemstart{EngBook}
\bibinfo{author}{Ryszard Engelking} (\bibinfo{year}{1978}):
  \emph{\bibinfo{title}{Dimension Theory}}.
\newblock \bibinfo{publisher}{North-Holland Publishing Co.,
  Amsterdam-Oxford-New York; PWN---Polish Scientific Publishers, Warsaw}.
\newblock \bibinfo{note}{{North-Holland Mathematical Library}, 19}.
\bibitemend

\bibitemstart{FedOsi08}
\bibinfo{author}{V.~V. Fedorchuk} \& \bibinfo{author}{E.~V. Osipov}
  (\bibinfo{year}{2008}): \emph{\bibinfo{title}{Certain classes of weakly
  infinite-dimensional spaces and topological games}}.
\newblock {\sl \bibinfo{journal}{Topology Appl.}}
  \bibinfo{volume}{156}(\bibinfo{number}{1}), pp. \bibinfo{pages}{61--69}.
\bibitemend

\bibitemstart{friedberg}
\bibinfo{author}{Richard~M. Friedberg} \& \bibinfo{author}{Hartley Rogers, Jr.}
  (\bibinfo{year}{1959}): \emph{\bibinfo{title}{Reducibility and completeness
  for sets of integers}}.
\newblock {\sl \bibinfo{journal}{Z. Math. Logik Grundlagen Math.}}
  \bibinfo{volume}{5}, pp. \bibinfo{pages}{117--125}.
\bibitemend

\bibitemstart{GanSos15}
\bibinfo{author}{Hristo~A. Ganchev} \& \bibinfo{author}{Mariya~I. Soskova}
  (\bibinfo{year}{2015}): \emph{\bibinfo{title}{Definability via {K}alimullin
  pairs in the structure of the enumeration degrees}}.
\newblock {\sl \bibinfo{journal}{Trans. Amer. Math. Soc.}}
  \bibinfo{volume}{367}(\bibinfo{number}{7}), pp. \bibinfo{pages}{4873--4893}.
\bibitemend

\bibitemstart{gregoriades3}
\bibinfo{author}{Vassilios Gregoriades} (\bibinfo{year}{2016}):
  \emph{\bibinfo{title}{Classes of {P}olish spaces under effective {B}orel
  isomorphism}}.
\newblock {\sl \bibinfo{journal}{Memoirs of the American Mathematical Society}}
  .
\bibitemend

\bibitemstart{GreKih}
\bibinfo{author}{Vassilios Gregoriades}, \bibinfo{author}{Takayuki Kihara} \&
  \bibinfo{author}{Keng~Meng Ng} (\bibinfo{year}{2014}).
\newblock \emph{\bibinfo{title}{Turing degrees in Polish spaces and
  decomposability of Borel functions}}.
\newblock \bibinfo{howpublished}{arXiv 1410.1052}.
\bibitemend

\bibitemstart{pauly-gregoriades-arxiv}
\bibinfo{author}{Vassilios Gregoriades}, \bibinfo{author}{Tam\'as Kisp\'eter}
  \& \bibinfo{author}{Arno Pauly} (\bibinfo{year}{2016}):
  \emph{\bibinfo{title}{A comparison of concepts from computable analysis and
  effective descriptive set theory}}.
\newblock {\sl \bibinfo{journal}{Mathematical Structures in Computer Science}}
  \urlprefix\url{http://arxiv.org/abs/1403.7997}.
\bibitemend

\bibitemstart{grubba3}
\bibinfo{author}{Tanja Grubba}, \bibinfo{author}{Matthias Schr{\"o}der} \&
  \bibinfo{author}{Klaus Weihrauch} (\bibinfo{year}{2007}):
  \emph{\bibinfo{title}{Computable metrization}}.
\newblock {\sl \bibinfo{journal}{MLQ Math. Log. Q.}}
  \bibinfo{volume}{53}(\bibinfo{number}{4-5}), pp. \bibinfo{pages}{381--395}.
\bibitemend

\bibitemstart{Har78}
\bibinfo{author}{Leo Harrington} (\bibinfo{year}{1978}):
  \emph{\bibinfo{title}{Analytic determinacy and {$0^{\sharp }$}}}.
\newblock {\sl \bibinfo{journal}{J. Symbolic Logic}}
  \bibinfo{volume}{43}(\bibinfo{number}{4}), pp. \bibinfo{pages}{685--693}.
\bibitemend

\bibitemstart{Hav74}
\bibinfo{author}{William~E. Haver} (\bibinfo{year}{1974}):
  \emph{\bibinfo{title}{A covering property for metric spaces}}.
\newblock In: {\sl \bibinfo{booktitle}{Topology {C}onference ({V}irginia
  {P}olytech. {I}nst. and {S}tate {U}niv., {B}lacksburg, {V}a., 1973)}},
  \bibinfo{publisher}{Springer, Berlin}, pp. \bibinfo{pages}{108--113. Lecture
  Notes in Math., Vol. 375}.
\bibitemend

\bibitemstart{kihara3}
\bibinfo{author}{K.~Higuchi} \& \bibinfo{author}{T.~Kihara}
  (\bibinfo{year}{2014}): \emph{\bibinfo{title}{Inside the {M}uchnik degrees
  {I}: {D}iscontinuity, learnability and constructivism}}.
\newblock {\sl \bibinfo{journal}{Ann. Pure Appl. Logic}}
  \bibinfo{volume}{165}(\bibinfo{number}{5}), pp. \bibinfo{pages}{1058--1114}.
\bibitemend

\bibitemstart{kihara3b}
\bibinfo{author}{K.~Higuchi} \& \bibinfo{author}{T.~Kihara}
  (\bibinfo{year}{2014}): \emph{\bibinfo{title}{Inside the {M}uchnik degrees
  {II}: {T}he degree structures induced by the arithmetical hierarchy of
  countably continuous functions}}.
\newblock {\sl \bibinfo{journal}{Ann. Pure Appl. Logic}}
  \bibinfo{volume}{165}(\bibinfo{number}{6}), pp. \bibinfo{pages}{1201--1241}.
\bibitemend

\bibitemstart{Hinman73}
\bibinfo{author}{Peter~G. Hinman} (\bibinfo{year}{1973}):
  \emph{\bibinfo{title}{Degrees of continuous functionals}}.
\newblock {\sl \bibinfo{journal}{J. Symbolic Logic}} \bibinfo{volume}{38}, pp.
  \bibinfo{pages}{393--395}.
\bibitemend

\bibitemstart{HKSS02}
\bibinfo{author}{Denis~R. Hirschfeldt}, \bibinfo{author}{Bakhadyr Khoussainov},
  \bibinfo{author}{Richard~A. Shore} \& \bibinfo{author}{Arkadii~M. Slinko}
  (\bibinfo{year}{2002}): \emph{\bibinfo{title}{Degree spectra and computable
  dimensions in algebraic structures}}.
\newblock {\sl \bibinfo{journal}{Ann. Pure Appl. Logic}}
  \bibinfo{volume}{115}(\bibinfo{number}{1-3}), pp. \bibinfo{pages}{71--113}.
\bibitemend

\bibitemstart{Hrba78}
\bibinfo{author}{Karel Hrb{\'a}{\v{c}}ek} (\bibinfo{year}{1978}):
  \emph{\bibinfo{title}{On the complexity of analytic sets}}.
\newblock {\sl \bibinfo{journal}{Z. Math. Logik Grundlag. Math.}}
  \bibinfo{volume}{24}(\bibinfo{number}{5}), pp. \bibinfo{pages}{419--425}.
\bibitemend

\bibitemstart{hurewicz}
\bibinfo{author}{Witold Hurewicz} \& \bibinfo{author}{Henry Wallman}
  (\bibinfo{year}{1941}): \emph{\bibinfo{title}{Dimension {T}heory}}.
\newblock \bibinfo{series}{Princeton Mathematical Series, v. 4}.
  \bibinfo{publisher}{Princeton University Press, Princeton, N. J.}
\bibitemend

\bibitemstart{Jayne74}
\bibinfo{author}{J.~E. Jayne} (\bibinfo{year}{1974}): \emph{\bibinfo{title}{The
  space of class {$\alpha $} {B}aire functions}}.
\newblock {\sl \bibinfo{journal}{Bull. Amer. Math. Soc.}} \bibinfo{volume}{80},
  pp. \bibinfo{pages}{1151--1156}.
\bibitemend

\bibitemstart{JayRog79a}
\bibinfo{author}{J.~E. Jayne} \& \bibinfo{author}{C.~A. Rogers}
  (\bibinfo{year}{1979}): \emph{\bibinfo{title}{Borel isomorphisms at the first
  level. {I}}}.
\newblock {\sl \bibinfo{journal}{Mathematika}}
  \bibinfo{volume}{26}(\bibinfo{number}{1}), pp. \bibinfo{pages}{125--156}.
\bibitemend

\bibitemstart{JayRog79b}
\bibinfo{author}{J.~E. Jayne} \& \bibinfo{author}{C.~A. Rogers}
  (\bibinfo{year}{1979}): \emph{\bibinfo{title}{Borel isomorphisms at the first
  level. {II}}}.
\newblock {\sl \bibinfo{journal}{Mathematika}}
  \bibinfo{volume}{26}(\bibinfo{number}{2}), pp. \bibinfo{pages}{157--179}.
\bibitemend

\bibitemstart{jockusch2}
\bibinfo{author}{Carl~G. Jockusch, Jr.} (\bibinfo{year}{1968}):
  \emph{\bibinfo{title}{Semirecursive sets and positive reducibility}}.
\newblock {\sl \bibinfo{journal}{Trans. Amer. Math. Soc.}}
  \bibinfo{volume}{131}, pp. \bibinfo{pages}{420--436}.
\bibitemend

\bibitemstart{KiNg}
\bibinfo{author}{T.~Kihara} \& \bibinfo{author}{K.~M. Ng}:
  \emph{\bibinfo{title}{A Generalization of the {Shore-Slaman} join theorem in
  {P}olish spaces}}.
\newblock \bibinfo{note}{In preparation}.
\bibitemend

\bibitemstart{edegrees}
\bibinfo{author}{Takayuki Kihara}, \bibinfo{author}{Steffen Lempp},
  \bibinfo{author}{Keng~Meng Ng} \& \bibinfo{author}{Arno Pauly}
  (\bibinfo{year}{201X}).
\newblock \emph{\bibinfo{title}{Enumeration degrees and non-metrizable
  topology}}.
\newblock \bibinfo{howpublished}{in preparation}.
\bibitemend

\bibitemstart{kulesza}
\bibinfo{author}{John Kulesza} (\bibinfo{year}{1990}):
  \emph{\bibinfo{title}{The dimension of products of complete separable metric
  spaces}}.
\newblock {\sl \bibinfo{journal}{Fundamenta Mathematicae}}
  \bibinfo{volume}{135}(\bibinfo{number}{1}), pp. \bibinfo{pages}{49--54}.
\bibitemend

\bibitemstart{lelek}
\bibinfo{author}{A.~Lelek} (\bibinfo{year}{1965}): \emph{\bibinfo{title}{On the
  dimensionality of remainders in compactifications}}.
\newblock {\sl \bibinfo{journal}{Dokl. Akad. Nauk SSSR}} \bibinfo{volume}{160},
  pp. \bibinfo{pages}{534--537}.
\bibitemend

\bibitemstart{lutz2}
\bibinfo{author}{Jack~H. Lutz} (\bibinfo{year}{2003}):
  \emph{\bibinfo{title}{The dimensions of individual strings and sequences}}.
\newblock {\sl \bibinfo{journal}{Inform. and Comput.}}
  \bibinfo{volume}{187}(\bibinfo{number}{1}), pp. \bibinfo{pages}{49--79}.
\bibitemend

\bibitemstart{Maul76}
\bibinfo{author}{R.~Daniel Mauldin} (\bibinfo{year}{1976}):
  \emph{\bibinfo{title}{On nonisomorphic analytic sets}}.
\newblock {\sl \bibinfo{journal}{Proc. Amer. Math. Soc.}} \bibinfo{volume}{58},
  pp. \bibinfo{pages}{241--244}.
\bibitemend

\bibitemstart{medvedev}
\bibinfo{author}{Yu.~T. Medvedev} (\bibinfo{year}{1955}):
  \emph{\bibinfo{title}{Degrees of difficulty of the mass problem}}.
\newblock {\sl \bibinfo{journal}{Dokl. Akad. Nauk SSSR (N.S.)}}
  \bibinfo{volume}{104}, pp. \bibinfo{pages}{501--504}.
\bibitemend

\bibitemstart{vMBook2}
\bibinfo{author}{Jan van Mill} (\bibinfo{year}{2001}):
  \emph{\bibinfo{title}{The Infinite-Dimensional Topology of Function Spaces}},
  {\sl \bibinfo{series}{North-Holland Mathematical
  Library}}~\bibinfo{volume}{64}.
\newblock \bibinfo{publisher}{North-Holland Publishing Co., Amsterdam}.
\bibitemend

\bibitemstart{miller2}
\bibinfo{author}{Joseph~S. Miller} (\bibinfo{year}{2004}):
  \emph{\bibinfo{title}{Degrees of unsolvability of continuous functions}}.
\newblock {\sl \bibinfo{journal}{J. Symbolic Logic}}
  \bibinfo{volume}{69}(\bibinfo{number}{2}), pp. \bibinfo{pages}{555--584}.
\bibitemend

\bibitemstart{Monta13}
\bibinfo{author}{Antonio Montalb{\'a}n} (\bibinfo{year}{2013}):
  \emph{\bibinfo{title}{A computability theoretic equivalent to {V}aught's
  conjecture}}.
\newblock {\sl \bibinfo{journal}{Adv. Math.}} \bibinfo{volume}{235}, pp.
  \bibinfo{pages}{56--73}.
\bibitemend

\bibitemstart{Mont14}
\bibinfo{author}{Antonio Montalb\'an} (\bibinfo{year}{2014}):
  \emph{\bibinfo{title}{Computability theoretic classifications for classes of
  structures}}.
\newblock {\sl \bibinfo{journal}{Proceedings of ICM 2014}} , pp.
  \bibinfo{pages}{79--101}.
\bibitemend

\bibitemstart{MRos13}
\bibinfo{author}{Luca Motto~Ros} (\bibinfo{year}{2013}):
  \emph{\bibinfo{title}{On the structure of finite level and
  {$\omega$}-decomposable {B}orel functions}}.
\newblock {\sl \bibinfo{journal}{J. Symbolic Logic}}
  \bibinfo{volume}{78}(\bibinfo{number}{4}), pp. \bibinfo{pages}{1257--1287}.
\bibitemend

\bibitemstart{schlicht}
\bibinfo{author}{Luca Motto~Ros}, \bibinfo{author}{Philipp Schlicht} \&
  \bibinfo{author}{Victor Selivanov} (\bibinfo{year}{2014}):
  \emph{\bibinfo{title}{Wadge-like reducibilities on arbitrary quasi-Polish
  spaces}}.
\newblock {\sl \bibinfo{journal}{Mathematical Structures in Computer Science}}
  , pp.
  \bibinfo{pages}{1--50}\urlprefix\url{http://journals.cambridge.org/article_S0960129513000339}.
\newblock \bibinfo{note}{ArXiv 1204.5338}.
\bibitemend

\bibitemstart{muchnik}
\bibinfo{author}{A.~A. Mu{\v{c}}nik} (\bibinfo{year}{1963}):
  \emph{\bibinfo{title}{On strong and weak reducibility of algorithmic
  problems}}.
\newblock {\sl \bibinfo{journal}{Sibirsk. Mat. \v Z.}} \bibinfo{volume}{4}, pp.
  \bibinfo{pages}{1328--1341}.
\bibitemend

\bibitemstart{nies}
\bibinfo{author}{Andr{\'e} Nies} (\bibinfo{year}{2009}):
  \emph{\bibinfo{title}{Computability and Randomness}}, {\sl
  \bibinfo{series}{Oxford Logic Guides}}~\bibinfo{volume}{51}.
\newblock \bibinfo{publisher}{Oxford University Press, Oxford}.
\bibitemend

\bibitemstart{OdiBook1}
\bibinfo{author}{P.~G. Odifreddi} (\bibinfo{year}{1999}):
  \emph{\bibinfo{title}{Classical Recursion Theory. {V}ol. {II}}}, {\sl
  \bibinfo{series}{Studies in Logic and the Foundations of Mathematics}}
  \bibinfo{volume}{143}.
\newblock \bibinfo{publisher}{North-Holland Publishing Co., Amsterdam}.
\bibitemend

\bibitemstart{OdiBook}
\bibinfo{author}{Piergiorgio Odifreddi} (\bibinfo{year}{1989}):
  \emph{\bibinfo{title}{Classical Recursion Theory}}, {\sl
  \bibinfo{series}{Studies in Logic and the Foundations of Mathematics}}
  \bibinfo{volume}{125}.
\newblock \bibinfo{publisher}{North-Holland Publishing Co., Amsterdam}.
\bibitemend

\bibitemstart{pauly-overview-arxiv}
\bibinfo{author}{Arno Pauly} (\bibinfo{year}{2014}).
\newblock \emph{\bibinfo{title}{The descriptive theory of represented spaces}}.
\newblock \bibinfo{howpublished}{arXiv:1408.5329}.
\bibitemend

\bibitemstart{pauly-synthetic-arxiv}
\bibinfo{author}{Arno Pauly} (\bibinfo{year}{2016}): \emph{\bibinfo{title}{On
  the topological aspects of the theory of represented spaces}}.
\newblock {\sl \bibinfo{journal}{Computability}}
  \bibinfo{volume}{5}(\bibinfo{number}{2}), pp. \bibinfo{pages}{159--180}.
\newblock \urlprefix\url{http://arxiv.org/abs/1204.3763}.
\bibitemend

\bibitemstart{pauly-descriptive}
\bibinfo{author}{Arno Pauly} \& \bibinfo{author}{Matthew de~Brecht}.
\newblock \emph{\bibinfo{title}{Towards Synthetic Descriptive Set Theory: An
  instantiation with represented spaces}}.
\newblock \bibinfo{howpublished}{arXiv 1307.1850}.
\bibitemend

\bibitemstart{paulydebrecht}
\bibinfo{author}{Arno Pauly} \& \bibinfo{author}{Matthew de~Brecht}
  (\bibinfo{year}{2014}): \emph{\bibinfo{title}{Non-deterministic Computation
  and the {J}ayne {R}ogers Theorem}}.
\newblock {\sl \bibinfo{journal}{Electronic Proceedings in Theoretical Computer
  Science}} \bibinfo{volume}{143}.
\newblock \bibinfo{note}{{DCM 2012}}.
\bibitemend

\bibitemstart{PaSa12}
\bibinfo{author}{Janusz Pawlikowski} \& \bibinfo{author}{Marcin Sabok}
  (\bibinfo{year}{2012}): \emph{\bibinfo{title}{Decomposing {B}orel functions
  and structure at finite levels of the {B}aire hierarchy}}.
\newblock {\sl \bibinfo{journal}{Ann. Pure Appl. Logic}}
  \bibinfo{volume}{163}(\bibinfo{number}{12}), pp. \bibinfo{pages}{1748--1764}.
\bibitemend

\bibitemstart{EPol96}
\bibinfo{author}{El{\.z}bieta Pol} (\bibinfo{year}{1996}):
  \emph{\bibinfo{title}{On infinite-dimensional {C}antor manifolds}}.
\newblock {\sl \bibinfo{journal}{Topology Appl.}}
  \bibinfo{volume}{71}(\bibinfo{number}{3}), pp. \bibinfo{pages}{265--276}.
\bibitemend

\bibitemstart{PoZa12}
\bibinfo{author}{R.~Pol} \& \bibinfo{author}{P.~Zakrzewski}
  (\bibinfo{year}{2012}): \emph{\bibinfo{title}{On {B}orel mappings and
  {$\sigma$}-ideals generated by closed sets}}.
\newblock {\sl \bibinfo{journal}{Adv. Math.}}
  \bibinfo{volume}{231}(\bibinfo{number}{2}), pp. \bibinfo{pages}{651--663}.
\bibitemend

\bibitemstart{pol2}
\bibinfo{author}{Roman Pol} (\bibinfo{year}{1981}): \emph{\bibinfo{title}{A
  weakly infinite-dimensional compactum which is not countable-dimensional}}.
\newblock {\sl \bibinfo{journal}{Proc. Amer. Math. Soc.}}
  \bibinfo{volume}{82}(\bibinfo{number}{4}), pp. \bibinfo{pages}{634--636}.
\bibitemend

\bibitemstart{pourel}
\bibinfo{author}{Marian~B. Pour-El} \& \bibinfo{author}{J.~Ian Richards}
  (\bibinfo{year}{1989}): \emph{\bibinfo{title}{Computability in Analysis and
  Physics}}.
\newblock \bibinfo{series}{Perspectives in Mathematical Logic}.
  \bibinfo{publisher}{Springer-Verlag, Berlin}.
\bibitemend

\bibitemstart{Radul}
\bibinfo{author}{T.~M. Radul} (\bibinfo{year}{2006}): \emph{\bibinfo{title}{On
  classification of sigma hereditary disconnected spaces}}.
\newblock {\sl \bibinfo{journal}{Mat. Stud.}}
  \bibinfo{volume}{26}(\bibinfo{number}{1}), pp. \bibinfo{pages}{97--100}.
\bibitemend

\bibitemstart{Rich81}
\bibinfo{author}{Linda~Jean Richter} (\bibinfo{year}{1981}):
  \emph{\bibinfo{title}{Degrees of structures}}.
\newblock {\sl \bibinfo{journal}{J. Symbolic Logic}}
  \bibinfo{volume}{46}(\bibinfo{number}{4}), pp. \bibinfo{pages}{723--731}.
\bibitemend

\bibitemstart{Rub80}
\bibinfo{author}{Leonard~R. Rubin} (\bibinfo{year}{1980}):
  \emph{\bibinfo{title}{Noncompact hereditarily strongly infinite-dimensional
  spaces}}.
\newblock {\sl \bibinfo{journal}{Proc. Amer. Math. Soc.}}
  \bibinfo{volume}{79}(\bibinfo{number}{1}), pp. \bibinfo{pages}{153--154}.
\bibitemend

\bibitemstart{schori}
\bibinfo{author}{Leonard~R. Rubin}, \bibinfo{author}{R.~M. Schori} \&
  \bibinfo{author}{John~J. Walsh} (\bibinfo{year}{1979}):
  \emph{\bibinfo{title}{New dimension-theory techniques for constructing
  infinite-dimensional examples}}.
\newblock {\sl \bibinfo{journal}{General Topology Appl.}}
  \bibinfo{volume}{10}(\bibinfo{number}{1}), pp. \bibinfo{pages}{93--102}.
\bibitemend

\bibitemstart{SacksBook}
\bibinfo{author}{Gerald~E. Sacks} (\bibinfo{year}{1990}):
  \emph{\bibinfo{title}{Higher recursion theory}}.
\newblock \bibinfo{series}{Perspectives in Mathematical Logic}.
  \bibinfo{publisher}{Springer-Verlag, Berlin}.
\newblock \urlprefix\url{http://dx.doi.org/10.1007/BFb0086109}.
\bibitemend

\bibitemstart{schroder5}
\bibinfo{author}{Matthias Schr\"oder} (\bibinfo{year}{2002}):
  \emph{\bibinfo{title}{Admissible Representations for Continuous
  Computations}}.
\newblock \bibinfo{type}{Ph.D. thesis}, \bibinfo{school}{FernUniversit\"at
  Hagen}.
\bibitemend

\bibitemstart{schroder}
\bibinfo{author}{Matthias Schr{\"o}der} (\bibinfo{year}{2002}):
  \emph{\bibinfo{title}{Extended admissibility}}.
\newblock {\sl \bibinfo{journal}{Theoret. Comput. Sci.}}
  \bibinfo{volume}{284}(\bibinfo{number}{2}), pp. \bibinfo{pages}{519--538}.
\bibitemend

\bibitemstart{schroder6}
\bibinfo{author}{Matthias Schr{\"o}der} (\bibinfo{year}{2004}):
  \emph{\bibinfo{title}{Spaces allowing type-2 complexity theory revisited}}.
\newblock {\sl \bibinfo{journal}{MLQ Math. Log. Q.}}
  \bibinfo{volume}{50}(\bibinfo{number}{4-5}), pp. \bibinfo{pages}{443--459}.
\bibitemend

\bibitemstart{selman}
\bibinfo{author}{Alan~L. Selman} (\bibinfo{year}{1971}):
  \emph{\bibinfo{title}{Arithmetical reducibilities. {I}}}.
\newblock {\sl \bibinfo{journal}{Z. Math. Logik Grundlagen Math.}}
  \bibinfo{volume}{17}, pp. \bibinfo{pages}{335--350}.
\bibitemend

\bibitemstart{SoareBook}
\bibinfo{author}{Robert~I. Soare} (\bibinfo{year}{1987}):
  \emph{\bibinfo{title}{Recursively Enumerable Sets and Degrees}}.
\newblock \bibinfo{series}{Perspectives in Mathematical Logic}.
  \bibinfo{publisher}{Springer-Verlag, Berlin}.
\bibitemend

\bibitemstart{CTopBook}
\bibinfo{author}{Lynn~Arthur Steen} \& \bibinfo{author}{J.~Arthur Seebach, Jr.}
  (\bibinfo{year}{1978}): \emph{\bibinfo{title}{Counterexamples in topology}}.
\newblock \bibinfo{publisher}{Springer-Verlag, New York-Heidelberg},
  \bibinfo{edition}{second} edition.
\bibitemend

\bibitemstart{Stuk07}
\bibinfo{author}{A.~I. Stukachev} (\bibinfo{year}{2007}):
  \emph{\bibinfo{title}{Degrees of presentability of structures. {I}}}.
\newblock {\sl \bibinfo{journal}{Algebra and Logic}}
  \bibinfo{volume}{46}(\bibinfo{number}{6}), pp. \bibinfo{pages}{419--432}.
\bibitemend

\bibitemstart{taylor}
\bibinfo{author}{Paul Taylor} (\bibinfo{year}{2010}): \emph{\bibinfo{title}{A
  lambda calculus for real analysis}}.
\newblock {\sl \bibinfo{journal}{J. Log. Anal.}} \bibinfo{volume}{2}, pp.
  \bibinfo{pages}{1--115}.
\bibitemend

\bibitemstart{taylor2}
\bibinfo{author}{Paul Taylor} (\bibinfo{year}{2011}):
  \emph{\bibinfo{title}{Foundations {\it for} computable topology}}.
\newblock In: {\sl \bibinfo{booktitle}{Foundational theories of classical and
  constructive mathematics}}, {\sl \bibinfo{series}{West. Ont. Ser. Philos.
  Sci.}}~\bibinfo{volume}{76}, \bibinfo{publisher}{Springer, Dordrecht}, pp.
  \bibinfo{pages}{265--310}.
\bibitemend

\bibitemstart{weihrauchd}
\bibinfo{author}{Klaus Weihrauch} (\bibinfo{year}{2000}):
  \emph{\bibinfo{title}{Computable Analysis: An Introduction}}.
\newblock \bibinfo{series}{Texts in Theoretical Computer Science. An EATCS
  Series}. \bibinfo{publisher}{Springer-Verlag, Berlin}.
\bibitemend

\bibitemstart{weihrauchm}
\bibinfo{author}{Klaus Weihrauch} (\bibinfo{year}{2013}):
  \emph{\bibinfo{title}{Computably regular topological spaces}}.
\newblock {\sl \bibinfo{journal}{Log. Methods Comput. Sci.}}
  \bibinfo{volume}{9}(\bibinfo{number}{3}), pp. \bibinfo{pages}{3:5, 24}.
\bibitemend

\bibitemstart{weihrauchf}
\bibinfo{author}{Klaus Weirauch} (\bibinfo{year}{2003}):
  \emph{\bibinfo{title}{Computational complexity on computable metric spaces}}.
\newblock {\sl \bibinfo{journal}{MLQ Math. Log. Q.}}
  \bibinfo{volume}{49}(\bibinfo{number}{1}), pp. \bibinfo{pages}{3--21}.
\bibitemend

\bibitemstart{Zap14}
\bibinfo{author}{Jind{\v{r}}ich Zapletal} (\bibinfo{year}{2014}):
  \emph{\bibinfo{title}{Dimension theory and forcing}}.
\newblock {\sl \bibinfo{journal}{Topology Appl.}} \bibinfo{volume}{167}, pp.
  \bibinfo{pages}{31--35}.
\bibitemend

\bibliographyend
\end{thebibliography}
\end{document}